\documentclass[10pt]{article}
\usepackage{amsmath,amssymb,amsthm,mathtools,microtype}
\usepackage{slashed}
\usepackage{fullpage}
\usepackage{enumitem}
\usepackage{needspace}
\usepackage{xcolor}
\usepackage{mathrsfs}
\usepackage[hidelinks]{hyperref}
\setlist[enumerate]{label=\textup{(\roman*)},nosep,leftmargin=*}
\newtheorem{theorem}{Theorem}[section]

\newtheorem{lemma}[theorem]{Lemma}

\newtheorem{corollary}[theorem]{Corollary}
\newtheorem{conjecture}[theorem]{Conjecture}
\newtheorem{question}[theorem]{Question}
\newtheorem{claim}{Claim}
\theoremstyle{definition}
\newtheorem{definition}[theorem]{Definition}
\newtheorem{remark}[theorem]{Remark}

\newcommand{\cG}{\mathcal G}
\newcommand{\cK}{\mathcal K}

\newcommand{\cF}{\mathcal F}
\newcommand{\cA}{\mathcal A}
\newcommand{\cT}{\mathcal T}
\newcommand{\cH}{\mathcal H}
\newcommand{\cQ}{\mathcal Q}
\newcommand{\cP}{\mathcal P}

\newcommand{\pu}{\color{violet}}
\title{Ore-degree conditions for matchings in uniform hypergraphs
	\thanks{The work is partially supported by NNSF of China (Grant No. 12271251) and Fundamental Research Funds for the Central Universities (Grant No. NJ2026048).\newline
		Email addresses: weilihua@nuaa.edu.cn (L. Wei), kexxu1221@126.com (K. Xu).
		$\dagger$ Corresponding author.}}
\author{Lihua Wei $^{a,b}$, Kexiang Xu $^{a,b\dagger}$
	\\
	\small $^{a}$ School of Mathematics, Nanjing University of
	Aeronautics and Astronautics,\\
	\small Nanjing, Jiangsu 210016, PR China\\
	\small $^{b}$ MIIT Key Laboratory of Mathematical Modelling and High Performance\\
	\small Computing of Air Vehicles, Nanjing, Jiangsu 210016, China\\
}
\date{}

\begin{document}
	\maketitle
	\begin{abstract}
		
		Let $\sigma_r(\cH)$ be the minimum sum of vertex degrees over all
		nonedges of an $r$-uniform hypergraph $\cH$ on $n$ vertices.
		Balogh, Palmer and Raeisi conjectured that $n>2rs$ and
		$\sigma_r(\cH)>r[\binom{n-1}{r-1}-\binom{n-s}{r-1}]$ imply
		the existence of a matching of size $s$ in $\cH$.
		We prove that, for every fixed $r\ge3$ and all sufficiently large $s$, $\sigma_r(\cH)>r[\binom{n-1}{r-1}-\binom{n-s}{r-1}]$ forces an $s$-matching whenever $n\ge rs+\lfloor(r-1)^{r-1}/r^{r-2}\rfloor$.
		A complete $r$-graph on $rs -1$ vertices
		together with $\lfloor(r-1)^{r-1}/r^{r-2}\rfloor$ isolated vertices shows that the vertex bound is
		sharp for all sufficiently large $s$.
	 By adding new vertices to $\cH$, retaining all its edges,
		and making every $r$-set containing at least one new vertex as a new edge,
		 we also deduce that $n\ge rs+C$ suffices as a bound of the above conjecture for every $s\ge1$,
        where $C$ is a finite constant depending only on $r$.
		For $r=3$, we provide the explicit sufficient order
		$n\ge\lceil(7s+11)/2\rceil$ for $s\ge6$,  together with slightly better bounds
		for $2\le s\le5$, establishing the conjecture in this uniformity.
		For $s=2$, we improve the known sufficient order $2r+2$ to
        $2r+1$ for every $r\ge3$ and show that this order is best possible.
	\end{abstract}
	
	\section{Introduction}\label{sec-intro}
	An {\it $r$-uniform hypergraph} (\emph{$r$-graph}) $\cH$ is a
	family of $r$-subsets (the \emph{edges} of $\cH$) of a \emph{vertex} set $V$.
	Unless otherwise stated, an edge means a hyperedge throughout this paper.
	For a given hypergraph $\cH$, denote the vertex set by $V(\cH)$ and, for simplicity, we usually suppose that $V(\cH)$ is $[n]:=\{1,2,\dots, n\}$.  We identify $\cH$ with its edge set, writing $F\in\cH$
	for an edge of $\cH$ and $v\in F$ for a vertex of $F$.
	When useful, we also denote the vertex set of $F$ by $V(F)$. Thus $\lvert \cH\rvert $ is the number of edges in $\cH$.
	A hypergraph $\cG$ is a \emph{subhypergraph} of $\cH$ if
	$V(\cG)\subseteq V(\cH)$ and $\cG\subseteq \cH$.
	Unless otherwise specified, a subfamily $\cF\subseteq\cH$
	is regarded as the subhypergraph with vertex set $V(\cH)$
	and edge set $\cF$.
	For $Y\subseteq V(\cH)$, let $
	\cH[Y]:=\cH\cap\binom{Y}{r}$
	be the induced subhypergraph on vertex set $Y$.
	For $Z\subseteq V(\cH)$, write
	$\cH-Z:=\cH[V(\cH)\setminus Z]$ for the hypergraph obtained by
	removing the vertices in $Z$ and all edges incident with them.
	In particular, for $v\in V(\cH)$, we write $\cH-v:=\cH-\{v\}$.
	When a singleton vertex set is used in a set expression, we omit its braces;
	thus $X\cup v:=X\cup\{v\}$ and $X\setminus v:=X\setminus\{v\}$.
	We use the same convention for singleton sets of indices.
	For an $r$-graph $\cF$ on $V(\cH)$, the hypergraph $\cH-\cF$
	has vertex set $V(\cH)$ and edge set $\cH\setminus\cF$.
	The degree $d_{\cH}(v)$ of
	$v\in V(\cH)$ is the number of edges containing $v${\pu ;} we write $d(v)$
	when $\cH$ is clear. For an $r$-graph $\cH$, write $\delta(\cH)$ and
	$\Delta(\cH)$ for its minimum and maximum vertex degrees. For $1\le d\le r-1$, the \emph{minimum $d$-degree} is $
	\delta_d(\cH):=\min_{X\in\binom{V(\cH)}d}
	\lvert \{F\in\cH\mid X\subseteq F\}\rvert $.
	In particular, $\delta_1(\cH)=\delta(\cH)$.
	More generally, for any
	hypergraph $\cH$, a \emph{matching} is a subfamily of pairwise
	disjoint members of $\cH$, and $\nu(\cH)$ denotes its
	maximum size. An edge $F$ \emph{covers} a vertex $v$ if
	$v\in F$. For a matching $M$, write
	$V(M):=\bigcup_{F\in M}V(F)$ for the set of vertices it covers.
	A matching $M$ in $\cH$ is \emph{perfect} if it covers all vertices of $\cH$, that is, $V(M)=V(\cH)$. For an $r$-graph on $n$ vertices, this is equivalent to $\lvert M\rvert =n/r$.
	For the complete $r$-graph $\binom Vr$ on a finite vertex set $V$,
the formula in~\cite{GZM2017} gives
	\begin{equation}\label{eq-complete-match}
		\nu\left(\binom Vr\right)=\left\lfloor\frac{\lvert V\rvert}{r}\right\rfloor.
	\end{equation}
	In particular, the complete $r$-graph on $V$ has a perfect matching
	if and only if $r\mid\lvert V\rvert$.
	
	A fractional matching assigns weights $f(F)\in[0,1]$ to edges $F\in \cH$
	so that $\sum_{F\ni v}f(F)\le1$ for every vertex $v$ in the hypergraph $\cH$. A fractional vertex cover
	assigns weights $w(v)\in[0,1]$ so that $\sum_{v\in F}w(v)\ge1$
	for every edge $F$. By linear-programming duality, the maximum
	total weight of a fractional matching equals the minimum total
	weight of a fractional vertex cover; see~\cite{Schrijver2003}.
	A fractional matching is perfect
	if its edge weights sum to one at every vertex.
A family is \emph{intersecting} if every two distinct members
    have nonempty intersection.
	Families $\cF$ and $\cG$ are \emph{cross-intersecting} if
	$F\cap G\ne\emptyset$ for all $F\in\cF$ and $G\in\cG$.
	For a graph $G$, write $E(G)$ for its edge set and $d_G(v)$ for
	the degree of a vertex $v$. A nonempty graph is a \emph{star} if
	all its edges are incident with a common vertex, called its \emph{center}.
	For $x\in[n]$, the $r$-graph
	$\{F\in\binom{[n]}r\mid x\in F\}$ is called a \emph{$1$-star}~\cite{BPR2026}
	with center $x$. The \emph{Kneser graph} $KG(n,r)$ is the graph with vertex set
	$\binom{[n]}r$ in which two distinct vertices are adjacent if and only
	if the corresponding $r$-sets are disjoint.  For $n\ge2r$, the neighbors of $S\in\binom{[n]}r$ are precisely
	the $r$-subsets of $[n]\setminus S$. Thus every vertex of $KG(n,r)$
	has degree $\binom{n-r}{r}$, so $KG(n,r)$ is
	$\binom{n-r}{r}$-regular.
	
	For an $r$-set $S\subseteq V(\cH)$ of the hypergraph $\cH$, its degree sum is
	$\sum_{v\in S}d_{\cH}(v)$.  The \emph{Ore-degree} of $\cH$ is
	\[
	\sigma_r(\cH):=
	\min_{S\in\binom{V(\cH)}r\setminus\cH}\sum_{v\in S}d_{\cH}(v).
	\]
	The complete hypergraph satisfies all matching conclusions considered below by~\eqref{eq-complete-match}, so $\sigma_r(\cH)$ will only be used when $\cH$ has a nonedge.
	For integers $r>0$ and $N\ge0$, we write $r\mid N$ if
	$N=rt$ for some integer $t\ge0$.
	For $n\ge r$  and $1\leq s\leq n$, define $B_r^s(n)=\binom{n-s}{r-1}$. In particular, write $B_r^1(n)=B_r(n)$.
	For $0\le k\le n-1$, define
	\[
	D_{k,r}(n):=B_r(n)-B_r^{k+1}(n)=\binom{n-1}{r-1}-\binom{n-k-1}{r-1},
	\]
	where $
	\binom{a}{b}=0$ whenever $0\le a<b$. We write $D_k(n)=D_{k,r}(n)$
	when the uniformity $r$ is clear from context.
	
	We write $\mathbb P(\mathcal E)$ for the probability of an event
	$\mathcal E$ and $\mathbb E( X)$ for the expectation of a random variable $X$.
	We describe the asymptotic notation for functions $f=f(m)$ and
	$g=g(m)>0$ as $m\to\infty$.
	The relation $f=O_{r,\varepsilon}(g)$ means that, for each fixed choice
	of $r$ and $\varepsilon$, there are constants $C>0$ and $m_0$, depending
	only on $r$ and $\varepsilon$, such that $\lvert f(m)\rvert\le Cg(m)$
	for all $m\ge m_0$. Other subscripts specify the permitted dependence
	on fixed parameters in the same way.
	When no subscript is given, the constants may depend on the fixed
	parameters of the argument, but not on $m$.
	The notation $f=o(g)$ means that $f(m)/g(m)\to0$ as $m\to\infty$.
	The parameters specified as fixed, such as $r$ and $\varepsilon$, do not vary
	with $m$.
	For example, $O_r(m)$ denotes a quantity bounded in absolute value by
	$C_r m$, while $O_{r,\varepsilon}(1)$ denotes a quantity bounded by a constant
	depending only on $r$ and $\varepsilon$. If a variable other than $m$ tends to infinity,
	we specify that variable in the argument.
	
	For a finite set $\Omega$, $\mathbb R^\Omega$ denotes the real
	vector space of column vectors with coordinates indexed by $\Omega$.
	For a function $f:\Omega\to\mathbb R$, we write
	$\mathbf{f}=(f(S))_{S\in\Omega}$ for the column vector of its values.
	More explicitly, after fixing an ordering
	$\Omega=\{S_1,\ldots,S_t\}$, this means $
	\mathbf{f}=
	\begin{pmatrix}
		f(S_1)\\
		\vdots\\
		f(S_t)
	\end{pmatrix}$ and $
	(\mathbf{f})_{S_i}=f(S_i)$ for $ i\in[t]$.	We use bold symbols for vectors and ordinary symbols
	for the corresponding functions; $(\mathbf{f})_S$ denotes the
	coordinate of $\mathbf{f}$ indexed by $S$.
	For $\mathcal E\subseteq\Omega$, its indicator function $\chi_{\mathcal E}$
	equals one on $\mathcal E$ and zero elsewhere. We also write
	$\chi_{\{P\}}$ for one if the condition $P$ holds and zero otherwise.
	We reserve $\mathbf{1}$ for the all-one vector in the ambient coordinate
	space; in $\mathbb R^\Omega$, this means $(\mathbf{1})_S=1$ for every
	$S\in\Omega$.  We use the
	unnormalized inner product $
	\langle\mathbf{f},\mathbf{g}\rangle
	=\sum_{S\in\Omega}(\mathbf{f})_S(\mathbf{g})_S
	=\sum_{S\in\Omega}f(S)g(S)$ and its associated norm
	$\lVert \mathbf{f}\rVert =\langle\mathbf{f},\mathbf{f}\rangle^\frac{1}{2}$.
	For a matrix $A=(A_{S,R})_{S,R\in\Omega}$
	whose rows and columns are indexed by $\Omega$,
	the product $A\mathbf{f}$ is defined coordinatewise by $
	(A\mathbf{f})_S=\sum_{R\in\Omega}A_{S,R}(\mathbf{f})_R
	=\sum_{R\in\Omega}A_{S,R}f(R)$ for $S\in\Omega$.
	We denote by  $\mathscr{A}:\mathbb R^\Omega\to\mathbb R^\Omega$
	the linear transformation defined by $\mathscr{A}(\mathbf{x})=A\mathbf{x}$.
	A subspace $L\subseteq\mathbb R^\Omega$ is \emph{$\mathscr{A}$-invariant}
	if $\mathscr{A}(\mathbf{x})\in L$ for every $\mathbf{x}\in L$, that is,
	$\mathscr{A}(L)\subseteq L$.
	The span of a family of vectors is the subspace consisting of all
	their linear combinations.
	For subspaces $L_1,\ldots,L_t\subseteq\mathbb R^\Omega$, their sum
	is a \emph{direct sum}, denoted by $L_1\oplus\cdots\oplus L_t$, if
	every vector in the sum has a unique representation
	$\mathbf{x}_1+\cdots+\mathbf{x}_t$ with $\mathbf{x}_i\in L_i$.
	The direct sum is \emph{orthogonal} if
	$\langle\mathbf{x}_i,\mathbf{x}_j\rangle=0$ whenever
	$\mathbf{x}_i\in L_i$, $\mathbf{x}_j\in L_j$, and $i\ne j$.
	Vectors $\mathbf{f},\mathbf{g}\in\mathbb R^\Omega$ are orthogonal if
	$\langle\mathbf{f},\mathbf{g}\rangle=0$.
	A vector is orthogonal to a subspace
	if it is orthogonal to every vector in that subspace.
	A scalar $\lambda$ is an eigenvalue of $A$ if
	$A\mathbf{f}=\lambda\mathbf{f}$ for some nonzero vector $\mathbf{f}$,
	called an eigenvector.
	The corresponding eigenspace is
	$\{\mathbf{f}\in\mathbb R^\Omega
	\mid A\mathbf{f}=\lambda\mathbf{f}\}$.

	The extremal problem underlying our results is the Erd\H{o}s matching
	conjecture~\cite{Erdos1965}. For $n\ge rs$, it asserts that every
	$r$-graph $\cH$ on $n$ vertices with $\nu(\cH)\le s-1$ satisfies
	\[
	\lvert \cH\rvert \le\max\left\{\binom{rs-1}{r},
	\binom nr-\binom{n-s+1}{r}\right\}.
	\]
	The first term is attained by a complete $r$-graph on $rs-1$ vertices,
	with the remaining vertices isolated; the second is attained by all
	$r$-sets intersecting a fixed $(s-1)$-set. The first construction has too
	few nonisolated vertices to contain an $s$-matching. In the second,
	disjoint edges must use distinct vertices of the fixed set, which
	again excludes an $s$-matching. The clique construction can be larger
	near $n=rs$, whereas the fixed-cover construction is larger when
	$n$ is sufficiently large relative to $s$. The graph case follows
	from the Erd\H{o}s--Gallai theorem (Lemma~\ref{lem-eg}). For $s=2$, the
	problem is the study of intersecting families: the
	Erd\H{o}s--Ko--Rado theorem gives the maximum size and the
	Hilton--Milner theorem gives the corresponding bound for $n>2r$
	when there is no common vertex~\cite{EKR1961,HM1967}.
	
Bounds for general uniformity in the above conjecture were developed by
    Bollob\'as, Daykin and Erd\H{o}s~\cite{BDE1976} and
    Huang, Loh and Sudakov~\cite{HLS2012}.
    Frankl~\cite{Frankl2013} proved the fixed-cover bound when
    $n\ge(2r-1)s-r+1$.
	For fixed $r$ and sufficiently large $s$, Frankl and Kupavskii \cite{FK2022}
	established the range $n\ge(5r-2)(s-1)/3$.
	More recently, Cao, Liu and Zhang \cite{CLZ2026} proved the exact fixed-cover bound,
	including its equality case, for $n\ge(r+1)(s-1)$ and sufficiently
	large $s$. Their stability theorem also describes
	families whose size is close to this bound. These results concern
	edge counts under a matching-number constraint.
	
	The range close to $rs$ has a different extremal character.
	Frankl proved a range adjoining the clique construction;
	Kolupaev and Kupavskii extended this direction for almost perfect
	matchings~\cite{Frankl2017Range,KK2023}. In uniformity three,
	Frankl, R\"odl and Ruci\'nski and \L{}uczak and Mieczkowska obtained
	substantial ranges, before Frankl proved the full
	conjecture~\cite{FRR2012,LM2014,Frankl2017Max}.
	In uniformity four, Frankl, Lu, Ma and Wu proved the range
	$n\ge5(s-1)$ for sufficiently large $s$, together with a stability
	theorem~\cite{FLMW2026}. Hou, Hu and Liu subsequently proved the
	full four-uniform conjecture for $s\ge6005$ in the present
	notation~\cite{HHL2026}. The distinction between the clique and
	fixed-cover constructions is therefore central to the search for
	an exact range of orders.
	
	Shifting is a classical tool that imposes order on a set family while preserving
	its size and not increasing its matching number~\cite{Frankl1987}.
	The challenge for a degree-sum problem is to use this structural
	information while preserving the degree inequalities needed in the
	argument. An $r$-graph $\cH$ on the vertex set $
	[n]$ is \emph{fully shifted} if, for every
	$F\in\cH$ and every pair of vertices $i,j\in[n]$ satisfying $
	j\in F$, $i\notin F$ and $ i<j$,
	we have $
	(F\setminus j)\cup i\in\cH$.
	
	Degree conditions replace a global edge count by information about
	individual vertices. Huang and Zhao proved that, when
	$n\ge3r^2s$, minimum vertex degree greater than $D_{s-1}(n)$ forces
	an $s$-matching~\cite{HZ2017}; degree refinements of the
	Hilton--Milner theorem were obtained by Frankl, Han, Huang and
	Zhao~\cite{FHHZ2018}. Ore-degree conditions instead control sums
	of degrees on nonedges and need not impose the same lower bound
	on each vertex. The terminology recalls Ore's degree-sum criterion
	for Hamiltonian graphs~\cite{Ore1960}.
	Balogh, Palmer and Raeisi proved that
	$\sigma_r(\cH)>rD_{s-1}(n)$ forces an $s$-matching when
	$n\ge3r^2(s-1)$ and conjectured that $n>2rs$ suffices~\cite{BPR2026}.
	Their degree-sum threshold is attained by the fixed-cover
	construction. Our aim is to retain this threshold while reducing
	the required number of vertices. We state their conjecture below.

	\begin{conjecture}[\cite{BPR2026}]\label{conj-bpr}
		Let $r\ge3$ and $s\ge1$.  If $\cH$ is an $r$-graph on $n>2rs$ vertices
		and $\sigma_r(\cH)>rD_{s-1}(n)$, then $\nu(\cH)\ge s$.
	\end{conjecture}
	
	The degree threshold is best possible: the family of all $r$-sets
	intersecting a fixed $(s-1)$-set has matching number $s-1$ and Ore-degree
	$rD_{s-1}(n)$, while this example does not determine the best condition on $n$. Our first result below considers the range of $n$ for every uniformity $r\ge3$.
	
	\begin{theorem}\label{thm-main}
		For every integer $r\ge3$, there exists an integer $s_0=s_0(r)\ge3$
		such that, for all $s\ge s_0$, if $\cH$ is an $r$-graph
		on  $n\ge rs+\lfloor\frac{(r-1)^{r-1}}{r^{r-2}}\rfloor$ vertices and $\sigma_r(\cH)>rD_{s-1}(n)$, then $\nu(\cH)\ge s$.
	\end{theorem}
	
	The construction in Section~\ref{sec-additive} shows
	that the order bound is sharp for all sufficiently large $s$.
	For every integer $s\ge1$, the sufficient order
	$n\ge rs+C(r)$ holds with
	$C(r)=\lfloor\frac{(r-1)^{r-1}}{r^{r-2}}\rfloor+(r-1)(s_0(r)-3)$; see
	Corollary~\ref{cor-all-s}.
	
	For $r=3$, we obtain the following explicit bound on $n$, better than that stated in Conjecture~\ref{conj-bpr}.
	\begin{theorem}\label{thm-r3}
		Let $s\ge2$ and $\cH$ be a $3$-graph on $n$ vertices with
		$\sigma_3(\cH)>3D_{s-1}(n)$. If
		\[
		n\ge
		\begin{cases}
			7, & s=2,\\
			3s+9, & 3\le s\le5,\\
			\left\lceil(7s+11)/2\right\rceil, & s\ge6,
		\end{cases}
		\]
		then $\nu(\cH)\ge s$.
		In particular, Conjecture~\ref{conj-bpr} holds for $r=3$.
	\end{theorem}

	The paper is organized as follows.
	In Section~\ref{sec-prelim} we list some known results and establish the structural
	lemmas for use in the subsequent proofs. Section~\ref{sec-main} proves
	Theorem~\ref{thm-main} and Section~\ref{sec-r3} proves
	Theorem~\ref{thm-r3}. Section~\ref{sec-additive}
 shows that the order bound in Theorem~\ref{thm-main} is sharp,
	gives a sufficient order for Conjecture \ref{conj-bpr} of the form $rs+C(r)$ for every $s\ge1$,
 and asks whether the same sharp additive term suffices uniformly for every $s\ge1$.
	
	\section{Preliminaries}\label{sec-prelim}\label{sec-tools}
	In this section we recall the results needed in the proofs and establish some structural lemmas.
	We first recall two Ore-degree results.
	\begin{theorem}[{\cite{BPR2026}}]
		\label{thm-bpr}
		Let $r\ge3$ be an integer and $\cH$ be an $r$-graph on $n$ vertices.
		\begin{enumerate}
			\item If $s\ge2$, $n\ge3r^2(s-1)$ and
			$\sigma_r(\cH)>rD_{s-1}(n)$, then $\nu(\cH)\ge s$.
			\item If $n\ge2r+2$ and $\cH$ is intersecting, then
			$\sigma_r(\cH)\le r\binom{n-2}{r-2}$,
			with equality if and only if $\cH$ is a $1$-star.
		\end{enumerate}
	\end{theorem}
	
	\begin{lemma}[\cite{EG1959}]\label{lem-eg}
		Let $b\ge0$ be an integer, and let $G$ be a graph on $q\ge2b+1$ vertices.
		If $\nu(G)\le b$, then
		\begin{equation}\label{eq-eg}
			\lvert E(G)\rvert\le
			\max\left\{\binom{2b+1}{2},\ bq-\frac{b(b+1)}2\right\}.
		\end{equation}
	\end{lemma}
	
	For $r\ge3$, let $
	p_r=1-\left(1-\frac1r\right)^{r-1}$, $
	\alpha_r=\frac{(r-1)^{r-1}}{r^{r-2}}=r(1-p_r)$ and $
	c_r=\lfloor\alpha_r\rfloor$. Since $(1+1/(r-1))^{r-1}>2$,
	$1<2(r-1)/r\leq\alpha_r<r/2<r-1$.
	Thus $p_r>1/2$ and $1\le c_r\le r-2$.
	Moreover, since
	$\gcd(r-1,r)=1$ and $r^{r-2}>1$, $\alpha_r$ is not an integer.
	
	\begin{lemma}\label{lem-binomial}
		For $r\ge3$, $s\ge2$ and $n\ge rs$, $
		D_{s-1}(n)>p_r\binom{rs-2}{r-1}
		>\frac12\binom{rs-2}{r-1}$.
	\end{lemma}
	\begin{proof}
		Since $D_{s-1}(n)$ is increasing in $n$, it suffices to consider $n=rs$.
		Write $z=r(s-1)$ and $q=(r-1)/r$. Let $
		Q(u)=\prod_{j=1}^{r-1}(u+j)$ and $
		P(z)=z\prod_{j=1}^{r-2}(z+j)$.
		Then $
		(r-1)!D_{s-1}(rs)=Q(z)-Q(qz)$ and $
		(r-1)!\binom{rs-2}{r-1}=P(z)$.
		For $0<u\le z$, each term of $Q'(u)$ is at least
		$\prod_{j=1}^{r-2}(u+j)$ and some are strictly larger. Also
		$u+j\ge (u/z)(z+j)$ for every $j\ge1$. Thus
		\[
		Q'(u)>(r-1)\prod_{j=1}^{r-2}(u+j)
		\ge (r-1)\left(\frac uz\right)^{r-2}\frac{P(z)}z.
		\]
		Integrating over $qz\le u\le z$ gives
		$Q(z)-Q(qz)>(1-q^{r-1})P(z)=p_rP(z)$, which implies $D_{s-1}(rs)=\dfrac{Q(z)-Q(qz)}{(r-1)!}>\dfrac{p_rP(z)}{(r-1)!}=p_r\binom{rs-2}{r-1}$.
		Finally the result follows from $p_r>1/2$ for $r\geq3$.
	\end{proof}
	
	\begin{definition}\label{def-nonext}
		Let $\eta>0$ be a constant. An $N$-vertex $r$-graph $\cH$ is $\eta$-nonextremal if $
		\lvert \cH[B]\rvert \ge\eta N^r$ for every $B\subseteq V(\cH)$ with $
		\lvert B\rvert =\left\lfloor(1-1/r)N\right\rfloor$.
	\end{definition}
	\begin{remark}\label{rem-monotone}
		Let $\eta>0$ and $\cH\subseteq\cK$ be $r$-graphs on the same
		vertex set $V$ of size $N$. If $\cH$ is $\eta$-nonextremal,
		then $\cK$ is also $\eta$-nonextremal. Indeed, for every
		$B\subseteq V$ with $\lvert B\rvert=\lfloor(1-1/r)N\rfloor$,
		we have $\cH[B]\subseteq\cK[B]$ and hence $
		\lvert\cK[B]\rvert\ge\lvert\cH[B]\rvert\ge\eta N^r$.
		The conclusion follows from Definition~\ref{def-nonext}.
	\end{remark}
	\begin{lemma}\label{lem-basic}
		Let $\cH$ be an $r$-graph on a vertex set $V$ of size $n$ with
		$n\ge r\ge3$. Let $Z\subseteq V$ with $z=\lvert Z\rvert$ such that $n-z\ge r$. Then $\cH$ has the following properties.
		\begin{enumerate}
			\item Let $D>0$. If $\sigma_r(\cH)>rD$, then
			$\{v\in V\mid d_\cH(v)\le D\}$ induces a complete $r$-graph.
			If $\sigma_r(\cH)\ge rD$, the same holds for
			$\{v\in V\mid d_\cH(v)<D\}$.
			\item For every $v\in V\setminus Z$, $
			d_{\cH-Z}(v)\ge d_\cH(v)-z\binom{n-2}{r-2}$.
			If $\cH-Z$ has a nonedge, then $
			\sigma_r(\cH-Z)\ge\sigma_r(\cH)-rz\binom{n-2}{r-2}$.
			\item Let $\eta>0$ and $\cF$ be an $r$-graph on $V$.
			If $\cH$ is $\eta$-nonextremal and $
			(z+1)B_r(n)+\lvert\cF\rvert\le\frac\eta2 n^r$,
			then $(\cH-\cF)-Z$ is $(\eta/2)$-nonextremal.
			In particular, for fixed $r,\eta$, there exists $\varepsilon>0$
			such that, if $\cH$ is $\eta$-nonextremal,
			$z\le\varepsilon n$ and $\lvert\cF\rvert=o(n^r)$, then the same conclusion holds for all sufficiently large
			$n$.
		\end{enumerate}
	\end{lemma}
	\begin{proof}
		If every vertex of an $r$-set $S$ has degree at most
		$D$, then $\sum_{v\in S}d_\cH(v)\le rD$.
		Thus $S$ cannot be a nonedge when $\sigma_r(\cH)>rD$.
		If all these degrees are strictly less than $D$, then
		$\sum_{v\in S}d_\cH(v)<rD$, so $S$ cannot be a nonedge when
		$\sigma_r(\cH)\ge rD$. This proves both assertions of (i).
		For $v\in V\setminus Z$, $
		d_\cH(v)-d_{\cH-Z}(v)=\left\lvert\{F\in\cH\mid v\in F,\ F\cap Z\ne\emptyset\}\right\rvert\le\sum_{u\in Z}
		\left\lvert\{F\in\cH\mid\{u,v\}\subseteq F\}\right\rvert
		\le z\binom{n-2}{r-2}$.
Suppose that $\cH-Z$ has a nonedge. For every
        $S\in\binom{V\setminus Z}{r}\setminus(\cH-Z)$,
        we have $S\notin\cH$ and hence $
		\sum_{v\in S}d_{\cH-Z}(v)
		\ge\sum_{v\in S}d_\cH(v)-rz\binom{n-2}{r-2}
		\ge\sigma_r(\cH)-rz\binom{n-2}{r-2}$. It follows that
		$
		\sigma_r(\cH-Z)\ge\sigma_r(\cH)-rz\binom{n-2}{r-2}$.
		
		For (iii), let $\cH'=(\cH-\cF)-Z$ and
		$X\subseteq V\setminus Z$ with
		$\lvert X\rvert=\lfloor(1-1/r)(n-z)\rfloor$.
		Extend $X$ to a set $B\subseteq V$ of size
		$\lfloor(1-1/r)n\rfloor$. Then $
		\lvert B\setminus X\rvert
		=\left\lfloor\left(1-\frac1r\right)n\right\rfloor
		-\left\lfloor\left(1-\frac1r\right)(n-z)\right\rfloor
		\le\left(1-\frac1r\right)z+1\le z+1$.
		Every edge in $\cH[B]\setminus\cH[X]$ intersects $B\setminus X$,
		so $
		\lvert\cH[B]\rvert-\lvert\cH[X]\rvert
		\le\sum_{v\in B\setminus X}d_\cH(v)
		\le(z+1)B_r(n)$.
		Since $X\cap Z=\emptyset$, we have
		$\cH'[X]=\cH[X]\setminus\cF$.
		Since $\cH$ is $\eta$-nonextremal and $(z+1)B_r(n)+\lvert\cF\rvert\le\frac\eta2 n^r$, we obtain
		\[
		\begin{aligned}
			\lvert\cH'[X]\rvert
			&\ge\lvert\cH[B]\rvert-(z+1)B_r(n)-\lvert\cF\rvert\ge\eta n^r-(z+1)B_r(n)-\lvert\cF\rvert\ge\frac\eta2 n^r\ge\frac\eta2(n-z)^r,
		\end{aligned}
		\]
		which implies $\cH'$ is $(\eta/2)$-nonextremal.
		For the final assertion, choose
		$\varepsilon=\min\{1/2,\eta(r-1)!/4\}$. Since $B_r(n)\le n^{r-1}/(r-1)!$,
		if $z\le\varepsilon n$ and $\lvert\cF\rvert=o(n^r)$, then
		\[
		\begin{aligned}
			(z+1)B_r(n)+\lvert\cF\rvert
			&\le\frac{\varepsilon}{(r-1)!}n^r
			+\frac{n^{r-1}}{(r-1)!}+o(n^r)\le\frac\eta4 n^r+o(n^r)
			\le\frac\eta2 n^r
		\end{aligned}
		\]
		for sufficiently large $n$. Thus the hypothesis of~(iii) holds and hence $\cH'$ is $(\eta/2)$-nonextremal.
	\end{proof}
	
	\begin{definition}\label{def-good} Let $Y,W$ be two disjoint sets with $\lvert Y\rvert +\lvert W\rvert=N$ and $\lambda>0$ be a constant. Suppose that $\cH\subseteq\binom{Y\cup W}{r}$ and $\cT(Y,W)=\{F\in\binom{Y\cup W}{r}\mid F\cap Y\neq\emptyset,\;F\cap W\neq\emptyset\}$.
		A vertex $v\in V(\cH)$
		is called $\lambda$-good with respect to $\cT(Y,W)$ if $\lvert\{F\in\cT(Y,W)\setminus\cH\mid v\in F\}\rvert\leq\lambda N^{r-1}$.
	\end{definition}

	\begin{theorem}[\cite{LYY2021}]\label{thm-lyy}
		Let $r\geq 3$ and $0<\lambda<(8^{r-1}r^{5(r-1)}r!)^{-1}$. Let $\cH\subseteq\binom{Y\cup W}{r}$ be an $r$-graph with $Y\cap W=\emptyset$ and $\lvert Y\rvert +\lvert W\rvert =N\ge8r^6$ such that
		$N/(2r^5)<\lvert W\rvert\le N/r$. If every vertex $v\in V(\cH)$ is
		$\lambda$-good with respect to $\cT(Y,W)$, then
		$\nu(\cH)\ge \lvert W\rvert$.
	\end{theorem}
	
	\begin{lemma}[{\cite{HLWY2026}}]
		\label{lem-han-degree}
		Fix integers $r\ge3$ and $1\le d\le r-2$ such that
		$(1-1/r)^{r-d}<1/2$, and let $\eta,\gamma>0$.
		There exist $\alpha,\beta>0$ and $N_0$ such that the following
		holds for every $N\ge N_0$ with $r\mid N$.
		Let $\cH$ be a fully shifted $r$-graph on $[N]$ satisfying
		$\delta(\cH)>\gamma N^{r-1}$ and $
		\lvert \cH[B]\rvert >\eta\binom Nr$ for every $B\subseteq[N]$ with $\lvert B\rvert =(1-1/r)N$.
		If there is a set $S\subseteq[N]$ with $\lvert S\rvert \le\beta N$ and $
		\delta_d(\cH-S)\ge
		\binom{N-d}{r-d}-\binom{N-N/r}{r-d}-\alpha N^{r-d}$,
		then $\cH$ has a perfect matching.
	\end{lemma}
	The following corollary is the vertex-degree specialization of
	Lemma~\ref{lem-han-degree}.
	\begin{corollary}
		\label{cor-han-shift}
		Fix $r\ge3$ and $\eta,\gamma>0$. There are positive constants
		$\alpha,\beta$ and $N_0$ such that the following holds for all
		$N\ge N_0$ with $r\mid N$. Let $\cH$ be an $r$-graph on $N$ vertices such that $\cH$ is fully shifted,
		$\eta$-nonextremal and $ \delta(\cH)\ge\gamma N^{r-1}$.
		If some $S\subseteq V(\cH)$ satisfies $\lvert S\rvert \le\beta N$ and $\delta(\cH-S)\ge(p_r-\alpha)B_r(N)$ where $
		p_r=1-\left(1-\frac1r\right)^{r-1}$,
		then $\cH$ has a perfect matching.
	\end{corollary}
	
	\begin{proof}
		We apply Lemma~\ref{lem-han-degree} with $d=1$ and
		$\gamma/2$ in place of $\gamma$, and denote its constants by
		$\alpha_0,\beta,N_0$. The condition
		$(1-1/r)^{r-1}<1/2$ holds for $r\ge3$.
		Moreover, for every
		$B\subseteq V(\cH)$ with $\lvert B\rvert=(1-1/r)N$, $
		\lvert\cH[B]\rvert\ge\eta N^r>\eta\binom Nr$.
		For fixed $r$ as $N\to\infty$,
		\begin{align}
			B_r(N)
			&=\frac{1}{(r-1)!}\prod_{j=1}^{r-1}(N-j)=\frac{N^{r-1}}{(r-1)!}+O_r(N^{r-2}),\label{eq-bin-asym}\\
			\binom{N-N/r}{r-1}
			&=\frac{1}{(r-1)!}
			\prod_{j=0}^{r-2}(N-\frac Nr-j)\nonumber\\
			&=\frac{1}{(r-1)!}
			\prod_{j=0}^{r-2}
		\bigg((1-\frac1r)N-j\bigg)=\frac{(1-\frac1r)^{r-1}
				N^{r-1}}{(r-1)!}
			+O_r(N^{r-2}).\nonumber
		\end{align}
		Subtracting the two expansions with $p_r=1-\left(1-\frac1r\right)^{r-1}$, we obtain
		\begin{align*}
			B_r(N)
			-\binom{N-N/r}{r-1}	&=\left[1-\left(1-\frac1r\right)^{r-1}\right]
			\frac{N^{r-1}}{(r-1)!}
			+O_r(N^{r-2})\\
			&=p_r\frac{N^{r-1}}{(r-1)!}
			+O_r(N^{r-2})=p_rB_r(N)
			+O_r(N^{r-2}).
		\end{align*}
		The last equality follows from $
		\frac{N^{r-1}}{(r-1)!}
		=B_r(N)+O_r(N^{r-2})$.
		Moreover, by \eqref{eq-bin-asym} we have $
		N^{r-1}
		=(r-1)!B_r(N)+O_r(N^{r-2})$.
		It follows that
		\begin{align}\label{eq-deg-asym}
			B_r(N)-\binom{N-N/r}{r-1}
			-\alpha_0 N^{r-1}&=\bigl(p_r-\alpha_0(r-1)!\bigr)
			B_r(N)
			+O_{r,\alpha_0}(N^{r-2}).
		\end{align}
		Choose $\alpha>0$ such that $
		\alpha<\min\left\{p_r,\frac{\alpha_0(r-1)!}{2}\right\}$.
		Since $
		N^{r-2}
		=o\left(B_r(N)\right)$, we may enlarge $N_0$ so that the error term of \eqref{eq-deg-asym}
		is at most $
		\bigl(\alpha_0(r-1)!-\alpha\bigr)
		B_r(N)$
		for all $N\ge N_0$. Hence $
		B_r(N)-\binom{N-N/r}{r-1}
		-\alpha_0 N^{r-1}\le(p_r-\alpha)B_r(N)$,
		as required.
		All assumptions of Lemma~\ref{lem-han-degree} with $d=1$ hold, so $\cH$ has a perfect matching.
	\end{proof}
	\begin{theorem}[{\cite{HLWY2026}}]
		\label{thm-han-nonext}
		Fix $r\ge3$ and $\eta>0$. There exist $\zeta>0$ and $N_0$ such
		that the following holds for every $N\ge N_0$ with $r\mid N$.
		Let $\cH$ be an $\eta$-nonextremal $r$-graph on $N$ vertices.
		If $
		\delta(\cH)\ge(p_r-\zeta)B_r(N)$,
		then $\cH$ has a perfect matching.
	\end{theorem}
	
	\begin{theorem}[{\cite{Morris2025}}]\label{thm-morris}
		Let $r\ge3$ be an integer and $0<\theta<1/(2r^2)$.
		There exists an integer $N_0$ such that the following holds for
		all $N\ge N_0$.
	Let $\cH$ be an $r$-graph on $N$ vertices. For a graph $J$ on $V(\cH)$, let $
		\cH_J=\left\{F\in\cH\mid\binom F2\cap E(J)\ne\emptyset\right\}$
		be an $r$-graph on $V(\cH)$. If $\cH-\cH_J$ has a perfect
		fractional matching for every graph $J$ on $V(\cH)$ with
		$\Delta(J)\le N^{r^2\theta}$, then $\cH$ has a matching $M$
		such that $N-\lvert V(M)\rvert\le N^{1-\theta}$.
	\end{theorem}
	
	We record the spectral description of the Kneser graph.
	\begin{lemma}[\cite{HZ2017}]
		\label{lem-kg-spec}
		Let $r\ge2$, $n\ge2r+1$ be integers and $\Omega=\binom{[n]}r$.
		Let $A$ be the adjacency matrix of the Kneser graph $KG(n,r)$,
		so $A_{S,R}=\chi_{\{S\cap R=\emptyset\}}$ for $S,R\in\Omega$.
		For each $v\in[n]$, let $u_v:\Omega\to\mathbb R$ be the function
		$u_v(S)=\chi_{\{v\in S\}}$, and let
		$\mathbf{u}_v=(u_v(S))_{S\in\Omega}$ be its column vector of values.
		Then the eigenvalues of $A$ are $
		\lambda_j=(-1)^j\binom{n-r-j}{r-j}$ for $0\le j\le r$
		with respective multiplicities $\binom nj-\binom n{j-1}$,
		where $\binom n{-1}=0$.
		The eigenspaces $W_0$ and $W_1$ corresponding to $\lambda_0$ and
		$\lambda_1$ satisfy $
		W_0=\operatorname{span}\{\mathbf{1}\}$ and $
		W_1=\bigl\{\mathbf{f}\in\operatorname{span}\{\mathbf{u}_v\mid v\in[n]\}
		\mid \langle \mathbf{f},\mathbf{1}\rangle=0\bigr\}$, respectively.
	\end{lemma}

	We also recall the Erd\H{o}s--Ko--Rado theorem, including
	its equality characterization for $n>2r$.
	\begin{theorem}[\cite{EKR1961}]
		\label{thm-ekr}
		Let $r\ge2$, $n\ge2r$ be integers and $\cH$ be an
		intersecting $r$-graph on $n$ vertices. Then $
		\lvert \cH\rvert \le B_r(n)$.
		If $n>2r$, equality holds if and only if $\cH$ is a $1$-star.
	\end{theorem}

	\begin{lemma}\label{lem-wt-shift}
		Let $n\ge r\ge2$ be integers and let $w:[n]\to[0,1]$ satisfy $
		w(1)\ge\cdots\ge w(n)$ and $
		\sum_{v=1}^n w(v)<\frac nr$. Let $\cH=\left\{F\in\binom{[n]}r\mid\sum_{v\in F}w(v)\ge1\right\}$ be an $r$-graph  on $[n]$. Then $\cH$ is fully shifted and has no perfect matching.
	\end{lemma}
	\begin{proof}
		Let $F\in\cH$ and $i,j\in[n]$ satisfy
		$j\in F$, $i\notin F$ and $i<j$. Since $w(i)\ge w(j)$,
		\[
		\sum_{v\in(F\setminus j)\cup i}w(v)
		=\sum_{v\in F}w(v)-w(j)+w(i)\ge1.
		\]
		Thus $(F\setminus j)\cup i\in\cH$ for each edge $F\in\cH$, so $\cH$ is fully shifted.
		If $M$ is a perfect matching of $\cH$, then
		\[
		\frac nr=\lvert M\rvert
		\le\sum_{F\in M}\sum_{v\in F}w(v)
		=\sum_{v=1}^n w(v)<\frac nr,
		\]
		contrary to $
		\sum_{v=1}^n w(v)<\frac nr$.
	\end{proof}

	We use the following forms of the Chernoff bounds.
	\begin{lemma}[{\cite{MU2005}}]
		\label{lem-chernoff}
		Let $m\ge0$ be an integer and let $X_1,\ldots,X_m$ be independent
		random variables taking values in $\{0,1\}$. Let
		$X=\sum_{i=1}^m X_i$ and $\mu=\mathbb E(X)$. Then
		\begin{enumerate}
			\item for every $0<a<1$, $
			\mathbb P\bigl(\lvert X-\mu\rvert\ge a\mu\bigr)
			\le2\exp\left(-\frac{a^2\mu}{3}\right)$;
			\item for every $t>0$ with $t\ge\mu$, $
			\mathbb P(X\ge t)\le\left(\frac{\mathrm e\mu}{t}\right)^t$.
		\end{enumerate}
	\end{lemma}
	
	\section{Proof of Theorem~\ref{thm-main}}
	\label{sec-main}\label{sec-sharp}
	Throughout this section, let $r\geq3$ be an integer. Unless otherwise stated,
	asymptotic estimates are taken as $n\to\infty$, with all parameters
	already chosen held fixed. Recall that $
	p_r=1-\left(1-\frac1r\right)^{r-1}$ and $B_r(n)=\binom{n-1}{r-1}$.
	\begin{lemma}\label{lem-robust-frac}
		Let $r\geq3$ be an integer. For every $\eta,\gamma>0$,
		there are $\zeta,\varepsilon>0$ and $n_0$ such that, for all $n\ge n_0$
		with $r\mid n${\pu ,} the following holds. If $\cH$ is an
		$\eta$-nonextremal $n$-vertex $r$-graph satisfying
		$\delta(\cH)\ge\gamma n^{r-1}$ and
		$\sigma_r(\cH)\ge r(p_r-\zeta)B_r(n)$,
		then  $\cH-\cF$ has a perfect
		fractional matching for every $r$-graph $\cF$ on $V(\cH)$ with
		$\Delta(\cF)\le\varepsilon n^{r-1}$.
	\end{lemma}
	
	\begin{proof}
		Fix $r,\eta,\gamma$, and choose $\alpha,\beta,N_0$ from
		Corollary~\ref{cor-han-shift} with $\eta/2$ and $\gamma/2$
		in place of $\eta$ and $\gamma$.
		Decreasing $\alpha$ only strengthens the minimum-degree hypothesis
		of that corollary, so we may assume $0<\alpha<p_r$. Let $
		a=\frac{\alpha}{4}$, $\zeta=\frac a2$, $
		b=\min\left\{\frac\beta2,
		\frac{\alpha}{8(r-1)},\frac{1-p_r}{4(r-1)}\right\}$
		and $
		\varepsilon=\frac12\min\left\{\frac\gamma4,\frac{r\eta}4,
		\frac{a}{8(r-1)!},\frac{b^{r-1}}{4(r-1)!}\right\}$.
		Then $a,\zeta,b,\varepsilon$ are positive
		and depend only on $r,\eta,\gamma$.
		Choose an integer $n_0\ge\max\{N_0,r\}$ sufficiently large in
		terms of these constants for the estimates below. Thus $n_0$
		also depends only on $r,\eta,\gamma$ and can be fixed before
		choosing $\cH$ or $\cF$. Let $n\ge n_0$ satisfy $r\mid n$ and $\cH$ be an $n$-vertex $r$-graph satisfying the assumptions.
		
		Suppose to the contrary that there exists an $r$-graph $\cF$ on $V(\cH)$ with
		$\Delta(\cF)\le\varepsilon n^{r-1}$ such that $\cH-\cF$ has no perfect fractional matching. For any fractional matching $f$ in $\cH-\cF$, we have
		\begin{align}\label{eq-frac-sum}
			\sum_{v\in V(\cH)}
			\sum_{\substack{F\in\cH-\cF\\v\in F}}f(F)&=
			\sum_{F\in\cH-\cF}
			\sum_{v\in F}f(F)=
			r\sum_{F\in\cH-\cF}f(F).
		\end{align}
		Since $
		\sum_{\substack{F\in\cH-\cF\\v\in F}}f(F)\le1$ for every $v\in V(\cH)$, the left-hand side of \eqref{eq-frac-sum} is at most $\sum_{v\in V(\cH)}1=n$. It follows that $\sum_{F\in\cH-\cF}f(F)\le\frac nr$, with equality if and only if $
		\sum_{\substack{F\in\cH-\cF\\v\in F}}f(F)=1$ for every $v\in V(\cH)$, which is the definition of a perfect
		fractional matching.
		Since $\cH-\cF$ has no perfect fractional matching,
		its maximum fractional matching weight is strictly less than
		$n/r$. By linear programming duality, there exists a fractional
		vertex cover $w:V(\cH)\to[0,1]$ such that $
		\sum_{v\in F}w(v)\ge1$ for every $F\in\cH-\cF$ and $
		\sum_{v\in V(\cH)}w(v)<\frac nr$.
		Relabel the vertices so that $w(1)\ge\cdots\ge w(n)$ and define $
		\cK:=\left\{F\in\binom{[n]}r\mid
		\sum_{v\in F}w(v)\ge1\right\}$.
		Then $\cH-\cF\subseteq\cK$. For every $v\in V(\cK)=[n]$,
		\[
		\begin{aligned}
			d_{\cK}(v)
			&\ge d_{\cH-\cF}(v)
			\ge d_{\cH}(v)-d_{\cF}(v)\ge\delta(\cH)-\Delta(\cF)
			\ge(\gamma-\varepsilon)n^{r-1}
			>\frac\gamma2 n^{r-1},
		\end{aligned}
		\]
		where the last inequality follows from $\varepsilon<\gamma/4$.
		Hence $\delta(\cK)>\frac\gamma2 n^{r-1}$. Moreover, $
		r\lvert\cF\rvert=\sum_{v\in[n]}d_{\cF}(v)\le n\Delta(\cF)
		\le\varepsilon n^r$. Thus $\lvert\cF\rvert\le\frac{\varepsilon}{r}n^r<\frac{\eta}{4}n^r$, where the last inequality holds since $
		\varepsilon<\frac{r\eta}{4}$. Since $r$ is fixed, we have $
		B_r(n)=\binom{n-1}{r-1}
		=\frac{n^{r-1}}{(r-1)!}+O_r(n^{r-2})
		=o(n^r)$. By Lemma~\ref{lem-basic}~(iii) with
		$Z=\emptyset$, $\cH-\cF$ is $(\eta/2)$-nonextremal
		for sufficiently large $n$. Since $\cH-\cF\subseteq\cK$ and
		$V(\cH-\cF)=V(\cK)=[n]$, $\cK$ is also $(\eta/2)$-nonextremal by Remark~\ref{rem-monotone}.
		
		Let $
		L=\{v\in[n]\mid d_\cK(v)<(p_r-a)B_r(n)\}$. We claim that $\binom{L}{r}\subseteq\cH$.
		Since $\cH-\cF\subseteq\cK$, for every $v\in L$,
		we have
		\[
		\begin{aligned}
			d_{\cH}(v)
			&\le d_{\cH-\cF}(v)+d_{\cF}(v)\le d_{\cK}(v)+\Delta(\cF)<(p_r-a)B_r(n)+\varepsilon n^{r-1}.
		\end{aligned}
		\]
		Since $
		B_r(n)
		=\frac{n^{r-1}}{(r-1)!}+O_r(n^{r-2})$, $
		\frac{B_r(n)}{n^{r-1}}
		\longrightarrow\frac{1}{(r-1)!}$ as $n\to\infty$.
		Combining this with $\varepsilon<\frac{a}{4(r-1)!}$ and $\zeta=\frac a2$, we have $d_{\cH}(v)
		<(p_r-\zeta)B_r(n)$ for every $v\in L$.
	Suppose to the contrary,
	$\binom{L}{r}\nsubseteq\cH$.
		Then there exists an $r$-subset
		$S\subseteq L$ such that $S\notin\cH$ and $
		\sum_{v\in S}d_{\cH}(v)<\sum_{v\in S}(p_r-\zeta)B_r(n)=r(p_r-\zeta)B_r(n)$.
		However, by assumption and $S\notin\cH$, $
		\sum_{v\in S}d_{\cH}(v)
		\ge\sigma_r(\cH)
		\ge r(p_r-\zeta)B_r(n)$,
		a contradiction. Therefore,
		\begin{equation}\label{eq-low-clique}
			\binom Lr\subseteq \cH.
		\end{equation}
		We next show the following claim.
		\begin{claim}\label{clm-low-set}
			$\lvert L\rvert<bn$.
		\end{claim}
		\noindent\textbf{Proof of Claim~\ref{clm-low-set}.}
		By Lemma~\ref{lem-wt-shift},
		$\cK$ is fully shifted. Then $
		d_{\cK}(1)\ge d_{\cK}(2)\ge\cdots\ge d_{\cK}(n)$  by the definition of a fully shifted $r$-graph.
	Suppose to the contrary,
	$\lvert L\rvert\ge bn$. Let $
		m=\lceil bn\rceil$ and $
		Q=\{n-m+1,\ldots,n\}$. Recall that $
		L=\left\{
		v\in[n]\mid
		d_{\cK}(v)<(p_r-a)B_r(n)
		\right\}$.
		It follows
		that if $v\in L$ and $u>v$, then $u\in L$. This implies $Q\subseteq L$.
		We first prove that  $\cK[Q]$ has no edge. Otherwise, there exists an edge
		$F\in\cK$ with $F\subseteq Q$.
		Let $v=\min F$. Then $
		v\ge n-m+1$. Since $\cK$ is fully shifted, $v\cup R\in\cK$ for every $R\in\binom{[v-1]}{r-1}$.
		Consequently, $
		d_{\cK}(v)\ge\binom{v-1}{r-1}\ge\binom{n-m}{r-1}$.
		Since $m=\lceil bn\rceil$, we have
		$m/n\longrightarrow b$ as $n\to\infty$.
		As $r$ is fixed,
		\[
		\begin{aligned}
			\frac{\binom{n-m}{r-1}}{B_r(n)}
			&=
			\frac{\binom{n-m}{r-1}}
			{\binom{n-1}{r-1}}=
			\prod_{j=0}^{r-2}
			\frac{n-m-j}{n-1-j}\longrightarrow (1-b)^{r-1}.
		\end{aligned}
		\]
		By Bernoulli's inequality and $b\le\frac{1-p_r}{4(r-1)}$, $
		(1-b)^{r-1}\ge1-(r-1)b\ge1-\frac{1-p_r}{4}>p_r>p_r-a$.
		Thus, for sufficiently large $n$, $
		d_{\cK}(v)\ge\binom{n-m}{r-1}
		>(p_r-a)B_r(n)$,
		contrary to $v\in Q\subseteq L$.
		Hence $\cK[Q]$ has no edge.
		
		On the other hand,
		since $Q\subseteq L$, $
		\binom{Q}{r}\subseteq\cH$ by \eqref{eq-low-clique}.
		Since $\cK[Q]$ has no edge and $
		\cH-\cF\subseteq\cK$, every $r$-subset of $Q$
		must belong to $\cF$, that is, $
		\binom{Q}{r}\subseteq\cF$.
		Since
		$\lvert Q\rvert=m$, for every $v\in Q$
		and every
		$S\in\binom{Q\setminus v}{r-1}$,
		we have $v\cup S\in\cF$.
		Thus
		\[
		\Delta(\cF)\geq d_{\cF}(v)
		\ge
		\left\lvert
		\binom{Q\setminus v}{r-1}
		\right\rvert
		=\binom{m-1}{r-1}.
		\]
		Since $m/n\longrightarrow b$, we have
		\[
		\frac{\binom{m-1}{r-1}}{n^{r-1}}
		=
		\frac{1}{(r-1)!}
		\prod_{j=1}^{r-1}\frac{m-j}{n}
		\longrightarrow
		\frac{b^{r-1}}{(r-1)!}.
		\]
		Combining this with $
		\varepsilon<\frac{b^{r-1}}{4(r-1)!}$, we therefore have $
		\binom{m-1}{r-1}
		>\varepsilon n^{r-1}$ for sufficiently large $n$.
		It follows that $
		\Delta(\cF)>\varepsilon n^{r-1}$,
		contradicting the assumption
		$\Delta(\cF)\le\varepsilon n^{r-1}$. This completes the proof of Claim~\ref{clm-low-set}.
		
		Recall that $
		L=\{v\in[n]\mid d_\cK(v)<(p_r-a)B_r(n)\}$. For every $v\in V(\cK)\setminus L$ and every $u\in L$, since $\lvert\{F\in\cK\mid v\in F,\;u\in F\}\rvert\leq \binom{n-2}{r-2}$, we have\[
		\begin{aligned}
			d_{\cK-L}(v)
			&\ge d_{\cK}(v)
			-\lvert L\rvert\binom{n-2}{r-2}\ge(p_r-a)B_r(n)
			-bn\binom{n-2}{r-2},
		\end{aligned}
		\]
		where the last inequality holds by $v\in V(\cK)\setminus L$ and Claim~\ref{clm-low-set}.
		Using $
		\binom{n-2}{r-2}
		=\frac{r-1}{n-1}B_r(n)$ and $
		\frac n{n-1}\le2$ for sufficiently large $n$,
		we obtain
		\[
		\begin{aligned}
			d_{\cK-L}(v)
			&\ge
			\left(
			p_r-a-b(r-1)\frac{n}{n-1}
			\right)B_r(n)\ge
			\left(p_r-a-2(r-1)b\right)B_r(n).
		\end{aligned}
		\]
		Recall that $
		a=\frac{\alpha}{4}$ and $
		b\le\frac{\alpha}{8(r-1)}$.
		Thus $
		a+2(r-1)b\le\frac{\alpha}{4}
		+2(r-1)\frac{\alpha}{8(r-1)}=\frac{\alpha}{2}
		\le\alpha$.
		It follows that $d_{\cK-L}(v)\ge(p_r-\alpha)B_r(n)$
		for every $v\in V(\cK)\setminus L$.
		Hence $
		\delta(\cK-L)\ge(p_r-\alpha)B_r(n)$.
		
		We have now established that $\cK$ is
		fully shifted, $(\eta/2)$-nonextremal
		and satisfies $
		\delta(\cK)\ge\frac{\gamma}{2}n^{r-1}$ and $
		\lvert L\rvert<bn\le\beta n$ with $
		\delta(\cK-L)\ge(p_r-\alpha)B_r(n)$.
		Then $\cK$ has a perfect matching by
		Corollary~\ref{cor-han-shift}
		
		However, every edge $F\in\cK$ satisfies $
		\sum_{v\in F}w(v)\ge1$ by the definition of $\cK$.
		If $M$ is a perfect matching in $\cK$, then its edges partition $V(\cK)$, and hence
		\[
		\begin{aligned}
			\sum_{v\in V(\cK)}w(v)
			&=\sum_{F\in M}\sum_{v\in F}w(v)\ge\sum_{F\in M}1=\frac nr,
		\end{aligned}
		\]
		which contradicts $
		\sum_{v\in V(\cK)}w(v)<\frac nr$.
		This
		completes the proof.
	\end{proof}
	
	For fixed $r,\eta,\gamma$, Lemma~\ref{lem-robust-frac}'s
	conclusion holds simultaneously for every $n\ge n_0$ with
	$r\mid n$, every $\cH$ satisfying its assumptions and every
	$\cF$ with $\Delta(\cF)\le\varepsilon n^{r-1}$. Thus the constants $\zeta,\varepsilon,n_0$ can be fixed before any hypergraph is chosen.
	
	\begin{lemma}\label{lem-almost}
		For every integer $r\ge3$ and $\eta,\gamma>0$, there exists
		$\zeta=\zeta(r,\eta,\gamma)>0$ such that, for every $\rho>0$,
		there is an integer $n_1=n_1(r,\eta,\gamma,\rho)$ with the following
		property. If $n\ge n_1$ satisfies $r\mid n$ and $\cH$ is an
		$\eta$-nonextremal $n$-vertex $r$-graph satisfying $
		\delta(\cH)\ge\gamma n^{r-1}$ and $
		\sigma_r(\cH)\ge r(p_r-\zeta)B_r(n)$,
		then $\cH$ contains a matching $M$ with
		$n-\lvert V(M)\rvert\le\rho n$.
	\end{lemma}
	\begin{proof}
		Choose $\zeta,\varepsilon,n_0$ from Lemma~\ref{lem-robust-frac}
		and fix $\rho>0$.
	Let $N_0$ be supplied by Theorem~\ref{thm-morris} with
		$\theta=1/(4r^2)$, so $r^2\theta=1/4$. Note that $n^{1/4}\binom{n-2}{r-2}
		+\frac{n^{5/4}}2\binom{n-3}{r-3}
		\le\left(\frac1{(r-2)!}+\frac1{2(r-3)!}\right)n^{r-7/4}=o(n^{r-1})$. Thus we may choose an integer
		$n_1\ge\max\{n_0,N_0,r\}$ such that, for every $n\ge n_1$, $
		n^{1/4}\binom{n-2}{r-2}
		+\frac{n^{5/4}}2\binom{n-3}{r-3}
		\le\varepsilon n^{r-1}$ and $n^{-1/(4r^2)}\le\rho$, where
		the second inequality holds for sufficiently large $n$
		as $n^{-1/(4r^2)}\to0$. Thus $n_1$ can be chosen
		using only $r,\eta,\gamma,\rho$, while $\zeta$ is independent
		of $\rho$.
		
		Let $\cH$ satisfy the assumptions with $n\ge n_1$ and
		$J$ be any graph on $V(\cH)$ with $\Delta(J)\le n^{1/4}$.
		Define $
		\cH_J=\left\{F\in\cH\mid\binom F2\cap E(J)\ne\emptyset\right\}$.
		Fix $v\in V(\cH)$. Every edge of $\cH_J$ containing $v$
		contains some edge $f\in E(J)$. If $v\in f$, there are at
		most $\binom{n-2}{r-2}$ $r$-sets containing $f$.
		If $v\notin f$, there are at most $\binom{n-3}{r-3}$
		$r$-sets containing $f\cup v$. Since
		$\lvert E(J)\rvert\le n\Delta(J)/2\le n^{5/4}/2$, we obtain
		\[
		\begin{aligned}
			d_{\cH_J}(v)
			&\le\sum_{\substack{f\in E(J)\\v\in f}}\binom{n-2}{r-2}
			+\sum_{\substack{f\in E(J)\\v\notin f}}\binom{n-3}{r-3}\le d_J(v)\binom{n-2}{r-2}
			+\lvert E(J)\rvert\binom{n-3}{r-3}\\
			&\le n^{1/4}\binom{n-2}{r-2}
			+\frac{n^{5/4}}2\binom{n-3}{r-3}
			\le\varepsilon n^{r-1}.
		\end{aligned}
		\]
		It follows that $\Delta(\cH_J)\le\varepsilon n^{r-1}$.
		Lemma~\ref{lem-robust-frac} with $\cF=\cH_J$
		therefore gives a perfect fractional matching in $\cH-\cH_J$.
		This holds for every graph $J$ with $\Delta(J)\le n^{1/4}$. Thus
		Theorem~\ref{thm-morris}
		with $N=n$ and $\theta=1/(4r^2)$ yields a matching $M\subseteq\cH$ satisfying $
		n-\lvert V(M)\rvert\le n^{1-1/(4r^2)}
		=n\cdot n^{-1/(4r^2)}\le\rho n$.
	\end{proof}

	\begin{lemma}\label{lem-absorb}
		Let $r\geq3$ and $m\ge1$ be fixed integers with $r\mid m$, and let $a>0$. Let $\cH$ and $\cF$ be two $r$-graphs on vertex set $[n]$.
		If for every edge $S\in\cF$, there are
		at least $a n^m$ sets $Q\subseteq[n]\setminus S$ with $\lvert Q\rvert=m$ such that both
		$\cH[Q]$ and $\cH[Q\cup S]$ have perfect matchings,
		then for every sufficiently small constant $\lambda>0$ and all sufficiently
		large $n$, there is a set $A\subseteq V(\cH)$ such that
		\begin{enumerate}
			\item $\cH[A]$ has a perfect matching with $\lvert A\rvert\le2m\lambda n$;
			\item the hypergraph
			$\cH[A\cup S_1\cup\cdots\cup S_t]$ has a perfect matching for pairwise disjoint edges $S_1,\ldots,S_t\in\cF$ with $t\le a\lambda n/8$ and $S_i\cap A=\emptyset$ for each $i\in[t]$.
		\end{enumerate}
	\end{lemma}
	
	\begin{proof}
		Fix $0<\lambda<\min\{1/2,a/(32m^2)\}$ and let $n$ be
		sufficiently large in terms of $r,m,a,\lambda$.
		For each $S\in\cF$, let
		\[
		\cG(S)=\left\{Q\in\binom{[n]\setminus S}{m}
		\mid \cH[Q]\text{ and }\cH[Q\cup S]
		\text{ have perfect matchings}\right\}.
		\]
		Then $\lvert\cG(S)\rvert\ge an^m$.
		Select each $Q\in\binom{[n]}m$ independently with probability $
		p=\frac{\lambda n}{\binom nm}
		=(1+o(1))\lambda m!n^{1-m}$, where $0<p\leq1$ for sufficiently large $n$.
		Let $\mathcal X$ be the family of selected sets.
		Write $X=\lvert\mathcal X\rvert$ and
		$X_S=\lvert\mathcal X\cap\cG(S)\rvert$ for $S\in\cF$.
		These are sums of independent Bernoulli random variables, with $
		\mathbb E(X)=p\binom nm=\lambda n$ and $
		\mathbb E(X_S)=p\lvert\cG(S)\rvert
		\ge\frac{a\lambda n^{m+1}}{\binom nm}\ge a\lambda n$.
		Lemma~\ref{lem-chernoff}~(i), with parameter $1/2$, gives
		\[
		\begin{aligned}
			\mathbb P(X>2\lambda n)
			&\le\mathbb P\left(\lvert X-\mathbb E(X)\rvert
			\ge\mathbb E(X)/2\right)\le2e^{-\lambda n/12},\\
			\mathbb P(X_S<a\lambda n/2)
			&\le\mathbb P\left(\lvert X_S-\mathbb E(X_S)\rvert
			\ge\mathbb E(X_S)/2\right)\le2e^{-a\lambda n/12}.
		\end{aligned}
		\]
		Since $\lvert\cF\rvert\le \binom{n}{r}\leq n^r$,
		\[
		\begin{aligned}
			\mathbb P\left(\{X>2\lambda n\}\cup
			\bigcup_{S\in\cF}\{X_S<a\lambda n/2\}\right)
			&\le\mathbb P(X>2\lambda n)
			+\sum_{S\in\cF}\mathbb P(X_S<a\lambda n/2)\le2e^{-\lambda n/12}+2n^r e^{-a\lambda n/12}=o(1).
		\end{aligned}
		\]
		Thus, with probability $1-o(1)$, for $S\in\cF$,
		\begin{equation}\label{eq-abs-random}
			\lvert\mathcal X\rvert\le2\lambda n,\qquad
			\lvert\mathcal X\cap\cG(S)\rvert\ge a\lambda n/2
		\end{equation}
		
		Let $
		y=\left\lvert\left\{\{Q,Q'\}\in\binom{\mathcal X}{2}
		\mid Q\cap Q'\ne\emptyset\right\}\right\rvert$.
		For each $Q\in\binom{[n]}m$, $
		\left\lvert\left\{Q'\in\binom{[n]}m
		\mid Q'\ne Q,\ Q'\cap Q\ne\emptyset\right\}\right\rvert
		\le\sum_{v\in Q}\binom{n-1}{m-1}
		=\frac{m^2}{n}\binom nm$.
		Since distinct sets $Q,Q'$ are selected independently, $ \mathbb P(\{Q,Q'\}\subseteq\mathcal X)=p^2$. Hence
		\[
		\begin{aligned}
			\mathbb E(y)
			&=\frac12\sum_{Q\in\binom{[n]}m}
			\sum_{\substack{Q'\in\binom{[n]}m\\ Q'\cap Q\ne\emptyset}}
			\mathbb P(\{Q,Q'\}\subseteq\mathcal X)\le\frac12\binom nm\frac{m^2}{n}\binom nm\,p^2
			=\frac{m^2\lambda^2n}{2}.
		\end{aligned}
		\]
		Markov's inequality gives
		$\mathbb P(y>a\lambda n/8)\le4m^2\lambda/a<1/8$.
		Consequently, for sufficiently large $n$,
		\[
		\begin{aligned}
			\mathbb P\left(\{X>2\lambda n\}\cup
			\bigcup_{S\in\cF}\{X_S<a\lambda n/2\}
			\cup\{y>a\lambda n/8\}\right)\le2e^{-\lambda n/12}+2n^r e^{-a\lambda n/12}
			+\frac{4m^2\lambda}{a}<1.
		\end{aligned}
		\]
		Hence there exists a family $\mathcal X$ satisfying~\eqref{eq-abs-random} with $y\le a\lambda n/8$. Fix such a family $\mathcal X$.

		Start with $
		\cQ=\{Q\in\mathcal X\mid\cH[Q]\text{ has a perfect matching}\}$.
		Initially,
		$\cQ\cap\cG(S)=\mathcal X\cap\cG(S)$ for every $S\in\cF$.
		If there exists $\{Q,Q'\}\in\binom{\cQ}{2}$ with
		$Q\cap Q'\ne\emptyset$, choose such a pair $(Q,Q')$ and replace
		$\cQ$ by $\cQ\setminus\{Q\}$.
		Repeat this operation until $
		Q\cap Q'=\emptyset$ for every $\{Q,Q'\}\in\binom{\cQ}{2}$.
		Each step removes one set and decreases the number of intersecting
		pairs by at least one.  Since $y=\left\lvert\left\{\{Q,Q'\}\in\binom{\mathcal X}{2}
		\mid Q\cap Q'\ne\emptyset\right\}\right\rvert$,
		the process ends after at most $y$ steps.
		Thus the resulting family $\cQ^*$ satisfies
			\begin{equation}\label{eq-abs-reserve}
				\lvert\cQ^*\cap\cG(S)\rvert
				\ge\lvert\mathcal X\cap\cG(S)\rvert-y
				\ge\frac{a\lambda n}{2}-\frac{a\lambda n}{8}
				=\frac{3a\lambda n}{8}
			\end{equation}
			for $S\in\cF$. Set $A=\bigcup_{Q\in\cQ^*}Q$. Then $
			\lvert A\rvert=m\lvert\cQ^*\rvert
			\le m\lvert\mathcal X\rvert\le2m\lambda n$.
			For each $Q\in\cQ^*$, choose a perfect matching
			$M_Q\subseteq\cH[Q]$. Since $
			V(M_Q)\cap V(M_{Q'})=Q\cap Q'=\emptyset$ for $\{Q,Q'\}\in\binom{\cQ^*}{2}$,
			$\bigcup_{Q\in\cQ^*}M_Q$ is a matching with $
			V\left(\bigcup_{Q\in\cQ^*}M_Q\right)
			=\bigcup_{Q\in\cQ^*}V(M_Q)=\bigcup_{Q\in\cQ^*}Q=A$ from $
			V(M_Q)=Q$.
			This proves~(i).
		
		For~(ii), let $S_1,\ldots,S_t$ satisfy the stated conditions.
		The case $t=0$ follows from~(i), so assume $t\ge1$.
		By~\eqref{eq-abs-reserve},
		$\lvert\cQ^*\cap\cG(S_1)\rvert\ge3a\lambda n/8>0$,
		so choose $Q_1\in\cQ^*\cap\cG(S_1)$.
		For $2\le i\le t$, suppose that distinct sets
		$Q_j\in\cQ^*\cap\cG(S_j)$ have been chosen for $j\in[i-1]$.
		Since $i\le t\le a\lambda n/8$, we have
		\[
		\left\lvert(\cQ^*\cap\cG(S_i))
		\setminus\{Q_1,\ldots,Q_{i-1}\}\right\rvert
		\ge\frac{3a\lambda n}{8}-(i-1)
		\ge\frac{a\lambda n}{4}+1>0.
		\]
		Therefore choose $
		Q_i\in(\cQ^*\cap\cG(S_i))\setminus\{Q_1,\ldots,Q_{i-1}\}$.
		Since the members of $\cQ^*$ are pairwise disjoint,
		$Q_i\cap Q_j=\emptyset$ for every $j\in[i-1]$.
		Repeat this step until $Q_t$ has been chosen.
		For every $i\in[t]$, we have
		$Q_i\subseteq A$, $S_i\cap A=\emptyset$ and $Q_i\in\cG(S_i)$.
		Together with $S_i\cap S_j=\emptyset$ for $i\ne j$, this gives $(Q_i\cup S_i)\cap(Q_j\cup S_j)=\emptyset$ for $1\le i<j\le t$ and $(Q_i\cup S_i)\cap Q=\emptyset$ with $i\in[t]$ and $ Q\in\cQ^*\setminus\{Q_1,\ldots,Q_t\}$.
		By the definition of $\cG(S_i)$, choose a perfect matching
		$M_i\subseteq\cH[Q_i\cup S_i]$ for each $i\in[t]$.
		Then $
		M':=\left(\bigcup_{i=1}^t M_i\right)
		\cup\left(\bigcup_{Q\in\cQ^*\setminus\{Q_1,\ldots,Q_t\}}M_Q\right)$
		is a matching in $\cH$. Moreover, $
		V(M')
		=\left(\bigcup_{i=1}^t(Q_i\cup S_i)\right)
		\cup\left(\bigcup_{Q\in\cQ^*\setminus\{Q_1,\ldots,Q_t\}}Q\right)=A\cup S_1\cup\cdots\cup S_t$.
		Thus $M'$ is a perfect matching in
		$\cH[A\cup S_1\cup\cdots\cup S_t]$, proving~(ii).
	\end{proof}
	
	We keep the notations $
	p_r=1-\left(1-\frac1r\right)^{r-1}>\frac{1}{2}$ for $r\ge3$, $
	\alpha_r=\frac{(r-1)^{r-1}}{r^{r-2}}=r(1-p_r)$ and $
	c_r=\lfloor\alpha_r\rfloor$ in the following result.
	\begin{lemma}\label{lem-nonext}
	Fix integers $r\ge3$ and $K\ge c_r$, and let $\eta>0$.
		For all sufficiently large $n$, depending only on $r,K,\eta$,
		the following holds for all integers $s\ge1$ and $c_r\le b\le K$ with $n=rs+b$.
		If an $r$-graph $\cH$ on $n$ vertices satisfies $
		\sigma_r(\cH)>rD_{s-1}(n)$ with $\lvert \cH[B]\rvert \ge\eta n^r$ for every set $B\subseteq V(\cH)$ of size
		$\lfloor(1-1/r)n\rfloor$,
		then
		$\nu(\cH)\ge s$.
	\end{lemma}
	\begin{proof}
		Write $V=V(\cH)$ and $D=D_{s-1}(n)=B_r(n)-B_r^{s}(n)$.  Let $n$ be sufficiently large. Since $n=rs+b$ and
		$c_r\le b\le K$, $B_r^{s}(n)=\binom{n-s}{r-1}
		=\frac{(1-1/r)^{r-1}}{(r-1)!}n^{r-1}+O_{r,K}(n^{r-2})$ and $B_r(n)=\frac{n^{r-1}}{(r-1)!}+O_r(n^{r-2})$.
		Hence $D=(p_r+O_{r,K}(1/n))B_r(n)$.
		Let $
		L=\{v\in V\mid d_\cH(v)\le D\}$ and $ T=V\setminus L$.
		Then $\binom Lr\subseteq\cH$ by Lemma~\ref{lem-basic}~(i).
		Write $q=c_r+1$ and choose $
		0<\epsilon<\min\left\{p_r,1-\frac{\alpha_r}{q}\right\}$.
		Then $q(1-\epsilon)>\alpha_r$ and $\epsilon B_r(n)<D$ for
		sufficiently large $n$. Choose constants $0<\zeta_0<p_r/2$
		and $N_0$ for which Theorem~\ref{thm-han-nonext} holds with
		$\eta/2$ in place of $\eta$ and let $
		\beta=\min\left\{
		\frac1{4r},\frac{\epsilon}{8r(r-1)},
		\frac{\eta(r-1)!}{8r},\frac{\zeta_0}{8r(r-1)}
		\right\}$.
		We distinguish two cases according to the value of $\lvert L\rvert$.
		
		\medskip\noindent
		\textbf{Case~1.} $\lvert L\rvert<\beta n$.
		
		Let $Z=\{v\in V\mid d_\cH(v)\le\epsilon B_r(n)\}$. Then $Z\subseteq L$. For every $S\subseteq V$ with $\lvert S\rvert=r$ such that $\lvert S\cap Z\rvert\ge q$, if $S\notin\cH$,
		then $
		rD<\sigma_r(\cH)
		\le\sum_{v\in S}d_\cH(v)
		\le\bigl(q\epsilon+r-q\bigr)B_r(n)<rD$,
		a contradiction. Thus $
		\left\{S\in\binom Vr\mid\lvert S\cap Z\rvert\ge q\right\}
		\subseteq\cH$.
		Assume that $\bigcup_{i=0}^t Z_i$ is a partition of $Z$ with $
		\lvert Z_i\rvert=q$ for $i\in[t]$ and $
		\lvert Z_0\rvert<q$.
		Since $
		(r-q)t\le r\lvert Z\rvert<r\beta n
		\le\frac n4<\lvert T\rvert$,
		we may choose $R_1,\ldots,R_t\in\binom T{r-q}$ with
		$R_i\cap R_j=\emptyset$ for $i\ne j$.  Then
		$M_Z:=\{Z_i\cup R_i\mid i\in[t]\}\subseteq\cH$
		is a matching, since its edges are pairwise disjoint. It satisfies $
		V(M_Z)\cap L=Z\setminus Z_0$ and $
		V(M_Z)\cap Z_0=\emptyset$.
		Here $M_Z=\emptyset$ if $t=0$.
		
		Start with $M=M_Z$. If
		$L\setminus(Z\cup V(M))\ne\emptyset$, choose
		$v\in L\setminus(Z\cup V(M))$ and let $X=V(M)\cup Z_0$.
		Then $v\notin X$ and $d_\cH(v)>\epsilon B_r(n)$.
		Every edge
		of $M$ intersects $L\setminus Z_0$, so $
		\lvert M\rvert
		\le\sum_{F\in M}\lvert F\cap(L\setminus Z_0)\rvert
		=\lvert V(M)\cap(L\setminus Z_0)\rvert
		\le\lvert L\rvert$.
		Consequently, $
		\lvert X\rvert=r\lvert M\rvert+\lvert Z_0\rvert
		\le r\lvert L\rvert+c_r$.
		For sufficiently large $n$,
		\begin{align}
			\left\lvert\{F\in\cH\mid v\in F,\ F\cap X\ne\emptyset\}\right\rvert
			&\le\sum_{w\in X}
			\left\lvert\{F\in\cH\mid\{v,w\}\subseteq F\}\right\rvert\le(r\lvert L\rvert+c_r)\binom{n-2}{r-2}\nonumber\\
			&<(r\beta n+c_r)\frac{r-1}{n-1}B_r(n)\le4r(r-1)\beta B_r(n)
			\le\frac\epsilon2 B_r(n)<d_\cH(v).\label{e8}
		\end{align}
		Thus we choose an edge $F\in\cH$ with $
		v\in F$ and $F\cap X=F\cap(V(M)\cup Z_0)=\emptyset$.
		Then $M\cup\{F\}$ is a matching with $
		V(M\cup\{F\})\cap Z_0=\emptyset$ and $
		F'\cap(L\setminus Z_0)\ne\emptyset$ for every $F'\in M\cup\{F\}$.
		Consequently, $\lvert M\cup\{F\}\rvert\le\lvert L\setminus Z_0\rvert
		\le\lvert L\rvert$.
		Moreover, since $v\in L\setminus(Z\cup V(M))$ and $v\in F$, $
		\left\lvert L\setminus\bigl(Z\cup V(M\cup\{F\})\bigr)\right\rvert
		\le\lvert L\setminus(Z\cup V(M))\rvert-1$.  Let $M_1=M\cup\{F\}$.
For $i\ge1$, if $L\setminus(Z\cup V(M_i))\ne\emptyset$, choose
$u\in L\setminus(Z\cup V(M_i))$ and let $X=V(M_i)\cup Z_0$.
Since
$V(M_i)\cap Z_0=\emptyset$ with $\lvert M_i\rvert\le\lvert L\rvert$, \eqref{e8} holds with $u,M_i$
in place of $v,M$, respectively. Hence we may choose $F'\in\cH$ with
$u\in F'$ and $F'\cap X=\emptyset$, and set
$M_{i+1}=M_i\cup\{F'\}$.
The preceding argument shows that $M_{i+1}$ is a matching with
$\lvert M_{i+1}\rvert\le\lvert L\rvert$,
$V(M_{i+1})\cap Z_0=\emptyset$ and $
 \lvert L\setminus(Z\cup V(M_{i+1}))\rvert
 \le\lvert L\setminus(Z\cup V(M_i))\rvert-1$.
Continue the above process until $L\setminus(Z\cup V(M_i))=\emptyset$,
and denote the last matching by $M^*$. (If no edge needs to be added to $M_Z$, take $M^*=M_Z$.)
Note that at most $\lvert L\setminus Z\rvert$ edges are added in total. We have $
 \lvert M^*\rvert
 \le\lvert M_Z\rvert+\lvert L\setminus Z\rvert
 =\frac{\lvert Z\setminus Z_0\rvert}{q}+\lvert L\setminus Z\rvert
 \le\lvert L\setminus Z_0\rvert\le\lvert L\rvert$.
Since $Z\setminus Z_0\subseteq V(M_Z)\subseteq V(M^*)$ and
$L\setminus Z\subseteq V(M^*)$, we also have
$L\setminus Z_0=(Z\setminus Z_0)\cup(L\setminus Z)\subseteq V(M^*)$
and $V(M^*)\cap Z_0=\emptyset$.
		
		Let $\cH'=\cH-(V(M^*)\cup Z_0)$ and choose $W\subseteq V(\cH')$ with
		$\lvert W\rvert<r$ and $r\mid N$, where $N=\lvert V(\cH'-W)\rvert$.
		For sufficiently large $n$, $
		n-N+1=r\lvert M^*\rvert+\lvert Z_0\rvert+\lvert W\rvert+1
		\le r\beta n+c_r+r\le2r\beta n$, $
		N\ge n/2\ge N_0$ and $
		(n-N+1)B_r(n)
		\le\frac{2r\beta}{(r-1)!}n^r\le\frac\eta4n^r$.
		Lemma~\ref{lem-basic}~(iii), with $\cF=\emptyset$, implies that
		$\cH'-W$ is $(\eta/2)$-nonextremal. Moreover, $
		(n-N)\binom{n-2}{r-2}
		=\frac{(n-N)(r-1)}{n-1}B_r(n)
		\le4r(r-1)\beta B_r(n)\le\frac{\zeta_0}{2}B_r(n)$.
		Since $V(\cH'-W)\subseteq T$ and
		$D\ge(p_r-\zeta_0/2)B_r(n)$ for sufficiently large $n$,
		$
		\delta(\cH'-W)\ge D-(n-N)\binom{n-2}{r-2}
		\ge(p_r-\zeta_0)B_r(n)\ge(p_r-\zeta_0)B_r(N)$ by Lemma~\ref{lem-basic}~(ii).
		Theorem~\ref{thm-han-nonext} yields a perfect matching
		$M'\subseteq \cH'-W$. Hence $M^*\cup M'\subseteq\cH$ is a matching and $
		\lvert M^*\rvert+\lvert M'\rvert
		=\frac{n-\lvert Z_0\rvert-\lvert W\rvert}{r}
		\ge\frac{rs+b-c_r-r+1}{r}>s-1$.
		Thus $\nu(\cH)\ge s$.
		
		\medskip\noindent
		\textbf{Case~2.} $\lvert L\rvert\ge\beta n$.
		
		Since $\cH[L]$ is complete, $
		\binom{\lvert L\rvert-1}{r-1}\le d_\cH(v)\le D$ for every $v\in L$.
		Using
		$\binom{\lvert L\rvert-1}{r-1}
		=\lvert L\rvert^{r-1}/(r-1)!+O_r(n^{r-2})$ and $
		D=(p_r+O_{r,K}(1/n))B_r(n)$,
		we obtain $
		\left(\frac{\lvert L\rvert}{n}\right)^{r-1}
		\le p_r+O_{r,K}(1/n)$. Then $
		\frac{\lvert L\rvert}{n}\le p_r^{1/(r-1)}+o(1)$.
		For sufficiently large $n$, we have $
		\lvert T\rvert=n-\lvert L\rvert\geq n(1-p_r^{1/(r-1)}-o(1))\ge\frac{1-p_r^{1/(r-1)}}2\,n$, $ d_\cH(v)\ge\binom{\lvert L\rvert-1}{r-1}
		\ge\frac{\beta^{r-1}}{2(r-1)!}n^{r-1}$ for every $v\in L$ and $ d_\cH(v)>D\ge\frac{p_r}{2(r-1)!}n^{r-1}$ for every $v\in T$.
		Thus $\delta(\cH)\ge\gamma n^{r-1}$ where
		$\gamma=\min\{\beta^{r-1},p_r\}/(2(r-1)!)>0$.
		For distinct $x,y\in V$, define $
		\mathcal L(x,y)=
		\left\{P\in\binom{V\setminus\{x,y\}}{r-1}
		\mid P\cup x, P\cup y\in\cH\right\}$.
		Note that $\cH[L]$ is complete. If $x,y\in L$, then $
		\binom{L\setminus\{x,y\}}{r-1}\subseteq\mathcal L(x,y)$ with $
		\lvert\mathcal L(x,y)\rvert
		\ge\binom{\lvert L\rvert-2}{r-1}
		\ge\frac{\beta^{r-1}}{2(r-1)!}n^{r-1}$
		for sufficiently large $n$.
		If $d_\cH(x)+d_\cH(y)>2D$, then $
		\lvert\mathcal L(x,y)\rvert
		\ge d_\cH(x)+d_\cH(y)-\binom n{r-1}
		>2D-\binom n{r-1}
		\ge\frac{2p_r-1}{2(r-1)!}n^{r-1}$
		for sufficiently large $n$. Hence, for
		$\mu=\min\{\beta^{r-1},2p_r-1\}/(2(r-1)!)>0$,
		\begin{equation}\label{eq-common-link}
			\lvert\mathcal L(x,y)\rvert\ge\mu n^{r-1}
			\quad\text{if }\{x,y\}\subseteq L
			\text{ or }d_\cH(x)+d_\cH(y)>2D.
		\end{equation}
		In particular, this holds for $\{x,y\}\subseteq T$.
		Choose $0<\zeta<p_r/2$ as in Lemma~\ref{lem-almost} with
		parameters $\eta/2$ and $\gamma/2$ and $
		0<\theta<
		\min\left\{\frac\beta4,\frac{\eta(r-1)!}{8},
		\frac{\gamma(r-2)!}{4},\frac{\zeta}{8(r-1)}\right\}$. We next prove the following claim.
        \begin{claim}\label{clm-delete}
			For sufficiently large $n$, let $X\subseteq V(\cH)$ with
			$\lvert X\rvert\le\theta n$ and
			$\cH'=\cH-X$ with $N=\lvert V(\cH')\rvert$.
			Then
			\begin{enumerate}
				\item $\cH'$ is $(\eta/2)$-nonextremal with $\delta(\cH')\ge(\gamma/2)N^{r-1}$.
				\item $\sigma_r(\cH')\ge r(p_r-\zeta)B_r(N)$ if $\cH'$ has a nonedge.
			\end{enumerate}
		\end{claim}
		\noindent\textbf{Proof of Claim~\ref{clm-delete}.}
		Fix $X$ as in the statement. Since
		$\theta<\beta/4\le1/(16r)<1/2$, we have
		$N=n-\lvert X\rvert\ge(1-\theta)n\ge n/2$.
		Since
		$B_r(n)=(1/(r-1)!+o(1))n^{r-1}$ with
		$\theta<\eta(r-1)!/8$, we obtain $
		(\lvert X\rvert+1)B_r(n)\le(\theta n+1)B_r(n)=\left(\frac{\theta}{(r-1)!}+o(1)\right)n^r
		\le\frac\eta4 n^r<\frac\eta2 n^r$
		for sufficiently large $n$.
		From Lemma~\ref{lem-basic}~(iii) with $Z=X$ and
		$\cF=\emptyset$, $\cH'$ is
		$(\eta/2)$-nonextremal.
		For each $v\in V(\cH')$ and each $x\in X$, at most
		$\binom{n-2}{r-2}$ edges contain both $v$ and $x$. Hence
		\[
		\begin{aligned}
			d_{\cH'}(v)
			&\ge d_\cH(v)-\lvert X\rvert\binom{n-2}{r-2}\ge\gamma n^{r-1}-\theta n\binom{n-2}{r-2}\ge\left(\gamma-\frac{\theta}{(r-2)!}\right)n^{r-1}
			\ge\frac{3\gamma}{4}n^{r-1}
			\ge\frac\gamma2 N^{r-1}.
		\end{aligned}
		\]
	The last two inequalities follow from $\theta<\gamma(r-2)!/4$
		and $N\le n$, respectively. This proves (i).
For~(ii), suppose that $\cH'$ has a nonedge, and let $S$ be any
		such nonedge. Then $S$ is also a nonedge of $\cH$. Thus $
		\sum_{v\in S}d_{\cH'}(v)\ge\sum_{v\in S}d_\cH(v)
		-r\lvert X\rvert\binom{n-2}{r-2}\ge\sigma_r(\cH)-r\lvert X\rvert\binom{n-2}{r-2}$.
		It follows that $
		\sigma_r(\cH')\ge\sigma_r(\cH)
		-r\lvert X\rvert\binom{n-2}{r-2}
		>rD-r\lvert X\rvert\binom{n-2}{r-2}$.
		Since $D=(p_r+O_{r,K}(1/n))B_r(n)$ for sufficiently large $n$, we have $
		D\ge\left(p_r-\frac\zeta4\right)B_r(n)$.
		Also, by $\theta<\zeta/(8(r-1))$, $
		\lvert X\rvert\binom{n-2}{r-2}\le\theta n\frac{r-1}{n-1}B_r(n)\le2(r-1)\theta B_r(n)\le\frac\zeta4 B_r(n)$.
		Combining these estimates gives
		\begin{equation}\label{eq-rem-ore}
			\sigma_r(\cH')>rD-\frac{r\zeta}{4}B_r(n)
			\ge r\left(p_r-\frac\zeta2\right)B_r(n)
			\ge r(p_r-\zeta)B_r(N).
		\end{equation}
		The last inequality follows from $B_r(N)\le B_r(n)$ and
		$p_r-\zeta>0$. This proves Claim~\ref{clm-delete}.
		
		Let $
		\cF=\left\{S\in\binom Vr\mid\lvert S\cap T\rvert=1\right\}$ and $
		\kappa=\frac{(1-p_r^{1/(r-1)})\beta^{r-1}}{8(r-1)!}$.
		For sufficiently large $n$, we have $
		\lvert\cF\rvert=\lvert T\rvert\binom{\lvert L\rvert}{r-1}
		\ge\frac{1-p_r^{1/(r-1)}}2 n\,
		\frac{\beta^{r-1}}{2(r-1)!}n^{r-1}
		=2\kappa n^r$.
		Since $\cF=(\cF\cap\cH)\cup(\cF\setminus\cH)$, $
		\max\{\lvert\cF\cap\cH\rvert,\lvert\cF\setminus\cH\rvert\}
		\ge\kappa n^r$. We now distinguish two subcases according to the numbers of
		edges and nonedges in $\cF$.
		
		\medskip\noindent
		\textbf{Case~2.1.} $\lvert\cF\cap\cH\rvert\ge\kappa n^r$.

		Fix $S=\{x_1,\ldots,x_r\}\in\cF$ with
		$x_1\in T$.
		Since $\lvert\cF\cap\cH\rvert\ge\kappa n^r$
		and at most $r\binom{n-1}{r-1}$ $r$-sets
		intersect $S$, \begin{equation}\label{e10}
			\left\lvert
			\left\{
			F\in\cF\cap\cH\mid F\cap S=\emptyset
			\right\}
			\right\rvert\ge
			\kappa n^r-r\binom{n-1}{r-1}\ge\frac{\kappa}{2}n^r
		\end{equation}
		for sufficiently large $n$.
		Choose such an edge
		$F=\{y_1,\ldots,y_r\}$ with $y_1\in T$.
		Since $x_1,y_1\in T$ and
		$\lvert S\cup F\rvert=2r$, \eqref{eq-common-link} gives $
		\left\lvert\{P\in\mathcal L(x_1,y_1)
		\mid P\cap(S\cup F)=\emptyset\}\right\rvert
		\ge\mu n^{r-1}-2r\binom{n-1}{r-2}
		\ge\frac\mu2 n^{r-1}>0$
		for sufficiently large $n$. Thus we may choose
		$P_1\in\mathcal L(x_1,y_1)$ with $P_1\cap(S\cup F)=\emptyset$.
	For $2\leq i\leq r$, suppose that $P_1,\ldots,P_{i-1}$ have
		been chosen so that they are pairwise disjoint and
		$P_j\cap(S\cup F)=\emptyset$ for $j\in[i-1]$.
		Since $\lvert P_j\rvert=r-1$ for each $j\in[i-1]$, $
		\left\lvert S\cup F\cup\bigcup_{j\in[i-1]}P_j\right\rvert
		=2r+(i-1)(r-1)\le r^2+r$.
		By \eqref{eq-common-link}, since $x_i,y_i\in L$, $$
	\lvert\{P\in\mathcal L(x_i,y_i)
		\mid P\cap(S\cup F\cup\bigcup_{j\in[i-1]}P_j)=\emptyset
		\}\rvert\ge\mu n^{r-1}-(r^2+r)\binom{n-1}{r-2}
		\ge\frac\mu2 n^{r-1}>0.$$
		We may therefore choose $
		P_i\in\mathcal L(x_i,y_i)$ such that $
		P_i\cap\left(S\cup F\cup\bigcup_{j\in[i-1]}P_j\right)=\emptyset$.
		Repeating this step for $i=2,\ldots,r$ ends with the choice of
		$P_r$. The resulting sets $P_1,\ldots,P_r$ are pairwise disjoint
		with $P_i\cap(S\cup F)=\emptyset$ for $i\in[r]$.
		
		Let $
		Q=F\cup\bigcup_{i=1}^rP_i$.
		Then $Q\cap S=\emptyset$
		and $
		\lvert Q\rvert
		=\lvert F\rvert+\sum_{i=1}^r\lvert P_i\rvert
		=r+r(r-1)=r^2$.
		Since $P_i\in\mathcal L(x_i,y_i)$,
		$P_i\cup x_i,P_i\cup y_i\in\cH$ for $i\in[r]$. Then $\{P_i\cup y_i\mid i\in[r]\}$ and $\{F\}\cup\{P_i\cup x_i\mid i\in[r]\}$ are perfect matchings of $\cH[Q]$ and $\cH[Q\cup S]$, respectively.
		For each $S\in\cF$, there are at least $(\kappa/2)n^r$
		choices for $F$ by \eqref{e10} and
		at least $(\mu/2)n^{r-1}$ choices for
		each $P_i$ for each chosen $F$. Hence the number of
		tuples $(F,P_1,\ldots,P_r)$
		is at least $
		\frac{\kappa}{2}n^r
		\left(\frac{\mu}{2}n^{r-1}\right)^r=
		\frac{\kappa}{2}
		\left(\frac{\mu}{2}\right)^r n^{r^2}$.
		Since the sets $F,P_1,\ldots,P_r$
		are pairwise disjoint and their sizes
		are fixed, each $Q$ of size $r^2$ arises from at
		most $(r^2)!$ such tuples.
		Consequently, $$
		\left\lvert
		\left\{
		Q\in\binom{V\setminus S}{r^2}
		\mid
		\nu(\cH[Q])=r,\
		\nu(\cH[Q\cup S])=r+1
		\right\}
		\right\rvert\ge
		\frac{\kappa}{2(r^2)!}
		\left(\frac{\mu}{2}\right)^r n^{r^2}
		=an^{r^2},$$
		where $
		a=\frac{\kappa}{2(r^2)!}
		\left(\frac{\mu}{2}\right)^r>0$.
		This estimate holds uniformly for
		every $S\in\cF$.
		
		We may therefore apply
		Lemma~\ref{lem-absorb} to $\cH$ and
		$\cF$ with $m=r^2$ and this constant $a$.
		Choose $\lambda>0$ sufficiently small
	 such that $
		2r^2\lambda\le\frac{\theta}{4}$.
		There exists an $A\subseteq V$ satisfying both conclusions of Lemma~\ref{lem-absorb} with $
		\lvert A\rvert\le2r^2\lambda n\le\theta n/4$.
	In particular, $r\mid\lvert A\rvert$ because $\cH[A]$ has a perfect matching.
		Choose $
		\omega\in (0,
		\min\left\{\frac{\theta}{4},
		\frac{\beta}{4}\right\})$.
		Since $\lvert L\rvert\ge\beta n$,
		we have $
		\lvert L\setminus A\rvert\ge\lvert L\rvert-\lvert A\rvert\ge\left(\beta-\frac{\theta}{4}\right)n>\omega n+r$
		for sufficiently large $n$.
		Thus we may choose
		$R\subseteq L\setminus A$ with $
		\omega n\le\lvert R\rvert<\omega n+r$ and $
		\lvert R\rvert\equiv n\pmod r$.
		Let $\cH'=\cH-(A\cup R)$ with
		$N=\lvert V(\cH')\rvert$.
		Since $r\mid\lvert A\rvert$ and
		$\lvert R\rvert\equiv n\pmod r$,
		we obtain $
		N=n-\lvert A\rvert-\lvert R\rvert
		\equiv0\pmod r$.
		Moreover, $
		\lvert A\cup R\rvert=\lvert A\rvert+\lvert R\rvert<\frac{\theta n}{4}+\omega n+r<\theta n$
		for sufficiently large $n$.
		Choose $
	\rho\in (0,\min\left\{\frac{a\lambda}{16},\frac{\omega}{2r}\right\})$.
		Apply Claim~\ref{clm-delete} with $X=A\cup R$.
		Note that $r\mid N$ with $N\ge n/2$ from $\theta<\frac{1}{2}$ for sufficiently large $n$. From
		Lemma~\ref{lem-almost} with parameters $\eta/2$, $\gamma/2$
		and $\rho$, $\cH'$ has a matching covering all but at most $\rho N$ vertices if $\cH'$ has a nonedge.
		If $\cH'$ is complete, \eqref{eq-complete-match} gives a perfect matching since $r\mid N$.
		Thus in both cases there is a matching $M\subseteq\cH'$ with $
		\lvert Y\rvert\le\rho N\le\rho n$ where $
		Y=V(\cH')\setminus V(M)$.
		Write $
		Y\cap T=\{y_1,\ldots,y_\ell\}$.
		Since $\ell\le\lvert Y\rvert\le\rho n$, $
		(r-1)\ell\le(r-1)\rho n<(r-1)\frac{\omega}{2r}n<\omega n\le\lvert R\rvert$.
		Hence we can choose pairwise disjoint
		sets $
		R_1,\ldots,R_\ell\in\binom R{r-1}$.
		For each $i\in[\ell]$, let $
		S_i=y_i\cup R_i$.
		Since $y_i\in T$, $R_i\subseteq R\subseteq L$ with $\lvert R_i\rvert=r-1$, we have $S_i\in\cF$  for $i\in[\ell]$.
		Furthermore, since $Y\subseteq V(\cH')$ and
		$R\cap A=\emptyset$, the sets $S_i$ are
		pairwise disjoint with $S_i\cap A=\emptyset$ for $i\in[\ell]$.
		Also, $
		\ell\le\rho n
		<\frac{a\lambda n}{8}$.
		Therefore Lemma~\ref{lem-absorb}~(ii)
		gives a matching $M_A$ in $\cH$
		with $
		V(M_A)=A\cup\bigcup_{i=1}^{\ell}S_i$.
		If $\ell=0$, take instead a perfect
		matching of $\cH[A]$, as guaranteed
		by Lemma~\ref{lem-absorb}~(i).
		
		Note that $M\subseteq\cH'$ with $V(M)\cap(A\cup R)=\emptyset$ and $
		Y=V(\cH')\setminus V(M)$. Then $
		V(M)\cap V(M_A)=\emptyset$.
		Let $M^*=M\cup M_A$. Since $R\subseteq L$, $
		V\setminus V(M^*)
		=\left(R\setminus\bigcup_{i=1}^{\ell}R_i\right)\cup(Y\cap L)
		\subseteq L$.
		Thus $\cH[V\setminus V(M^*)]$ is complete. A maximum matching
		in this induced subhypergraph can be added to $M^*$, so
		\eqref{eq-complete-match} gives $
		\nu(\cH)\ge\lvert M^*\rvert+
		\left\lfloor\frac{n-r\lvert M^*\rvert}{r}\right\rfloor
		=\left\lfloor\frac nr\right\rfloor
		=s+\left\lfloor\frac br\right\rfloor\ge s$.
		This completes Case~2.1.
		
		\medskip\noindent
		\textbf{Case~2.2.}
		$\lvert\cF\setminus\cH\rvert\ge\kappa n^r$.
		
		Define $
		\cP=
		\left\{
		(x,y)\in L\times T
		\mid d_\cH(x)+d_\cH(y)>2D
		\right\}$.
		Fix $
		S=\{y,x_1,\ldots,x_{r-1}\}
		\in\cF\setminus\cH$
		where $y\in T$ and
		$x_1,\ldots,x_{r-1}\in L$.
		If $(x_i,y)\notin\cP$
		for every $i\in[r-1]$,
		then $
		d_\cH(x_i)+d_\cH(y)\le2D$. It follows that $
		d_\cH(x_i)\le2D-d_\cH(y)$.
		Consequently,
		\[
		\begin{aligned}
			\sigma_r(\cH)
			&\le\sum_{v\in S}d_\cH(v)=d_\cH(y)+
			\sum_{i=1}^{r-1}d_\cH(x_i)\le d_\cH(y)
			+(r-1)(2D-d_\cH(y))=rD-(r-2)(d_\cH(y)-D)<rD,
		\end{aligned}
		\]
		where the last inequality follows
		from  $r\ge3$ and $d_\cH(y)>D$ for $y\in T$.
		This contradicts $\sigma_r(\cH)\geq rD$.
		Therefore, every member of
		$\cF\setminus\cH$ contains a pair
		$(x,y)\in\cP$. Since such a fixed pair
		belongs to at most
		$\binom{n-2}{r-2}$ $r$-sets, $
		\kappa n^r\le\lvert\cF\setminus\cH\rvert\le\sum_{(x,y)\in\cP}
		\left\lvert
		\left\{
		S\in\cF\setminus\cH
		\mid\{x,y\}\subseteq S
		\right\}
		\right\rvert\le\lvert\cP\rvert
		\binom{n-2}{r-2}\le\lvert\cP\rvert n^{r-2}$.
		Hence $
		\lvert\cP\rvert\ge\kappa n^2$.
		
		For distinct vertices $x,y\in V$,
		define $
		\mathcal C(x,y)=
		\left\{
		C\in\binom{V\setminus\{x,y\}}{3r-1}
		\mid
		\nu(\cH[C\cup x])
		=\nu(\cH[C\cup y])=3
		\right\}$.
		\begin{claim}\label{clm-c-count}
			There exists a constant $a_1>0$, independent of $x,y$ and $n$,
			such that $
			\lvert\mathcal C(x,y)\rvert\ge a_1n^{3r-1}$
			for all distinct $x,y\in V$.
		\end{claim}
		\noindent\textbf{Proof of Claim~\ref{clm-c-count}.}
		The following cases depend on the location of $(x,y)$ in $V$.
		First suppose that $x\in L$ and
		$y\in T$.
		Choose $(u,v)\in\cP$ such that $
		\{u,v\}\cap\{x,y\}=\emptyset$. For any fixed distinct vertices $x,y\in V(\cH)$,
		at most $2n$ pairs $(u,v)\in\mathcal P$ satisfy
		$\{u,v\}\cap\{x,y\}\ne\emptyset$.
		Hence, since $\lvert\cP\rvert\ge\kappa n^2$, \begin{equation}\label{e11}
			\left\lvert\left\{(u,v)\in\mathcal P
			\mid \{u,v\}\cap\{x,y\}=\emptyset\right\}\right\rvert\ge\kappa n^2-2n
			\ge\frac{\kappa}{2}n^2
		\end{equation}
		for sufficiently large $n$.
		Note that $
		\mathcal L(x,y)=
		\left\{P\in\binom{V\setminus\{x,y\}}{r-1}
		\mid P\cup x, P\cup y\in\cH\right\}$.
		Since $x,u\in L$, $v,y\in T$ and $(u,v)\in\cP$,
		\eqref{eq-common-link} gives $
		\min\{\lvert\mathcal L(x,u)\rvert,
		\lvert\mathcal L(u,v)\rvert,
		\lvert\mathcal L(v,y)\rvert\}
		\ge\mu n^{r-1}$.
		For any $X\subseteq V$ with $\lvert X\rvert\le2r+2$,
		at most $\lvert X\rvert\binom{n-1}{r-2}$ members in each
		of these three families intersect $X$. Thus each family has at least \begin{equation}\label{e12}
			\mu n^{r-1}-\lvert X\rvert\binom{n-1}{r-2}\ge\mu n^{r-1}-(2r+2)\binom{n-1}{r-2}=\mu n^{r-1}-O_r(n^{r-2})\ge\frac\mu2 n^{r-1}>0
		\end{equation}
		members disjoint from $X$, for sufficiently large $n$.
		We may therefore first choose $P_1\in\mathcal L(x,u)$ disjoint
		from $\{x,y,u,v\}$. Since $
		\lvert\{x,y,u,v\}\cup P_1\rvert=4+(r-1)
		\le4+2(r-1)$,
		we next choose $P_2\in\mathcal L(u,v)$ disjoint from
		$\{x,y,u,v\}\cup P_1$. Finally, $
		\lvert\{x,y,u,v\}\cup P_1\cup P_2\rvert=4+2(r-1)$,
		so we choose $P_3\in\mathcal L(v,y)$ disjoint from this set.
		There are at least $\mu n^{r-1}/2$ choices at each step by \eqref{e12}. Each $P_i$ with $i\in[3]$ satisfies $
		P_i\cap\{x,y,u,v\}=\emptyset$ and $P_i\cap P_j=\emptyset$ with $1\le i<j\le3$.
		Let $
		C=\{u,v\}\cup P_1\cup P_2\cup P_3$. Then $
		\lvert C\rvert
		=2+3(r-1)=3r-1$ and $C\cap\{x,y\}=\emptyset$.
		Since $P_1\in\mathcal L(x,u)$,
		$P_2\in\mathcal L(u,v)$ and
		$P_3\in\mathcal L(v,y)$,
		the three edges $
		P_1\cup x$, $
		P_2\cup u$ and $
		P_3\cup v$
		are pairwise disjoint and form a perfect matching of
		$\cH[C\cup x]$.
		Similarly, $
		P_1\cup u$, $
		P_2\cup v$ and $
		P_3\cup y$
		are pairwise disjoint edges and form a perfect matching of
		$\cH[C\cup y]$.
		Thus $C\in\mathcal C(x,y)$.
		The number of tuples
		$(u,v,P_1,P_2,P_3)$ is at least $
		\frac{\kappa}{2}n^2
		\left(\frac{\mu}{2}n^{r-1}\right)^3
		=
		\frac{\kappa}{2}
		\left(\frac{\mu}{2}\right)^3n^{3r-1}$ by \eqref{e11} and \eqref{e12}.
		The same argument applies when
		$x\in T$ and $y\in L$, since
		$\mathcal C(x,y)=\mathcal C(y,x)$.
		
		Now suppose that $x,y\in L$ or
		$x,y\in T$.
		Choose $P\in\mathcal L(x,y)$.
		Since $\cH[L]$ is complete with
		$\lvert L\rvert\ge\beta n$, we can
		choose two disjoint edges
		$F_1,F_2\in\binom Lr$ such that $
		(P\cup F_1\cup F_2)
		\cap\{x,y\}=\emptyset$
		and $
		P\cap F_1=P\cap F_2
		=\emptyset$. Moreover, when choosing $F_1,F_2$
		successively, the number of
	choices for each $F_i$, $i\in[2]$, is at least \begin{equation}\label{e13}
			\binom{\lvert L\rvert-(2r+1)}{r}\ge \frac{1}{r!}
			\prod_{j=0}^{r-1}(\beta n-2r-1-j)\ge \frac{\beta^r}{2r!}n^r
		\end{equation}
		for sufficiently large $n$.
		Let $
		C=P\cup F_1\cup F_2$.
	Then $
		\lvert C\rvert=(r-1)+2r=3r-1$. The families $
		\{P\cup x,F_1,F_2\}$
		and $
		\{P\cup y,F_1,F_2\}$
		are perfect matchings in
		$\cH[C\cup x]$ and
		$\cH[C\cup y]$, respectively.
		Hence $C\in\mathcal C(x,y)$. By \eqref{eq-common-link} and \eqref{e13},
		the number of tuples
		$(P,F_1,F_2)$ is at least $
		\frac{\mu}{2}
		\left(\frac{\beta^r}{2r!}\right)^2
		n^{3r-1}$.
		
	For a fixed $(3r-1)$-vertex set $C$, at most $(3r-1)!$ tuples in both construction yield $C$.
		Consequently, choosing $
		a_1=
		\frac{1}{(3r-1)!}
		\min\left\{
		\frac{\kappa}{2}
		\left(\frac{\mu}{2}\right)^3,\
		\frac{\mu}{2}
		\left(\frac{\beta^r}{2r!}\right)^2
		\right\}>0$,
		we obtain $
		\lvert\mathcal C(x,y)\rvert
		\ge a_1n^{3r-1}$
		for every pair of distinct vertices
		$x,y\in V$.
		This proves Claim~\ref{clm-c-count}.
		
		We now consider every $S\in\binom Vr$.
		Fix $
		S=\{x_1,\ldots,x_r\}\in\binom Vr$.
		Since $\cH[L]$ is complete and
		$\lvert L\rvert\ge\beta n$,
		there are at least $
		\binom{\lvert L\setminus S\rvert}{r}
		\ge\frac{\beta^r}{2r!}n^r$
		choices for an edge $
		F=\{y_1,\ldots,y_r\}\in\binom Lr$ with $F\cap S=\emptyset$.
		Fix orderings of $S$ and $F$.
		We first choose $C_1\in\mathcal C(x_1,y_1)$ with $C_1\cap(S\cup F)=\emptyset$.
		Note that $S\cup F$ has size $2r$, and each of its vertices belongs to at most $\binom{n-1}{3r-2}$ sets
		of size $3r-1$. By Claim~\ref{clm-c-count},  $
		\left\lvert\{C\in\mathcal C(x_1,y_1)
		\mid C\cap(S\cup F)=\emptyset\}\right\rvert\ge a_1n^{3r-1}-2r\binom{n-1}{3r-2}
		\ge\frac{a_1}{2}n^{3r-1}>0$ for sufficiently large $n$.
		Thus there are at least $(a_1/2)n^{3r-1}$ choices for $C_1$.
		Assume that
		$C_j\in\mathcal C(x_j,y_j)$  has been chosen for each $j<i$ with $2\le i\le r$,
		with these sets pairwise disjoint and disjoint from $S\cup F$.
		Since each $C_j$ has size $3r-1$, $
		\left\lvert S\cup F\cup\bigcup_{j\in[i-1]}C_j\right\rvert
		=2r+(i-1)(3r-1)$.
		At most $[2r+(i-1)(3r-1)]\binom{n-1}{3r-2}$ members of
		$\mathcal C(x_i,y_i)$ intersect $S\cup F\cup\bigcup_{j\in[i-1]}C_j$. Hence Claim~\ref{clm-c-count}
		gives
		\[
		\begin{aligned}
			\lvert\{C\in\mathcal C(x_i,y_i)
			\mid C\cap(S\cup F\cup\bigcup_{j<i}C_j)=\emptyset
			\}\rvert\ge a_1n^{3r-1}
			-[2r+(i-1)(3r-1)]\binom{n-1}{3r-2}
			\ge\frac{a_1}{2}n^{3r-1}>0.
		\end{aligned}
		\]
		The penultimate inequality holds since $r$ and $a_1$ are fixed for all $i\in\{2,\ldots,r\}$ when $n$
		is sufficiently large.
		We may therefore choose $C_i\in\mathcal C(x_i,y_i)$ disjoint
		from $S\cup F\cup\bigcup_{j\in[i-1]}C_j$ satisfying $
		C_i\cap(S\cup F)=\emptyset$ and $C_i\cap C_j=\emptyset$ with $j<i$.
		Repeat this step until $C_r$ has been chosen. Then
		$C_1,\ldots,C_r$ are pairwise disjoint and $C_i\cap(S\cup F)=\emptyset$ for $i\in[r]$,
		with at least $(a_1/2)n^{3r-1}$ choices at every step.
		
		Let $
		Q=F\cup\bigcup_{i=1}^r C_i${\pu .}
		Then $
		\lvert Q\rvert
		=r+r(3r-1)=3r^2$ and $Q\cap S=\emptyset$ by the above construction.
		By the definition of
		$\mathcal C(x_i,y_i)$,
		there exist perfect matchings
		$M_i^x$ with $V(M_i^x)=C_i\cup x_i$ and $M_i^y$ with $V(M_i^y)=C_i\cup y_i$ in
		$\cH[C_i\cup x_i]$ and
		$\cH[C_i\cup y_i]$,
		respectively.
		It follows that $
		\bigcup_{i=1}^rM_i^y$ and $
		\{F\}\cup\bigcup_{i=1}^rM_i^x$ are perfect matchings of $\cH[Q]$ and $\cH[Q\cup S]$, respectively. There are at least
		$\beta^r n^r/(2r!)$ choices for
		$F$ and, for each such choice,
		at least $(a_1/2)n^{3r-1}$
		choices for each $C_i$.
		Hence the number of tuples
		$(F,C_1,\ldots,C_r)$ is at least $\frac{\beta^r}{2r!}n^r
		\left(\frac{a_1}{2}n^{3r-1}\right)^r=
		\frac{\beta^r}{2r!}
		\left(\frac{a_1}{2}\right)^r
		n^{3r^2}$. For a fixed $3r^2$-vertex set $Q$, at most $(3r^2)!$ tuples yield $Q$.
		Consequently, choosing $
		a_2=
		\frac{\beta^r}{2r!(3r^2)!}
		\left(\frac{a_1}{2}\right)^r>0$,
		we obtain $
		\lvert\{
		Q\in\binom{V\setminus S}{3r^2}
		\mid
		\nu(\cH[Q])=3r,\
		\nu(\cH[Q\cup S])=3r+1\}\rvert\ge a_2n^{3r^2}$
		for every $S\in\binom Vr$.
		
		We may now apply
		Lemma~\ref{lem-absorb}
		with $\cF=\binom Vr$,
		$m=3r^2$ and $a=a_2$.
		Choose $\lambda>0$ sufficiently
		small with $
		6r^2\lambda\le\frac{\theta}{4}$.
		Then there exists an $A\subseteq V$ satisfying  both conclusions of Lemma~\ref{lem-absorb} with $
		\lvert A\rvert\le6r^2\lambda n\le\theta n/4$ and $
		r\mid\lvert A\rvert$.
		Choose $W\subseteq V\setminus A$
		such that $
		0\le\lvert W\rvert<r$ with $
		\lvert W\rvert\equiv n\pmod r$.
		Let $
		\cH'=\cH-(A\cup W)$ with $
		N=\lvert V(\cH')\rvert$.
		Since $r\mid\lvert A\rvert$,
		we have $
		N=n-\lvert A\rvert-\lvert W\rvert
		\equiv0\pmod r$.
		Moreover, $
		\lvert A\cup W\rvert
		\le\frac{\theta n}{4}+r-1
		\le\theta n$
		for sufficiently large $n$.
		Choose $0<\rho<a_2\lambda/16$ and apply Claim~\ref{clm-delete}
		with $X=A\cup W$. Note that $r\mid N$ and $N\ge n/2$.
		From Lemma~\ref{lem-almost} with parameters $\eta/2$, $\gamma/2$
		and $\rho$, for sufficiently large $n$, $\cH'$ has a matching covering all but at most $\rho N$ vertices if $\cH'$ has
		a nonedge. If $\cH'$ is complete, \eqref{eq-complete-match} gives a perfect matching since $r\mid N$. Consequently, in both cases there is a
		matching $M\subseteq\cH'$ satisfying  $
		\lvert Y\rvert\le\rho N\le\rho n$ where $
		Y=V(\cH')\setminus V(M)$.
		Since $r\mid N$ and $M$ is a
		matching in an $r$-graph, $
		\lvert Y\rvert=N-r\lvert M\rvert
		\equiv0\pmod r$.
		Thus $Y$ can be partitioned
		into pairwise disjoint $r$-sets $
		Y=S_1\cup\cdots\cup S_\ell$ with $
		\ell=\frac{\lvert Y\rvert}{r}\le\frac{\rho n}{r}
		<\frac{a_2\lambda n}{8}$.
		Moreover,
		each $S_i$ is disjoint from $A$
		because $Y\subseteq V(\cH')$.
		By Lemma~\ref{lem-absorb}~(ii),
		there exists a matching $M_A$
		in $\cH$ with $
		V(M_A)=A\cup Y$.
		If $\ell=0$, take a perfect
		matching of $\cH[A]$ instead.
		Since $M\subseteq\cH'$ and $
		Y=V(\cH')\setminus V(M)$, $V(M)\cap V(M_A)=\emptyset$.
		Moreover,
		$V(M\cup M_A)=V\setminus W$. Therefore, $
		\nu(\cH)\ge\lvert M\rvert+\lvert M_A\rvert=\frac{n-\lvert W\rvert}{r}=\left\lfloor\frac nr\right\rfloor\ge s$.
		This completes Case~2.2 and the proof of the lemma.
	\end{proof}
	
	\begin{lemma}\label{lem-extremal}
		Let $r\ge3$ and $C_0\ge0$ be integers.
		There exist $\eta>0$ and $n_0$ such that the following holds for all
		integers $s\ge1$ and $c$ with $0\le c\le C_0$.
		Let $\cH$ be an $r$-graph on $n=rs+c\ge n_0$ vertices with $
		\sigma_r(\cH)>rD_{s-1}(n)$.
		If there is a set $B\subseteq V(\cH)$ with $
		\lvert B\rvert =n-s+1$ and $ \lvert \cH[B]\rvert \le\eta n^r$,
		then $\nu(\cH)\ge s$.
	\end{lemma}
	\begin{proof}
		Fix a constant $
		0<\rho_0<\bigl(8^{r-1}r^{5(r-1)}r!\bigr)^{-1}$.
		Choose positive constants $\rho=\rho_0/2^r$, $ \varepsilon$, $ \beta$, $\xi$ and $\eta$ with $
		\varepsilon\le
		\min\{\frac1{12r^3},\frac{\rho(r-2)!}{8r}\}$,
		$\beta\le
		\min\{\frac\rho4,
		\frac{\varepsilon}{4\cdot3^{r-1}(r-2)!}\}$, $
		\xi\le\frac{\beta\varepsilon}{128r}$ and $
		\eta\le
		\min\{\frac{\xi}{4r},
		\frac{\xi^r}{4^r r!},
		\frac{\rho}{48r(r-1)}\}$.
		These constants depend only on $r$. Choose $n_0$ sufficiently
		large in terms of $r$ and $C_0$ such that $
		\varepsilon n\ge12$ and $n\ge\max\{4C_0,2r/\xi,16r^6\}$ for all $n\ge n_0$.
		Let $\cH$ and $B$ satisfy the assumptions. If $\cH$ is complete, \eqref{eq-complete-match} gives $\nu(\cH)=\lfloor n/r\rfloor\ge s$. So we assume that $\cH$
		has a nonedge in the following. Let $
		A=V(\cH)\setminus B$, $D=D_{s-1}(n)$, $
		t=\lfloor\varepsilon n\rfloor$ and $a=s-t$.
		Then $\lvert A\rvert=s-1$, $t\ge\varepsilon n/2\ge2$ and $
		a\ge\frac{n-C_0}{r}-\varepsilon n\ge\frac n{2r}>0$.
		
		By Definition~\ref{def-good}, $
		\mathcal T(A,B)
		=
		\left\{
		F\in\binom{V(\cH)}r\mid
		F\cap A\ne\emptyset,\ F\cap B\ne\emptyset
		\right\}$.
		We first consider the size of $\cT(A,B)\setminus\cH$.
		Let $
		L=\{v\in V(\cH)\mid d_\cH(v)\le D\}$.
		By Lemma~\ref{lem-basic}~(i), $\binom Lr\subseteq\cH$.
		For $v\in B$, define $
		m_B(v)=\left\lvert
		\{F\in\cT(A,B)\setminus\cH\mid v\in F\}
		\right\rvert$.
		The number of $r$-sets containing $v$ and intersecting $A$ is
		$B_r(n)-\binom{\lvert B\rvert-1}{r-1}=D$. Every edge of $\cH$
		containing $v$ either intersects $A$ or belongs to $\cH[B]$, so
		\begin{equation}\label{eq-ext-deg}
			d_\cH(v)=D-m_B(v)+d_{\cH[B]}(v).
		\end{equation}
		Then $m_B(v)<d_{\cH[B]}(v)$ for $v\in B\setminus L$ and $m_B(v)\le D$ for $v\in B\cap L$.
		Since $\cH[B\cap L]$ is complete, $
		\binom{\lvert B\cap L\rvert}{r}\le\lvert\cH[B]\rvert\le\eta n^r$.
		If $\lvert B\cap L\rvert\ge r$, the left-hand side is at least
		$(\lvert B\cap L\rvert-r+1)^r/r!$. If
		$\lvert B\cap L\rvert<r$, the following bound is obvious.
		Thus, in both cases, $
		\lvert B\cap L\rvert\le(r!\eta)^{1/r}n+r$.
		Using
		$D\le B_r(n)\le n^{r-1}$, we obtain
		\[
		\begin{aligned}
			\lvert\cT(A,B)\setminus\cH\rvert
			&\le\sum_{v\in B}m_B(v)\le\sum_{v\in B\setminus L}d_{\cH[B]}(v)
			+\lvert B\cap L\rvert D\le r\lvert\cH[B]\rvert+\bigl((r!\eta)^{1/r}n+r\bigr)n^{r-1}\\
			&\le r\eta n^r+\bigl((r!\eta)^{1/r}n+r\bigr)n^{r-1}\le\left(\frac\xi4+\frac\xi4+\frac rn\right)n^r
			\le\xi n^r.
		\end{aligned}
		\]
		Recall that a vertex $v\in V(\cH)$
		is $\beta$-good with respect to $\cT(A,B)$ if $\lvert\{F\in\cT(A,B)\setminus\cH\mid v\in F\}\rvert\leq\beta n^{r-1}$. Let $B_A=\{v\in A\mid v\;\text{is not }	\beta\text{-good  with respect to } \cT(A,B) \}$ with $b_A=\lvert B_A\rvert$. Since each vertex of $B_A$ belongs to at least $\beta n^{r-1}$ edges of $\cT(A,B)\setminus\cH$, $
		b_A\beta n^{r-1}
		\le r\lvert\cT(A,B)\setminus\cH\rvert
		\le r\xi n^r$.
		Therefore $
		b_A\le\frac{r\xi}{\beta}n
		\le\frac{\varepsilon n}{128}\le\frac t{64}\le t-1$. It follows that at least $(s-1)-(t-1)=a$ vertices of $A$ are $\beta$-good. Choose $G\subseteq\{v\in A\mid v\;\text{is }	\beta\text{-good with respect to } \cT(A,B)\}$ with $\lvert G\rvert=a$, and let $R=A\setminus G$.
		Then $\lvert R\rvert=t-1$.
		
		We next show that $\cF=\cH-G$ with $N=\lvert V(\cF)\rvert$
		contains a matching of size $t$. For each $v\in V(\cF)$, let $
		m_G(v)=\left\lvert
		\left\{F\in\binom{V(\cH)}r\setminus\cH
		\mid v\in F,\ F\cap G\ne\emptyset\right\}
		\right\rvert$.
		The number of $r$-sets containing
		$v$ and intersecting $G$ is $
		D_a(n)=\binom{n-1}{r-1}
		-\binom{n-a-1}{r-1}$, of which exactly $m_G(v)$ are not edges of $\cH$. Thus
		\begin{equation}\label{eq-ext-inherit}
			d_\cF(v)=d_\cH(v)-D_a(n)+m_G(v).
		\end{equation}
		Note that $a+t=s$ and every nonedge of $\cF$ is also a nonedge of $\cH$.
		Hence, if $\cF$ has a nonedge $S$, then $
		\sum_{v\in S}d_\cF(v)
		\ge\sum_{v\in S}d_\cH(v)-rD_a(n)>rD_{s-1}(n)-rD_a(n)
		=rD_{t-1}(N)$.
		Moreover, $
		N=n-s+t\ge(1-1/r)n$, and hence $
		3r^2(t-1)\le3r^2\varepsilon n\le\frac n{4r}\le N$.
		If $\cF$ has a nonedge, Theorem~\ref{thm-bpr}~(i) gives
		$\nu(\cF)\ge t$. If not, \eqref{eq-complete-match} gives
		$\nu(\cF)=\lfloor N/r\rfloor\ge t$, since $N-rt=(r-1)a+c\ge0$.
		
		Let $
		X=\{v\in V(\cF)\mid m_G(v)>\rho n^{r-1}\}$ and $b_B=\lvert X\cap B\rvert$.
		For $v\in X\cap B$, $m_G(v)>\rho n^{r-1}$ and every set counted by $m_G(v)$ belongs to $\cT(A,B)\setminus\cH$. Consequently, $
		b_B\rho n^{r-1}
		\le\sum_{v\in B}m_G(v)
		\le r\lvert\cT(A,B)\setminus\cH\rvert
		\le r\xi n^r$.
		Since $X\subseteq V(\cF)=R\cup B$ and $\beta\le\rho$, it follows that
		\begin{equation}\label{eq-ext-size}
			b_B\le\frac{r\xi}{\rho}n\le\frac{\varepsilon n}{128},
			\qquad
			\lvert X\rvert\le t-1+b_B.
		\end{equation}
		Let $M$ be a matching of $\cF$ of size $t$ such that $M$ maximizes
		$\lvert V(M)\cap X\rvert$. We prove the following claim.
		\begin{claim}\label{clm-fex}
			There is no edge
			$F\in\cF$ with $
			F\cap V(M)=\emptyset$ and $F\cap X\ne\emptyset$.
		\end{claim}
		\noindent\textbf{Proof of Claim~\ref{clm-fex}.}
		Otherwise, there is an edge $F\in\cF$ satisfying $
		F\cap V(M)=\emptyset$ and $F\cap X\ne\emptyset$. If there exists an edge $E\in M$ with $E\cap X=\emptyset$, then $M'=(M\setminus\{E\})\cup\{F\}$ is a matching of size $t$ in $\cF$ and $\lvert V(M')\cap X\rvert>\lvert V(M)\cap X\rvert$, contrary to the choice of $M$. Thus every edge of
		$M$ intersects $X$ and \eqref{eq-ext-size} gives \begin{equation}\label{e17}
			0\le\sum_{E\in M}(\lvert E\cap X\rvert-1)
			=\lvert V(M)\cap X\rvert-t\le b_B-1{\pu .}
		\end{equation}
		In particular, $b_B\ge1$. Let $M^*=\{E\in M\mid \lvert E\cap X\rvert=1\}$. Then $\lvert M^*\rvert\geq t-b_B+1$ by \eqref{e17}. Let $
		Q=\{v\in R\mid E\cap X=v\text{ for some }E\in M^*,\ v\text{ is }\beta\text{-good}\}$.
	Since each edge of $M^*$ contains exactly one vertex of $X$,
		$\left\lvert\left\{E\in M^*\mid E\cap X\subseteq B_A\cup B\right\}\right\rvert\le b_A+b_B$.
		Hence $
		\lvert Q\rvert\ge t-b_A-2b_B+1
		\ge t-\frac{3\varepsilon n}{128}+1\ge t/2$.
		
		Let $B_0=B\setminus(V(M)\cup X\cup F)$. Since
		$\lvert B\rvert\ge(1-1/r)n$, $
		\lvert B_0\rvert
		\ge\lvert B\rvert-rt-b_B-r
		\ge\left(1-\frac1r-r\varepsilon-\frac{\varepsilon}{128}\right)n-r
		\ge n/3$
		for sufficiently large $n$. For a fixed $u\in Q$, let $\cA=\{\{u,v\}\cup P\mid v\in Q\setminus u,\;P\in\binom{B_0}{r-2}\}$. Since $u,v\in R\subseteq A$ and $B_0\subseteq B$, $\cA\subseteq\cT(A,B)$. Since
		$\lvert Q\rvert-1\ge t/2-1\ge\varepsilon n/3$, we obtain
		\[
		\begin{aligned}
			\lvert\cA\rvert\geq (\lvert Q\rvert-1)\binom{\lvert B_0\rvert}{r-2}
			&\ge\frac{\varepsilon n}{3}
			\frac{n^{r-2}}{2\cdot3^{r-2}(r-2)!}=\frac{\varepsilon}{2\cdot3^{r-1}(r-2)!}n^{r-1}
			>\beta n^{r-1}.
		\end{aligned}
		\]
		Here the second inequality holds for sufficiently large $n$
		because $r$ is fixed and $\lvert B_0\rvert\ge n/3$.
		As $u$ is $\beta$-good, at most $\beta n^{r-1}$
		sets in $\cA$ belong to $\cA\setminus\cH$. We may therefore choose
		$v\in Q\setminus u$ and $P\in\binom{B_0}{r-2}$ with
		$F'=\{u,v\}\cup P\in\cH$. Since $F'\cap G=\emptyset$, $F'\in\cF$.
		Let $E_u,E_v$ be the distinct edges of $M$ containing $u,v$, respectively.
		By the definitions of $Q$ and $B_0$, $
		E_u\cap X=u$, $E_v\cap X=v$,  $F'\cap F=\emptyset$ and $
		F'\cap X=\{u,v\}=F'\cap V(M)$.
		It follows that
		$M'=(M\setminus\{E_u,E_v\})\cup\{F',F\}$ is a matching
		of size $t$ in $\cF$ with $
		\lvert V(M')\cap X\rvert
		=\lvert V(M)\cap X\rvert+\lvert F\cap X\rvert
		>\lvert V(M)\cap X\rvert$.
		This contradicts the choice of $M$ and proves the claim.
		
		We now prove that $X\subseteq V(M)$. Otherwise, there is a vertex
		$x\in X\setminus V(M)$. By Claim~\ref{clm-fex}, every edge of $\cF$
		containing $x$ intersects $V(M)$. Then $
		d_\cF(x)\le\sum_{w\in V(M)}
		\left\lvert\{F\in\cF\mid\{x,w\}\subseteq F\}\right\rvert\le rt\binom{N-2}{r-2}
		\le\frac{r\varepsilon}{(r-2)!}n^{r-1}
		\le\frac\rho2 n^{r-1}$.
		Since $a\le s-1$, we have $D_a(n)\le D$. Note that $x\in X$ with $m_G(x)>\rho n^{r-1}$. Combining
		\eqref{eq-ext-inherit} with this gives
		\begin{equation}\label{eq-ext-uncovered}
			d_\cH(x)=d_\cF(x)+D_a(n)-m_G(x)
			<D-\frac\rho2 n^{r-1}.
		\end{equation}
		Moreover, \eqref{eq-ext-deg} implies
		$d_\cH(v)\le D+d_{\cH[B]}(v)$ for $v\in B$. Let $
		B'=\left\{v\in B\mid
		d_{\cH[B]}(v)>
		\frac{\rho}{4(r-1)}n^{r-1}\right\}$.
		Since
		\[
		\begin{aligned}
			\lvert B'\rvert\frac{\rho}{4(r-1)}n^{r-1}
			&\mathrel{ \le}\sum_{v\in B'}d_{\cH[B]}(v)\le\sum_{v\in B}d_{\cH[B]}(v)=r\lvert\cH[B]\rvert\le r\eta n^r,
		\end{aligned}
		\]
		we obtain $
		\lvert B'\rvert\le\frac{4r(r-1)\eta}{\rho}n
		\le\frac n{12}$,
		where the last inequality follows from
		$\eta\le\rho/(48r(r-1))$.
		Consequently, at least $
		\lvert B\rvert-rt-1-n/12
		\ge\left(1-\frac1r-r\varepsilon-\frac1{12}\right)n-1
		\ge r-1$
		vertices of $B\setminus(V(M)\cup x)$ satisfy
		$d_\cH(v)\le D+\rho n^{r-1}/(4(r-1))$.
		Choose distinct such vertices $v_1,\ldots,v_{r-1}$ and
		$S=\{x,v_1,\ldots,v_{r-1}\}$. By \eqref{eq-ext-uncovered}, we have
		\[
		\begin{aligned}
			\sum_{v\in S}d_\cH(v)
			&<D-\frac\rho2 n^{r-1}
			+(r-1)\left(D+\frac{\rho}{4(r-1)}n^{r-1}\right)=rD-\frac\rho4 n^{r-1}<rD.
		\end{aligned}
		\]
		The assumption $\sigma_r(\cH)>rD$ therefore forces
		$S\in\cH$. Moreover, $S\cap G=\emptyset$, which implies $S\in\cF$. Note that $S\cap V(M)=\emptyset$ with
		$x\in S\cap X$, contradicting Claim~\ref{clm-fex}. Hence
		$X\subseteq V(M)$ as desired.

		Since $
		\lvert G\rvert=a$, $\lvert M\rvert=t$ and
		$G\cap V(M)=\emptyset$, $
		\lvert V(\cH)\setminus(G\cup V(M))\rvert
		=n-a-rt=(r-1)a+c\ge(r-1)a$.
		Choose a set $
		Y\subseteq V(\cH)\setminus(G\cup V(M))$ with $\lvert Y\rvert=(r-1)a$.
		Let $\cH'=\cH[G\cup Y]$ with $N'=\lvert G\cup Y\rvert=ra$.
		For $u\in Y$, $X\subseteq V(M)$ implies that $u\notin X$. Hence $
		m_G(u)\le\rho n^{r-1}$.
		Every edge of $\cT(G,Y)\setminus \cH'$ containing $u$ intersects $G$,
		so $
		\left\lvert\{F\in\cT(G,Y)\setminus \cH'\mid u\in F\}\right\rvert
		\le m_G(u)\le\rho n^{r-1}$.
		For $g\in G$, any edge of $\cT(G,Y)\setminus \cH'$ containing $g$
		and intersecting $B$ is also an edge of $\cT(A,B)\setminus\cH$.
		There are at most $\beta n^{r-1}$ such edges because $g$
		is $\beta$-good with respect to $\cT(A,B)$. Moreover, if an edge $F$ of $\cT(G,Y)\setminus \cH'$ containing $g$
		satisfies $F\cap B=\emptyset$, then $F\subseteq A$
		and $F$ intersects $Y\cap A\subseteq R$. For each $w\in R$, there
		are at most $\binom{n-2}{r-2}$ $r$-sets containing both $g$
		and $w$. Thus
		\[
		\begin{aligned}
			\left\lvert\{F\in\cT(G,Y)\setminus \cH'\mid g\in F\}\right\rvert
			&\le\beta n^{r-1}
			+\sum_{w\in R}\binom{n-2}{r-2}\le\beta n^{r-1}+t\binom{n-2}{r-2}\le\left(\beta+\frac{\varepsilon}{(r-2)!}\right)n^{r-1}
			\le\rho n^{r-1}.
		\end{aligned}
		\]
		As $a\ge n/(2r)$, we have $N'=ra\ge n/2$ and therefore $
		\rho n^{r-1}\le2^{r-1}\rho(N')^{r-1}
		\le\rho_0(N')^{r-1}$.
		Every vertex of $\cH'$ is consequently $\rho_0$-good with respect to
		$\cT(G,Y)$. From $
		N'\ge n/2\ge8r^6$, we have $
		\frac{N'}{2r^5}<\lvert G\rvert=a=\frac{N'}r$.
		From Theorem~\ref{thm-lyy} with $W=G$ and $\lambda=\rho_0$, we obtain a matching $M_1$ of size $a$ in $\cH'$.
		Since $V(\cH')\cap V(M)=\emptyset$, $M\cup M_1$
		is a matching in $\cH$ with $t+a=s$ edges. This completes the proof.
	\end{proof}

	\begin{proof}[\textbf{Proof of Theorem~\ref{thm-main}}]
		Fix $r\ge3$ and let $C_0=c_r+r-2$ with $c_r=\left\lfloor\frac{(r-1)^{r-1}}{r^{r-2}}\right\rfloor$. By~\eqref{eq-complete-match}, it suffices to consider hypergraphs with a nonedge.
		We first prove the conclusion for $c_r\le n-rs\le C_0$.
		Let $\eta_0>0$ be supplied by Lemma~\ref{lem-extremal} for
		this choice of $C_0$, and choose $0<\eta<\eta_0/2$.
		
	Let $N_1$ be the order bound in Lemma~\ref{lem-extremal} for $C_0,\eta_0$,
		and $N_2$ be the order bound in Lemma~\ref{lem-nonext} for $K=C_0$ and $\eta$.
		Both bounds depend only on $r$ and hold uniformly for $c\in[c_r,C_0]$.
		Choose an integer $
		n_0\ge\max\left\{r,N_1,N_2,
		\frac{\lceil C_0/r\rceil+1}{(\eta_0-\eta)(r-1)!}\right\}$.
		Since $B_r(n)\le n^{r-1}/(r-1)!$, this also gives
		$(\lceil C_0/r\rceil+1)B_r(n)\le(\eta_0-\eta)n^r$
		for every $n\ge n_0$.
		Now let $n=rs+c\ge n_0$ with $c_r\le c\le C_0$ and
		$\cH$ satisfy $\sigma_r(\cH)>rD_{s-1}(n)$.
		If $\cH$ is $\eta$-nonextremal, Lemma~\ref{lem-nonext}
		gives a matching of size $s$ in $\cH$.
		If not, there is a set $B'\subseteq V(\cH)$ such that $
		\lvert B'\rvert=\left\lfloor(1-1/r)n\right\rfloor$ and $\lvert\cH[B']\rvert<\eta n^r$.
		Since $n=rs+c$, $
		\left\lfloor(1-1/r)n\right\rfloor
		=n-\left\lceil\frac nr\right\rceil
		=n-s-\left\lceil\frac cr\right\rceil$.
		Thus we can extend $B'$ to a set $B$ of size $n-s+1$ by
		adding exactly $
		\lvert B\setminus B'\rvert
		=\left\lceil\frac cr\right\rceil+1
		\le\left\lceil\frac{C_0}{r}\right\rceil+1$
		vertices. Every edge in $\cH[B]\setminus\cH[B']$ contains
		one of these added vertices, so $
			\lvert\cH[B]\rvert
			\le\lvert\cH[B']\rvert
			+\sum_{v\in B\setminus B'}d_\cH(v)\le\eta n^r+
			\left(\left\lceil\frac{C_0}{r}\right\rceil+1\right)B_r(n)
			\le\eta_0n^r$.
		Lemma~\ref{lem-extremal} therefore also gives a matching of
		size $s$ in this case.
		
		Let $
		s_0(r)=\max\left\{3,\left\lceil\frac{n_0}{r}\right\rceil\right\}$.
		We have proved that, for every $s\ge s_0(r)$ and every
		integer $c_r\le c\le C_0$, any $r$-graph on $rs+c$ vertices
		with $\sigma_r(\cH)>rD_{s-1}(rs+c)$ has a matching of size $s$.
		All constants used here depend only on $r$.
		We now reduce the general case to this statement.
		Let $s\ge s_0(r)$ and $n\ge rs+c_r$. Let $
		q=\left\lfloor\frac{n-rs-c_r}{r-1}\right\rfloor$, $
		n'=n+q$ and $s'=s+q$.
		Then $q\ge0$ and $
		(r-1)q\le n-rs-c_r<(r-1)(q+1)$.
		It follows that $
		c_r\le n'-rs'=n-rs-(r-1)q<c_r+r-1$.
		As $n'-rs'$ is an integer, we obtain $
		c_r\le n'-rs'\le c_r+r-2=C_0$, $
		s'\ge s\ge s_0(r)$ and $n'\ge rs'+c_r\ge n_0$.
		Take a set $Z$ of $q$ new vertices disjoint from $V(\cH)$
		and define an $r$-graph on $V(\cH)\cup Z$ by $
		\cH^+=\cH\cup
		\left\{F\in\binom{V(\cH)\cup Z}r\mid F\cap Z\ne\emptyset\right\}$.
		Note that $\cH^+$ and $\cH$ have precisely the same nonedges.
		For each $v\in V(\cH)$, the added edges containing $v$ in $\cH^+$ are
		the $r$-sets containing $v$ and intersecting $Z$, whose number is
		$B_r(n')-B_r(n)$. Therefore $
		d_{\cH^+}(v)=d_\cH(v)+B_r(n')-B_r(n)$.
		Taking the minimum over the common nonedges, we obtain $
		\sigma_r(\cH^+)=\sigma_r(\cH)+r\bigl(B_r(n')-B_r(n)\bigr)$.
		As $n'-s'=(n+q)-(s+q)=n-s$,
		\[
		\begin{aligned}
			D_{s'-1}(n')
			&=B_r(n')-\binom{n'-s'}{r-1}=B_r(n')-\binom{n-s}{r-1}
			=D_{s-1}(n)+B_r(n')-B_r(n).
		\end{aligned}
		\]
		Consequently,
		$\sigma_r(\cH^+)-rD_{s'-1}(n')
			=\sigma_r(\cH)-rD_{s-1}(n)>0$.
		Thus $\cH^+$ has $n'$ vertices with
		$s'\ge s_0(r)$ and $c_r\le n'-rs'\le C_0$ and satisfies
		$\sigma_r(\cH^+)>rD_{s'-1}(n')$.
		The first part of the proof, applied with $s'$ and $n'-rs'$
		in place of $s$ and $c$, gives a matching $M^+$ of size
		$s'$ in $\cH^+$.
		Since the edges of $M^+$ are pairwise disjoint, $
		\left\lvert\{F\in M^+\mid F\cap Z\ne\emptyset\}\right\rvert
		\le\sum_{F\in M^+}\lvert F\cap Z\rvert
		=\lvert V(M^+)\cap Z\rvert\le q$.
		Every edge $F$ of $M^+$ with $F\cap Z=\emptyset$ belongs to $\cH$.
		Hence $M=\{F\in M^+\mid F\cap Z=\emptyset\}$ is a matching
		in $\cH$ with $\lvert M\rvert\ge s'-q=s$, as required. This completes the proof.
	\end{proof}

	\section{Proof of Theorem~\ref{thm-r3}}\label{sec-r3}
	We extend Theorem~\ref{thm-bpr}~(ii) to $n\ge2r+1$
	by proving the remaining case $n=2r+1$.
	\begin{lemma}
		\label{lem-ore-ekr}
		Let $\cH$ be an intersecting $r$-graph on $n\ge2r+1$ vertices with $r\ge3$. Then $
		\sigma_r(\cH)\le r\binom{n-2}{r-2}=rD_1(n)$,
		with equality if and only if $\cH$ is a $1$-star.
		Consequently, every $r$-graph $\cH$ on $n\ge2r+1$ vertices with
		$\sigma_r(\cH)>rD_1(n)$ has a matching of size two.
	\end{lemma}
	\begin{proof}
		By Theorem~\ref{thm-bpr}~(ii), it suffices to consider
		$n=2r+1$. Let $
		\Omega=\binom{[n]}r$, $N=B_r(n)$ and $D=\binom{n-2}{r-2}$.
		We first prove the upper bound and show that equality forces a $1$-star.
		Write $m=\lvert \cH\rvert $ and $d_v=d_\cH(v)$. If $m=0$, the conclusion
		is immediate since $\cH$ has no edge in this case, so assume $m>0$ in the following. By Theorem~\ref{thm-ekr}, $
			m\le N$.
		Fix an edge $F\in\cH$. Recall that $KG(n,r)$ is
		$\binom{n-r}{r}$-regular. Since $n=2r+1$, the $r+1$ neighbors of $F$ are precisely the
		$r$-subsets of $[n]\setminus F$. Since $\cH$ is intersecting,
		all these neighbors are nonedges of $\cH$. Hence each such set
		$S$ satisfies $\sum_{v\in S}d_v\ge\sigma_r(\cH)$.
	Summing over all these nonedges $S\in\binom{[n]\setminus F}{r}$, we have
		\begin{align}
			(r+1)\sigma_r(\cH)
			&\le
			\sum_{S\in\binom{V(\cH)\setminus F}{r}}
			\sum_{v\in S}d_v.\label{eq-ore-edge}
		\end{align}
		For each fixed $v\in V(\cH)\setminus F$,
		the number of $r$-subsets of $V(\cH)\setminus F$
		containing $v$ is $
		\binom{\lvert V(\cH)\setminus F\rvert -1}{r-1}
		=\binom r{r-1}=r$.
		Consequently, each $d_v$ with $v\notin F$
		appears exactly $r$ times on the right-hand side of \eqref{eq-ore-edge}. Hence
		\[
		\sum_{S\in\binom{V(\cH)\setminus F}{r}}
		\sum_{v\in S}d_v=\sum_{v\notin F}
		d_v\binom r{r-1}=r\sum_{v\notin F}d_v{.}\]
		It follows that $
		(r+1)\sigma_r(\cH)
		\le r\sum_{v\notin F}d_v$ for every $F\in\cH$.
		Summing \eqref{eq-ore-edge} over all $m$ edges
		$F\in\cH$, we obtain \begin{align}
			(r+1)m\sigma_r(\cH)
			&\le r\sum_{F\in\cH}\sum_{v\notin F}d_v.\label{eq-ore-sum}
		\end{align}
		Interchanging the order of summation in \eqref{eq-ore-sum} yields
		$\sum_{F\in\cH}\sum_{v\notin F}d_v
		=\sum_{v\in V(\cH)}
		d_v\lvert \{F\in\cH\mid v\notin F\}\rvert =\sum_{v\in V(\cH)}d_v(m-d_v)$.
		Combining this with $
		\sum_{v\in V(\cH)}d_v=rm$, we have
		\begin{align}
			(r+1)m\sigma_r(\cH)\le r\sum_{v\in V(\cH)}d_v(m-d_v)&=r\bigg(m\sum_vd_v-\sum_vd_v^2\bigg)=r\left(rm^2-\sum_vd_v^2\right).\label{eq-ore-count}
		\end{align}
		
		Next we prove that
		\begin{equation}\label{eq-deg-square}
			\sum_vd_v^2\ge\frac{3r-1}{4}m^2
			-\frac{r-2}{2}(N-D)m.
		\end{equation}
		Let $A$ be the adjacency matrix of the Kneser graph $KG(n,r)$,
		that is, $A_{S,R}=\chi_{\{S\cap R=\emptyset\}}$ for $S,R\in\Omega$.
		Since $KG(2r+1,r)$ is $(r+1)$-regular, each row of $A$
		sums to $r+1$. Hence
		\begin{equation}\label{eq-kg-regular}
			A\mathbf{1}=(r+1)\mathbf{1}.
		\end{equation}
		For each $v\in[n]$, let $u_v:\Omega\to\mathbb R$ be the function
		$u_v(S)=\chi_{\{v\in S\}}$, and write
		$\mathbf{u}_v=(u_v(S))_{S\in\Omega}$ for its column vector of values.
		We prove the following claim.
		\begin{claim}\label{clm-kg-bound}
			If a vector $\mathbf{h}\in\mathbb R^\Omega$ satisfies $\langle \mathbf{h},\mathbf{u}_v\rangle=0$
			for every $v\in[n]$, then
			\begin{equation}\label{eq-kg-bound}
				\langle \mathbf{h},A\mathbf{h}\rangle\ge-(r-2)\lVert \mathbf{h}\rVert ^2.
			\end{equation}
		\end{claim}
		\noindent\textbf{Proof of Claim~\ref{clm-kg-bound}.}
		By Lemma~\ref{lem-kg-spec}, the eigenvalues of $A$ are
		$\lambda_j=(-1)^j(r+1-j)$ for $0\le j\le r$.
		Let $W_j$ be the eigenspace corresponding to $\lambda_j$.
		Since $A$ is real and symmetric,
		$\mathbb R^\Omega=W_0\oplus\cdots\oplus W_r$ is an orthogonal direct sum.
		For each $S\in\Omega$, $\lvert S\rvert=r$, so
		$\sum_{v=1}^n\mathbf{u}_v=r\mathbf{1}$. Then $\mathbf{h}$ is orthogonal to
		$\mathbf{1}$ and to every $\mathbf{u}_v$ by $\langle \mathbf{h},\mathbf{u}_v\rangle=0$
		for every $v\in[n]$.
		By Lemma~\ref{lem-kg-spec}, $
		W_0=\operatorname{span}\{\mathbf{1}\}$ and $
		W_1=\bigl\{\mathbf{f}\in\operatorname{span}\{\mathbf{u}_v\mid v\in[n]\}
		\mid \langle \mathbf{f},\mathbf{1}\rangle=0\bigr\}$, which implies that
		$\mathbf{h}\in W_2\oplus\cdots\oplus W_r$.
		Write $\mathbf{h}=\sum_{j=2}^r\mathbf{h}_j$ with $\mathbf{h}_j\in W_j$.
		Orthogonality and $\lambda_j\ge-(r-2)$ for $2\le j\le r$ give
		\[
		\langle\mathbf{h},A\mathbf{h}\rangle
		=\sum_{j=2}^r\lambda_j\lVert \mathbf{h}_j\rVert ^2
		\ge-(r-2)\sum_{j=2}^r\lVert \mathbf{h}_j\rVert ^2
		=-(r-2)\lVert \mathbf{h}\rVert ^2.
		\]
		This completes the proof of Claim~\ref{clm-kg-bound}.
		
		Define the functions $f:\Omega\to\mathbb R$ by $
		f(S)=\chi_{\cH}(S)=
		\begin{cases}
			1,&S\in\cH,\\
			0,&S\notin\cH,
		\end{cases}$ and $\ell:\Omega\to\mathbb R$ by $
		\ell(S)=\sum_{v\in S}d_v$.
		Let $\mathbf{f}=(f(S))_{S\in\Omega}$ and
		$\boldsymbol{\ell}=(\ell(S))_{S\in\Omega}$ be their corresponding
		column vectors. For every $S\in\Omega$,
		\[
		\left(\sum_{v=1}^n d_v\mathbf{u}_v\right)_S
		=\sum_{v=1}^n d_vu_v(S)
		=\sum_{v\in S}d_v
		=\ell(S)=(\boldsymbol{\ell})_S.
		\]
		Hence
		$\boldsymbol{\ell}=\sum_{v=1}^n d_v\mathbf{u}_v$.
		Let $
		\mathbf{h}=\mathbf{f}-\frac{\boldsymbol{\ell}}{N-D}
		+\frac{(r-1)m}{2(N-D)}\mathbf{1}$. For each $v\in[n]$,\begin{align*}
			\langle \mathbf{f},\mathbf{u}_v\rangle&=\sum_{S\in\Omega}f(S)u_v(S)=\lvert \{S\in\cH\mid v\in S\}\rvert =d_v,\quad
			\langle\mathbf{1},\mathbf{u}_v\rangle=\sum_{S\in\Omega}u_v(S)=\lvert \{S\in\Omega\mid v\in S\}\rvert =N,\\
			\langle\boldsymbol{\ell},\mathbf{u}_v\rangle&
			=\sum_{w=1}^nd_w\langle\mathbf{u}_w,\mathbf{u}_v\rangle=Nd_v+D\sum_{w\ne v}d_w=Nd_v+D(rm-d_v)=(N-D)d_v+Drm.
		\end{align*}
		Using $rD=(r-1)N/2$, we obtain for every $v\in[n]$,
		\begin{align*}
			\langle \mathbf{h},\mathbf{u}_v\rangle&
			=\langle \mathbf{f},\mathbf{u}_v\rangle
			-\frac{\langle\boldsymbol{\ell},\mathbf{u}_v\rangle}{N-D}+\frac{(r-1)m}{2(N-D)}
			\langle\mathbf{1},\mathbf{u}_v\rangle
			=d_v
			-\frac{(N-D)d_v+Drm}{N-D}+\frac{(r-1)mN}{2(N-D)}\\
			&=\frac{-Drm+(r-1)mN/2}{N-D}=0.
		\end{align*}
		Let $L=\operatorname{span}\{\mathbf{u}_v\mid v\in[n]\}$.
		Since $\boldsymbol{\ell},\mathbf{1}\in L$, the definition of
		$\mathbf{h}$ gives $\mathbf{f}-\mathbf{h}\in L$.
		As $\mathbf{h}$ is orthogonal to $L$, we have
		$\langle\mathbf{h},\mathbf{f}-\mathbf{h}\rangle=0$.
		Consequently, $
		\lVert \mathbf{h}\rVert ^2=\langle\mathbf{h},\mathbf{h}\rangle
		=\langle\mathbf{h},\mathbf{f}\rangle=\langle\mathbf{f},\mathbf{f}\rangle
		-\frac{\langle\boldsymbol{\ell},\mathbf{f}\rangle}{N-D}
		+\frac{(r-1)m}{2(N-D)}
		\langle\mathbf{1},\mathbf{f}\rangle$.
		Moreover, $
		\langle\mathbf{f},\mathbf{f}\rangle
		=\sum_{S\in\Omega}f(S)^2=m$, $
		\langle\mathbf{1},\mathbf{f}\rangle=m$ and $
		\langle\boldsymbol{\ell},\mathbf{f}\rangle=\left\langle
		\sum_{v=1}^n d_v\mathbf{u}_v,\mathbf{f}
		\right\rangle=\sum_{v=1}^n d_v
		\langle\mathbf{u}_v,\mathbf{f}\rangle=\sum_{v=1}^n d_v^2$.
		Therefore,
		\begin{equation}
			\label{eq-h-norm}
			\lVert \mathbf{h}\rVert ^2
			=m-\frac{\sum_vd_v^2}{N-D}
			+\frac{(r-1)m^2}{2(N-D)}.
		\end{equation}
		We next compute $\langle\mathbf{h},A\mathbf{h}\rangle$.
		For $v\in[n]$ and $S\in\Omega$, $
		(A\mathbf{u}_v)_S
		=\sum_{R\in\Omega}A_{S,R}u_v(R)
		=\lvert \{R\in\Omega\mid R\cap S=\emptyset,\ v\in R\}\rvert $.
		If $v\in S$, then $
		(A\mathbf{u}_v)_S=0$. If $v\notin S$, then the other $r-1$ vertices can
		be chosen from the $r$-set $[n]\setminus(S\cup v)$,
		so $
		(A\mathbf{u}_v)_S=\binom r{r-1}=r$.
		Hence $A\mathbf{u}_v=r(\mathbf{1}-\mathbf{u}_v)$.
		Recall that $\boldsymbol{\ell}=\sum_vd_v\mathbf{u}_v$
		and $\sum_vd_v=rm$. Then $A\boldsymbol{\ell}=\sum_vd_vA\mathbf{u}_v=r\sum_vd_v(\mathbf{1}-\mathbf{u}_v)=r\left(rm\mathbf{1}-\boldsymbol{\ell}\right)$.

		Because $\mathbf{1}\in L$ and
		$\mathscr{A}(\mathbf{u}_v)=A\mathbf{u}_v=r(\mathbf{1}-\mathbf{u}_v)$,
		the subspace $L$ is $\mathscr{A}$-invariant.
		Since $\mathbf{f}-\mathbf{h}\in L$,
		$A(\mathbf{f}-\mathbf{h})=\mathscr{A}(\mathbf{f}-\mathbf{h})\in L$.
		Combining this with the fact that $\mathbf{h}$ is orthogonal to $L$, we have
		$\langle\mathbf{h},A(\mathbf{f}-\mathbf{h})\rangle=0$, and hence
		$\langle\mathbf{h},A\mathbf{h}\rangle
		=\langle\mathbf{h},A\mathbf{f}\rangle$.
		By the definition of $\mathbf{h}$,
		\begin{align}\label{eq-h-expand}
			\langle\mathbf{h},A\mathbf{h}\rangle
			&=\langle\mathbf{f},A\mathbf{f}\rangle
			-\frac{\langle\boldsymbol{\ell},A\mathbf{f}\rangle}{N-D}\quad+\frac{(r-1)m}{2(N-D)}
			\langle\mathbf{1},A\mathbf{f}\rangle.
		\end{align}
		Note that $\cH$ is intersecting. Thus $
		\langle\mathbf{f},A\mathbf{f}\rangle=\sum_{S,R\in\Omega}
		f(S)A_{S,R}f(R)=\sum_{\substack{S,R\in\cH\\S\cap R=\emptyset}}1
		=0$.
		Since $A$ is symmetric, we also have $
		\langle\boldsymbol{\ell},A\mathbf{f}\rangle=\langle A\boldsymbol{\ell},\mathbf{f}\rangle=r\langle
		rm\mathbf{1}-\boldsymbol{\ell},\mathbf{f}
		\rangle=r\left(
		rm\langle\mathbf{1},\mathbf{f}\rangle
		-\langle\boldsymbol{\ell},\mathbf{f}\rangle
		\right)=r\left(rm^2-\sum_vd_v^2\right)$ and $
		\langle\mathbf{1},A\mathbf{f}\rangle=\langle A\mathbf{1},\mathbf{f}\rangle=(r+1)\langle\mathbf{1},\mathbf{f}\rangle=(r+1)m$.
		Substituting these identities into \eqref{eq-h-expand} gives
		\begin{equation}
			\label{eq-h-quadratic}
			\langle\mathbf{h},A\mathbf{h}\rangle
			=\frac{r\sum_vd_v^2}{N-D}
			-\frac{(r^2+1)m^2}{2(N-D)}.
		\end{equation}
		
		By Claim~\ref{clm-kg-bound}, $
		\langle\mathbf{h},A\mathbf{h}\rangle
		\ge -(r-2)\lVert \mathbf{h}\rVert ^2$.
		Substituting
		\eqref{eq-h-norm} and
		\eqref{eq-h-quadratic}, we obtain
		\[
		\begin{aligned}
			\frac{r\sum_vd_v^2}{N-D}
			-\frac{(r^2+1)m^2}{2(N-D)}\quad\ge -(r-2)
			\left(
			m-\frac{\sum_vd_v^2}{N-D}
			+\frac{(r-1)m^2}{2(N-D)}
			\right).
		\end{aligned}
		\]
		Multiplying by $N-D>0$ and collecting terms gives
		\[
		2\sum_vd_v^2
		\ge\frac{r^2+1-(r-2)(r-1)}{2}m^2-(r-2)(N-D)m
		=\frac{3r-1}{2}m^2-(r-2)(N-D)m.
		\]
		Dividing both sides by $2$ proves \eqref{eq-deg-square}.

		Finally, combining
		\eqref{eq-ore-count} with
		\eqref{eq-deg-square} yields
		\[
		\begin{aligned}
			(r+1)m\sigma_r(\cH)
			&\le r\left(rm^2-\sum_vd_v^2\right)\le r\left(
			rm^2-\frac{3r-1}{4}m^2
			+\frac{r-2}{2}(N-D)m
			\right)=r\left(
			\frac{r+1}{4}m^2
			+\frac{r-2}{2}(N-D)m
			\right).
		\end{aligned}
		\]
		Since $m>0$, division by $(r+1)m$ gives $
		\sigma_r(\cH)
		\le\frac r4m+\frac{r(r-2)}{2(r+1)}(N-D)$.
		As $n=2r+1$, we have $
		\frac DN=\frac{r-1}{n-1}
		=\frac{r-1}{2r}$. Hence $
		N-D=\frac{r+1}{2r}N$ and $
		rD=\frac{r-1}{2}N$.
		It follows that
		\begin{align}\label{eq-ore-bound}
			\sigma_r(\cH)
			&\le\frac r4m
			+\frac{r(r-2)}{2(r+1)}
			\cdot\frac{r+1}{2r}N=\frac{rm+(r-2)N}{4}=rD-\frac r4(N-m)\le rD,
		\end{align}
		where the last inequality follows from
		$m\le N$. If $\sigma_r(\cH)=rD$, \eqref{eq-ore-bound} implies that $\lvert\cH\rvert=m=N=B_r(n)$. Since $n>2r$,
		$\cH$ is a $1$-star by
		Theorem~\ref{thm-ekr}.
		
		Conversely, let $\cH$ be a $1$-star with center $x$. Every nonedge $S$ satisfies $x\notin S$ and every vertex $v\ne x$ has degree
		$\binom{n-2}{r-2}=D$. Thus each nonedge has degree sum $rD$,
		so $\sigma_r(\cH)=rD$.
		Finally, an $r$-graph with no matching of size two is intersecting.
		The bound just proved, together with
		Theorem~\ref{thm-bpr}~(ii), therefore implies the stated
		matching conclusion for every $n\ge2r+1$.
	\end{proof}
	
	For the rest of this section, let $r=3$ and write
	$D_m(n)=B_3(n)-\binom{n-m-1}{2}=m(2n-m-3)/2$.
	\begin{lemma}\label{lem-r3-cap}
		Let $E$, $Y$ be disjoint sets with $E=\{1,2,3\}$, $\lvert Y\rvert=q$ and $K\subseteq Y$. Let
		$\cH\subseteq\binom{E\cup Y}{3}$ be an intersecting family all of whose members
		intersect $E$. For each $a\in E$, define $G_a=\{P\in\tbinom Y2\mid a\cup P\in\cH\}$.
		\begin{enumerate}
			\item If $\lvert K\rvert\ge5$, then $
			\sum_{z\in K}d_{\cH}(z)\le(q+1)\lvert K\rvert$.
			\item If $\lvert K\rvert=3$ and $q\ge5$, then $
			\sum_{z\in K}d_{\cH}(z)\le3q+6$.
			\item If $q\ge12$ and $\nu(G_a)\le4$ for every $a\in E$, then
			$\sum_{z\in Y}d_{\cH}(z)\le10q-20$.
		\end{enumerate}
	\end{lemma}
	\begin{proof}
		Let $
		k=\lvert K\rvert $.
		Since $\cH$ is intersecting, the graphs $G_1,G_2,G_3$ are pairwise
		cross-intersecting.  Otherwise, there exist $xy\in G_a$ and $uv\in G_b$ with
		$a\ne b$ and $\{x,y\}\cap\{u,v\}=\emptyset$. Then the edges $\{a,x,y\}$ and $\{b,u,v\}$ of $\cH$ are disjoint, contradicting the assumption that $\cH$ is intersecting.
		
		We next bound $\sum_{z\in K}d_{\cH}(z)$
		by classifying the edges of $\cH$
		according to their intersections with $E$. {For each edge $F\in\cH$}, $1\leq\lvert F\cap E\rvert\leq3$ since every edge of $\cH$ intersects $E$. Moreover, $
		\sum_{z\in K}d_{\cH}(z)
		=\sum_{F\in\cH}\lvert F\cap K\rvert$.
		For each $a\in E$, let $
		C(G_a):=\bigcap_{F\in G_a}F $ and
		$C(G_a):=Y$ if $G_a=\emptyset$.
		Then $\lvert C(G_a)\rvert=2$ if $G_a$ has a single edge, and $\lvert C(G_a)\rvert\le1$ if $G_a$ has at least two edges.
		
		Define $
		f(G_a):=
		\sum_{z\in K}d_{G_a}(z)
		+
		\lvert K\cap C(G_a)\rvert$.
		First, if $F\in\cH$ with
		$\lvert F\cap E\rvert=1$, then $
		F=a\cup P$ with $a\in E$ and $P\in G_a$.
		Since $K\subseteq Y$ and $E\cap Y=\emptyset$,
		we have $
		\lvert F\cap K\rvert=\lvert P\cap K\rvert$.
		Then
		\[
		\sum_{\substack{F\in\cH\\\lvert F\cap E\rvert=1}}\lvert F\cap K\rvert=\sum_{a\in E}\sum_{P\in G_a}
		\lvert P\cap K\rvert
		=\sum_{a\in E}\sum_{z\in K}d_{G_a}(z).
		\]
		{Second, if $F\in\cH$ with}
		$\lvert F\cap E\rvert=2$, then $
		F=(E\setminus a)\cup y$ with $a\in E$ and $y\in Y$.
		Since $K\subseteq Y$, such an edge contributes
		one to $\sum_{z\in K}d_{\cH}(z)$ if $y\in K$
		and zero otherwise.
		Fix $a\in E$ and suppose that $
		(E\setminus a)\cup y\in\cH$.
		For every $P\in G_a$, we also have
		$a\cup P\in\cH$.
		Since $\cH$ is intersecting, $\bigg((E\setminus a)\cup y\bigg)\cap(a\cup P)\neq\emptyset$. As $
		(E\setminus a)\cap(a\cup P)
		=\emptyset$,
		it follows that $y\in P$ for every $P\in G_a$.
		Hence $y\in C(G_a)$ for $G_a\ne\emptyset$.
		The same inclusion holds when $G_a=\emptyset$
		by the definition $C(G_a)=Y$.
		Consequently, for each fixed $a\in E$,
		\[
		\begin{aligned}
			{\sum_{\substack{F\in\cH\\F\cap E=E\setminus a}}\lvert F\cap K\rvert}=\left\lvert
			\left\{y\in K\mid
			(E\setminus a)\cup y\in\cH
			\right\}\right\rvert\le\lvert K\cap C(G_a)\rvert.
		\end{aligned}
		\]
		{Finally, if $F\in\cH$ with}
		$\lvert F\cap E\rvert=3$, then $F=E$. Since
		{
			$E\cap K=\emptyset$, we have $\lvert E\cap K\rvert=0$.
		}
		
		Therefore,
		\begin{align}\label{eq-cap-reduce}
			\sum_{z\in K}d_{\cH}(z)
			=\sum_{a\in E}
			\left(
			\sum_{P\in G_a}\lvert P\cap K\rvert
			+\left\lvert
			\left\{y\in K\mid
			(E\setminus a)\cup y\in\cH
			\right\}\right\rvert
			\right)\notag&\le\sum_{a\in E}
			\left(
			\sum_{z\in K}d_{G_a}(z)
			+\lvert K\cap C(G_a)\rvert
			\right)\notag\\
			&=\sum_{a\in E}f(G_a).
		\end{align}

		For parts~\textup{(i)} and~\textup{(ii)}, we have $q\ge5$ and $k\ge3$.
		We now distinguish two cases according to the matching numbers of
		$G_1,G_2,G_3$.
		
		\medskip
		\noindent\textbf{Case~1.} $\nu(G_a)\le1$
		for every $a\in E$.
		
		We first describe the possible structure of such a graph.  Let $G_a$
		contain at least two edges. Since $\nu(G_a)\le1$, any two edges of
		$G_a$ intersect. Thus we can choose $
		xy,xz\in E(G_a)$.
		If every edge of $G_a$ contains $x$, then $G_a$ is a star with centre $x$.  Otherwise, let $F\in G_a$ be an edge such that $x\notin F$.  Since $F\cap xy\neq \emptyset$ and $F\cap xz\neq \emptyset$, $
		F=yz$.
		Moreover, any edge of $G_a$ must intersect
		all three edges $
		xy, xz$ and $yz$, hence must be one of these three edges.  Thus $G_a$ is a triangle.
		Consequently, every graph $G_a$ with $\nu(G_a)\le1$ is either
		empty, a star, or a triangle.
		
		We next record the corresponding bounds for $f(G_a)$.  If
		$G_a=\emptyset$, by the convention $C(G_a)=Y$, $
		f(G_a)=k$.
		If $G_a=\{F\}$, then $f(G_a)=2\lvert F\cap K\rvert\le4$.
		If $G_a$ is a two-edge star, its total degree sum is $4$ and
		$\lvert C(G_a)\rvert\le1$, so $f(G_a)\le4+1=5$. If $G_a$ is a
		triangle, then $C(G_a)=\emptyset$ and hence $
		f(G_a)\le6$. More generally, suppose that $G_a$ is contained in a star with centre
		$x$ and has at least two edges.  Then $
		C(G_a)=x$.
		If $x\in K$, then $x$ contributes at most $q-1$ to
		$\sum_{z\in K}d_{G_a}(z)$, while every vertex of
		$K\setminus x$ has degree at most one in $G_a$.  Therefore $
		f(G_a)
		\le
		(q-1)+(k-1)+1
		=
		q+k-1$.
		If $x\notin K$, then $
		f(G_a)\le k\le q+k-1$.
		Thus we have, for $\nu(G_a)\le1$,
		\begin{equation}\label{eq-r3-star}
			f(G_a)\le q+k-1.
		\end{equation}

		Suppose that one of the graphs, say $G_1$, is a star with
		at least three edges with centre $x$.  Since $G_1,G_2,G_3$ are pairwise
		cross-intersecting, every edge of $G_2$ or
		$G_3$ must intersect each of these three star edges. {Hence each nonempty graph among
			$G_2,G_3$} is also contained in a star with centre $x$.
		By \eqref{eq-r3-star},
		\begin{equation}\label{eq-common-star}
			f(G_1)+f(G_2)+f(G_3)
			\le
			3(q+k-1).
		\end{equation}
		If no $G_a$ is a star with at least three edges, then each $G_a$ is empty, a single edge,
		a two-edge star, or a triangle.  Hence
		\begin{equation}\label{eq-small-graphs}
			f(G_1)+f(G_2)+f(G_3)
			\le
			3\max\{k,6\}
			=
			\max\{3k,18\}.
		\end{equation}
		
		\medskip
		\noindent\textbf{Case~2.} $\nu(G_a)\ge2$ for some $a\in E$.
		
		Without loss of generality, assume that $G_1$ contains two disjoint
		edges. Let $
		\{x_1,x_2\},\ \{x_3,x_4\}\in G_1$ such that $x_1,x_2,x_3,x_4$ are distinct. Every edge of $G_2$ and
		$G_3$ must intersect both $\{x_1,x_2\}$ and $\{x_3,x_4\}$.
		Therefore $
		G_2,G_3
		\subseteq
		\bigl\{
		\{x_1,x_3\},
		\{x_1,x_4\},
		\{x_2,x_3\},
		\{x_2,x_4\}
		\bigr\}$,
		which gives a cycle $C_4$. Let $
		j=
		\left\lvert
		\{a\in\{2,3\}\mid G_a\ne\emptyset\}
		\right\rvert$. Next we consider the cases based on the value of $j$.
		
		\medskip
		\noindent\textbf{Case~2.1.} $j=0$.
		
		In this case, $G_2=G_3=\emptyset$, and therefore $
		f(G_2)=f(G_3)=k$.
		Since $G_1$ contains two disjoint edges, $
		C(G_1)=\emptyset$.
		Each vertex of $K$ has degree at most $q-1$ in $G_1$, so $
		f(G_1)
		=
		\sum_{z\in K}d_{G_1}(z)
		\le
		k(q-1)$.
		Consequently,
		\begin{equation}\label{eq-r3-j0}
			f(G_1)+f(G_2)+f(G_3)
			\le
			k(q+1).
		\end{equation}
		
		\medskip
		\noindent\textbf{Case~2.2.} $j=1$.
		
		In this case, assume that $
		G_2\ne\emptyset$ and $
		G_3=\emptyset$.
		Then either $G_2$ consists of a single edge, $G_2$
		is a two-edge star, or $G_2$ contains two disjoint edges.
		
		First suppose that $G_2$ consists of a single edge $x_1x_3$.
		Hence every
		edge of $G_1$ contains $x_1$ or $x_3$. Since $G_1$ contains two
		disjoint edges, $
		C(G_1)=\emptyset$.
		For every $
		z\in K\setminus\{x_1,x_3\}$,
		an edge of $G_1$ containing $z$ must be either $
		\{x_1,z\}$ or $
		\{x_3,z\}$.
		Hence $
		d_{G_1}(z)\le2$.
		{
			The vertices $x_1$ and $x_3$, when they belong to $K$, each have
			degree at most $q-1$. Since $q\ge5$ and $k\ge3$, it follows that
			$f(G_1)\le2(q-1)+2(k-2)$.
		}
		
		Moreover, $
		f(G_2)\le4$ and $
		f(G_3)=k$.
		Consequently,
		\begin{equation}\label{eq-j1-single}
			f(G_1)+f(G_2)+f(G_3)
			\le
			2(q-1)+3k.
		\end{equation}

		Next suppose that $G_2$ is a two-edge star with $
		G_2=
		\bigl\{
		\{x_1,x_3\},
		\{x_1,x_4\}
		\bigr\}$.
		Every edge of $G_1$ must intersect both $\{x_1,x_3\}$ and
		$\{x_1,x_4\}$. Therefore every edge of $G_1$ either contains
		$x_1$ or equals $\{x_3,x_4\}$. Hence $
		G_1
		\subseteq
		\bigl\{
		\{x_1,v\}\mid v\in Y\setminus x_1
		\bigr\}
		\cup
		\bigl\{\{x_3,x_4\}\bigr\}$.
		Since $G_1$ contains two disjoint edges, $
		C(G_1)=\emptyset$.
		If $x_1\in K$, then $x_1$ contributes at most $q-1$ to
		$\sum_{z\in K}d_{G_1}(z)$, the vertices in
		$K\setminus x_1$ contribute at most $k-1$ through the star
		centred at $x_1$ and the edge $\{x_3,x_4\}$ contributes
		at most $2$. Thus
		$f(G_1)\le(k-1)+(q-1)+2=q+k$. If $x_1\notin K$, then the star contributes at most $k$ and
		$\{x_3,x_4\}$ contributes at most $2$, so $
		f(G_1)\le k+2\le q+k$.
		Moreover, $
		f(G_2)\le5$ and $
		f(G_3)=k$.
		It follows that
		\begin{equation}\label{eq-j1-star}
			f(G_1)+f(G_2)+f(G_3)
			\le
			q+2k+5.
		\end{equation}
		
		Finally, suppose that $G_2$ contains two disjoint edges. Since
		$G_2$ is a subgraph of the four-cycle, we may assume that $
		\{x_1,x_3\},\ \{x_2,x_4\}\in G_2$.
		Every edge of $G_1$ must intersect both of these edges. Hence $
		G_1
		\subseteq
		\bigl\{
		\{x_1,x_2\},
		\{x_1,x_4\},
		\{x_2,x_3\},
		\{x_3,x_4\}
		\bigr\}$,
		which is again a four-cycle. Since both $G_1$ and $G_2$ contain
		two disjoint edges, $
		C(G_1)=C(G_2)=\emptyset$.
		A four-cycle has four edges and total degree sum $8$, so $
		f(G_a)\le8$ with $a\in[2]$.
		Together with $f(G_3)=k$, this gives
		\begin{equation}\label{eq-j1-match}
			f(G_1)+f(G_2)+f(G_3)
			\le
			16+k.
		\end{equation}
		Combining the three possibilities in Case~2.2, we obtain
		\begin{equation}\label{eq-r3-j1}
			f(G_1)+f(G_2)+f(G_3)
			\le
			\max\{2(q-1)+3k,\ q+2k+5,\ 16+k\}.
		\end{equation}
		
		\medskip
		\noindent\textbf{Case~2.3.} $j=2$.
		
		In this case, $
		G_2\ne\emptyset$ and $
		G_3\ne\emptyset$.
		We first suppose that both $G_2$ and $G_3$ consist of a single
		edge. Since $G_2$ and $G_3$ are cross-intersecting, these two edges
		are either equal or intersect.
		Suppose first that they are equal. Assume that $
		G_2=G_3=
		\bigl\{
		\{x_1,x_3\}
		\bigr\}$.
		Every edge of $G_1$ must intersect $\{x_1,x_3\}$, so every edge of
		$G_1$ contains $x_1$ or $x_3$. Since $G_1$ contains two
		disjoint edges, $
		C(G_1)=\emptyset$.
		As in Case~2.2, $
		f(G_1)
		\le
		2(q-1)+2(k-2)$.
		Since $
		f(G_2)\le4$ and $
		f(G_3)\le4$,
		we obtain $
		f(G_1)+f(G_2)+f(G_3)
		\le
		2q+2k+2$.
		Now suppose that the single edges of $G_2$ and $G_3$ are distinct.
		Since they intersect, we may assume that
		$
		G_2=
		\bigl\{
		\{x_1,x_3\}
		\bigr\}$ and $
		G_3=
		\bigl\{
		\{x_1,x_4\}
		\bigr\}$.
		Every edge of $G_1$ must intersect both $\{x_1,x_3\}$ and
		$\{x_1,x_4\}$. Hence every edge of $G_1$ either contains $x_1$
		or equals $\{x_3,x_4\}$, so $
		f(G_1)\le q+k$.
		Combining this with $q+k+8\le2q+2k+2$, we have
		\begin{align}\label{eq-j2-single}
			f(G_1)+f(G_2)+f(G_3)&\le
			q+k+8\le
			2q+2k+2.
		\end{align}
		
		Next suppose that neither $G_2$ nor $G_3$ contains two disjoint
		edges, but at least one of them is a two-edge star. Without loss of generality,
		assume that $
		G_2=
		\bigl\{
		\{x_1,x_3\},
		\{x_1,x_4\}
		\bigr\}$.
		Every edge of $G_1$ must intersect both edges of $G_2$, so $
		G_1
		\subseteq
		\bigl\{
		\{x_1,v\}\mid v\in Y\setminus x_1
		\bigr\}
		\cup
		\bigl\{\{x_3,x_4\}\bigr\}$.
		As above, $
		f(G_1)\le q+k$.
		Since $G_3$ is a nonempty subgraph of the four-cycle with $\nu(G_3)\le1$,
		it is either a single edge or a two-edge star. Therefore $
		f(G_2)\le5$ and $
		f(G_3)\le5$.
		It follows that
		\begin{equation}\label{eq-j2-star}
			f(G_1)+f(G_2)+f(G_3)
			\le
			q+k+10.
		\end{equation}
		
		Finally, suppose that at least one of $G_2,G_3$ contains two
		disjoint edges. Without loss of generality, assume that $
		\{x_1,x_3\},\ \{x_2,x_4\}\in G_2$.
		Every edge of $G_1$ and every edge of $G_3$ must intersect both
		$\{x_1,x_3\}$ and $\{x_2,x_4\}$. Consequently, $
		G_1,G_3
		\subseteq
		\bigl\{
		\{x_1,x_2\},
		\{x_1,x_4\},
		\{x_2,x_3\},
		\{x_3,x_4\}
		\bigr\}$. Thus each of $G_1,G_2,G_3$ is contained in a
		four-cycle on the four vertices $x_1,x_2,x_3,x_4$. We claim that $
		f(G_a)\le\min\{8,2k\}$ with $a\in[3]$.
		If $G_a$ contains two disjoint edges, then $
		C(G_a)=\emptyset$.
		Since a four-cycle has total degree sum $8$ and maximum degree $2$,
		$f(G_a)
		=
		\sum_{z\in K}d_{G_a}(z)
		\le
		\min\{8,2k\}$.
		If $G_a$ consists of a single edge or is a two-edge star, then $
		f(G_a)\le4 $ or $
		f(G_a)\le5$,
		respectively. Since $k\ge3$,
		$f(G_a)\leq\min\{8,2k\}$. Hence
		\begin{equation}\label{eq-j2-match}
			f(G_1)+f(G_2)+f(G_3)
			\le
			\min\{24,6k\}.
		\end{equation}
		Combining
		\eqref{eq-j2-single}--\eqref{eq-j2-match}, we obtain
		\begin{equation}\label{eq-r3-j2}
			f(G_1)+f(G_2)+f(G_3)
			\le
			\max\{2q+2k+2,\ q+k+10,\ \min\{24,6k\}\}.
		\end{equation}
		
		We now complete the proof.
		For part~\textup{(i)}, we have $
		q\ge k\ge5$. Then the bounds in \eqref{eq-common-star}, \eqref{eq-small-graphs}, \eqref{eq-r3-j0}, \eqref{eq-r3-j1} and \eqref{eq-r3-j2} are at
		most $k(q+1)$.
		Therefore, by \eqref{eq-cap-reduce}, $
		\sum_{z\in K}d_{\cH}(z)
		\le
		k(q+1)
		=
		(q+1)\lvert K\rvert$.
		For part~\textup{(ii)}, we have $
		k=3 $ and $
		q\ge5$.
		Then the bounds in \eqref{eq-common-star}, \eqref{eq-small-graphs}, \eqref{eq-r3-j0}, \eqref{eq-r3-j1} and \eqref{eq-r3-j2} are at most $3q+6$ for $q\ge5$. Hence, again by
		\eqref{eq-cap-reduce}, $
		\sum_{z\in K}d_{\cH}(z)
		\le
		3q+6$.
		{
			This proves \textup{(i)} and \textup{(ii)}.
			For \textup{(iii)}, take $K=Y$, so $k=q$ and
			$f(G_a)=2\lvert G_a\rvert+\lvert C(G_a)\rvert$.
			If some $G_a$ contains three disjoint edges, cross-intersection forces
			the other two graphs to be empty. Since $\nu(G_a)\le4$ and $q\ge12$,
			Lemma~\ref{lem-eg} gives $\lvert G_a\rvert\le4q-10$ with
			$C(G_a)=\emptyset$. Hence $\sum_{a\in E}f(G_a)\le10q-20$.
			Otherwise, every $G_a$ has matching number at most two. If one of them,
			say $G_1$, contains two disjoint edges, Lemma~\ref{lem-eg} gives
			$f(G_1)=2\lvert G_1\rvert\le4q-6$.
			As in Case~2, $G_2$ and $G_3$ are subgraphs of a four-cycle. Then  $f(G_a)\le\max\{q,8\}=q$ for $a=2,3$.
			Thus $\sum_{a\in E}f(G_a)\le6q-6\le10q-20$.
			Finally, if every $G_a$ has matching number at most one,
			\eqref{eq-common-star} and \eqref{eq-small-graphs}, with $k=q\ge12$, give
			$\sum_{a\in E}f(G_a)\le6q-3\le10q-20$.
			Part~\textup{(iii)} now follows from \eqref{eq-cap-reduce}.
		}
	\end{proof}
	
	\begin{lemma}\label{lem-r3-linear}
		Let $\cH$ be a $3$-graph on $n$ vertices with $\nu(\cH)=m\ge1$.
		If
		\[
		n\ge
		\begin{cases}
			13,&m=1,\\
			\max\{3m+12,\,7m/2+9\},&m\ge2,
		\end{cases}
		\]
		then $\sigma_3(\cH)\le3D_m(n)$.
	\end{lemma}
	\begin{proof}
		Suppose to the contrary that $\sigma_3(\cH)>3D_m(n)$.
		Fix a maximum matching $M=\{E_1,\ldots,E_m\}$ with $A=V(M)$. Let
		$Y=V(\cH)\setminus A$ with $t=\lvert Y\rvert =n-3m$.
		Since $M$ is maximum,  $\cH[Y]$ has no edge and the assumption on $n$ gives $t\ge10$.
	 By the definition of $\sigma_3(\cH)$, we have
		\[
		\begin{aligned}
		3\binom{t}{3}D_m(n)
		&<\binom{t}{3}\sigma_3(\cH)\le\sum_{S\in\binom{Y}{3}}\sum_{z\in S}d_{\cH}(z)=\sum_{z\in Y}d_{\cH}(z)
		\left\lvert\{S\in\binom{Y}{3}\mid z\in S\}\right\rvert=\binom{t-1}{2}\sum_{z\in Y}d_{\cH}(z).
		\end{aligned}
		\]
		Since $3\binom{t}{3}=t\binom{t-1}{2}$ and $t\ge10$, dividing by
		$\binom{t-1}{2}>0$ gives
		\begin{equation}\label{eq-sum-lower}
			\sum_{z\in Y}d_{\cH}(z)>tD_m(n).
		\end{equation}
		
		If $m=1$, then   every edge in $\cH$ intersects $E_1$. By
		{Lemma~\ref{lem-r3-cap}~(i)} with $E=E_1$ and $K=Y$, we have $
		\sum_{z\in Y}d_{\cH}(z)\le(t+1)t=tD_1(n)$,
		contrary to \eqref{eq-sum-lower}.
		
		Next assume that $m\ge2$. Then $t\ge\max\{12,m/2+9\}$.
		For $x\in E_i$ with $i\in[m]$, define the graph $
		G_{i,x}=\{R\in\tbinom Y2\mid x\cup R\in\cH\}$. For each $i\in[m]$, the family $
		\cH[E_i\cup Y]$
		is intersecting. If not, there exist two disjoint edges $
		F_1,F_2\in \cH[E_i\cup Y]$.
		Then $
		\{F_1\}\cup\{F_2\}\cup\{E_j\mid j\in[m]\setminus i\}$
		is a matching in $\cH$ of size $m+1$,
		contradicting $\nu(\cH)=m$.
		Let $
		I=
		\left\{
		i\in[m]
		\mid
		\max_{x\in E_i}\nu(G_{i,x})\ge5
		\right\}$. We claim that for each $i\in I$, there is a unique vertex $c_i\in E_i$ with
		$\nu(G_{i,c_i})\ge5$.
		Otherwise, there exist two distinct vertices $x,y\in E_i$ such that $G_{i,x}\neq\emptyset$, $G_{i,y}\neq\emptyset$ and $G_{i,x}$ contains three pairwise disjoint edges $
		R_1,R_2,R_3$.
		Choose any edge $
		R\in G_{i,y}$.
		By the definition of $G_{i,x}$ and $G_{i,y}$, we have $
		x\cup R_j\in\cH$ with $j\in[3]$
		and $
		y\cup R\in\cH$.
		Since $\cH[E_i\cup Y]$ is intersecting and $x\ne y$, the edge
		$y\cup R$ must intersect each of $x\cup R_1,
		x\cup R_2,
		x\cup R_3$
		in a vertex of $Y$. Hence $
		R\cap R_j\ne\emptyset$ with $j\in[3]$.
		However, $R_1,R_2,R_3$ are pairwise disjoint, while $R$ has only two
		vertices, so $R$ cannot intersect all three of them. This proves the uniqueness of such $c_i$.
		Let $
		g=\lvert I\rvert$, $a=m-g$ and $C=\{c_i\mid i\in I,\;\nu(G_{i,c_i})\ge5\}$. We first show the following claim.
		\begin{claim}\label{clm-core}
			For every edge $F\in\cH$, if $F\cap Y\ne\emptyset$ and
			$F\cap A\subseteq\bigcup\limits_{i\in I}E_i$, then $F\cap C\ne\emptyset$.
		\end{claim}
		\noindent\textbf{Proof of Claim~\ref{clm-core}.}
		Suppose to the contrary that $F\cap C=\emptyset$.
		{Since $F\cap Y\ne\emptyset$ and $\cH[Y]$ has no edge, the assumption
			$F\cap A\subseteq\bigcup_{i\in I}E_i$ gives $
			1\le\lvert\{h\in I\mid F\cap E_h\ne\emptyset\}\rvert\le2$.}
		If $\{h\in I\mid F\cap E_h\ne\emptyset\}=i$,
		then $\lvert F\cap Y\rvert\le2$.
		At most two edges from the five pairwise disjoint edges of $G_{i,c_i}$
		intersect $F\cap Y$. Thus there is $R\in G_{i,c_i}$ with
		$R\cap(F\cap Y)=\emptyset$.
		The edges $F$ and $c_i\cup R$, together with
		$E_h$ for $h\in[m]\setminus i$, form an $(m+1)$-matching
		in $\cH$, contradicting $\nu(\cH)=m$.
		If $\{h\in I\mid F\cap E_h\ne\emptyset\}=\{i,j\}$ with $i\ne j$,
		then $\lvert F\cap Y\rvert=1$.
		Choose $R_i\in G_{i,c_i}$ with $R_i\cap(F\cap Y)=\emptyset$.
		Since $\lvert(F\cap Y)\cup R_i\rvert=3$ and
		$\nu(G_{j,c_j})\ge5$, there is $R_j\in G_{j,c_j}$ with
		$R_j\cap((F\cap Y)\cup R_i)=\emptyset$.
		The edges $F$, $c_i\cup R_i$ and $c_j\cup R_j$,
		together with $E_h$ for $h\in[m]\setminus\{i,j\}$,
		form an $(m+1)$-matching in $\cH$, again contradicting
		$\nu(\cH)=m$. {This completes the proof of Claim~\ref{clm-core}.
		}	
		By Claim~\ref{clm-core}, every $F\in\cH$ with
		$F\cap Y\ne\emptyset$ and $F\cap C=\emptyset$
		intersects some $E_i$ with $i\in[m]\setminus I$.
		Since $\lvert F\rvert=3$ and $F\cap Y\ne\emptyset$,
		such an edge intersects at most two members of $M$.
		We now estimate $
		\sum_{z\in Y}d_{\cH}(z)
		=\sum_{F\in\cH}\lvert F\cap Y\rvert$
		by classifying the edges $F\in\cH$ with
		$F\cap Y\ne\emptyset$ according to their
		intersections with $C$.
		
		\medskip
		\noindent\textbf{Case~1.}  $F\cap C\ne\emptyset$.
		
		For each $z\in Y$, $
		\lvert\{F\in\cH\mid z\in F,\;F\cap C\ne\emptyset\}\rvert
		\le B_3(n)-\binom{n-g-1}{2}=D_g(n)$.
		Consequently, the contribution of these edges to
		$\sum_{z\in Y}d_{\cH}(z)$ is
		\[
		\begin{aligned}
			\sum_{\substack{F\in\cH\\F\cap Y\ne\emptyset,\;F\cap C\ne\emptyset}}
			\lvert F\cap Y\rvert
			&=\sum_{z\in Y}
			\lvert\{F\in\cH\mid z\in F,\;F\cap C\ne\emptyset\}\rvert\le tD_g(n).
		\end{aligned}
		\]
		
		\medskip
		\noindent\textbf{Case~2.}  $F\cap C=\emptyset$.
		
		In this case, $F\cap E_i\neq\emptyset$
		for some $i\notin I$ by the preceding argument. Then every such edge $F$ satisfies  one of the following three subcases.
		
		\medskip
		\noindent\textbf{Case~2.1.} $F\subseteq E_i\cup Y$ with $i\in [m]\setminus I$.
		
		In this case, $\nu(G_{i,x})\le4$ for each $x\in E_i$ with $i\in [m]\setminus I$. Since $\cH[E_i\cup Y]$ is intersecting with $t\geq12$, {Lemma~\ref{lem-r3-cap}~(iii)}
		gives $
		\sum_{z\in Y}d_{\cH[E_i\cup Y]}(z)\le10t-20$.
		Each edge in this subcase intersects exactly one member of $M$.
		Thus the total contribution of these edges $F$ in this subcase to $\sum_{z\in Y}d_{\cH}(z)$ is
		\[
		\sum_{i\in[m]\setminus I}
		\sum_{\substack{F\in\cH[E_i\cup Y]\\F\cap Y\ne\emptyset}}
		\lvert F\cap Y\rvert
		=\sum_{i\in[m]\setminus I}\sum_{z\in Y}
		d_{\cH[E_i\cup Y]}(z)
		\le a(10t-20).
		\]
		
		\medskip
		\noindent\textbf{Case~2.2.} $\lvert F\cap E_i\rvert=\lvert F\cap E_j\rvert=1$ with distinct $ i,j\in[m]\setminus I$.
		
		Let $
		\mathcal N_{ij}=
		\left\{
		F\in\cH
		\mid
		\lvert F\cap Y\rvert
		=
		\lvert F\cap E_i\rvert
		=
		\lvert F\cap E_j\rvert
		=
		1
		\right\}$ with $i,j\in[m]\setminus I$.
		We claim that $
		\nu(\mathcal N_{ij})\le2$.
		Otherwise, there exist three pairwise disjoint members $
		F_1,F_2,F_3\in \mathcal N_{ij}$.
		Since $
		Y\cap A=\emptyset$ with $A=V(M)$,
		$F_x\cap E_h=\emptyset$ for every $x\in[3]$ and every $h\in[m]\setminus\{i,j\}$. Hence $
		\{F_1,F_2,F_3\}
		\cup
		\{E_h\mid h\in[m]\setminus\{i,j\}\}$
		is a matching in $\cH$ of size $
		3+(m-2)=m+1$,
		contradicting $\nu(\cH)=m$. Therefore $\nu(\mathcal N_{ij})\le2$. For distinct $i,j\in[m]$, let $
		\mathcal R_{ij}
		:=
		\left\{
		\{z,x,y\}
		\mid
		z\in Y,\ x\in E_i,\ y\in E_j
		\right\}$.
		Thus $\mathcal N_{ij}\subseteq\mathcal R_{ij}$ and $
		\lvert\mathcal R_{ij}\rvert
		=
		\lvert Y\rvert\lvert E_i\rvert\lvert E_j\rvert
		=
		9t$.
		Choose uniformly at random a matching $M'\subseteq\mathcal R_{ij}$ of size
		three. Since $
		\lvert E_i\rvert=\lvert E_j\rvert=3$,
		every such $M'$ uses all three vertices of $E_i$, all three vertices
		{
			of $E_j$ and three distinct vertices of $Y$.
			Since $M'$ is a matching and $\nu(\mathcal N_{ij})\le2$, we have $
			\lvert M'\cap\mathcal N_{ij}\rvert\le2$,
			and hence $
			\mathbb E (\lvert M'\cap\mathcal N_{ij}\rvert)\le2$.
			Note that every member of $\mathcal R_{ij}$ has the same probability of
			belonging to $M'$. Since $M'$ always contains exactly three members
			of $\mathcal R_{ij}$, for every $F\in\mathcal R_{ij}$ we have $\mathbb P(F\in M')
			=
			\frac{3}{\lvert\mathcal R_{ij}\rvert}
			=
			\frac{1}{3t}$.
			Therefore, by linearity of expectation,
			\[
			\mathbb E(\lvert M'\cap\mathcal N_{ij}\rvert)
			=
			\sum_{F\in \mathcal N_{ij}}\mathbb P(F\in M')
			=
			\frac{\lvert \mathcal N_{ij}\rvert}{3t}.
			\]
			Combining this with $\mathbb E (\lvert M'\cap\mathcal N_{ij}\rvert)\le2$, we obtain $
			\lvert \mathcal N_{ij}\rvert
			\le
			6t$.
		}
		
		Since $\lvert F\cap Y\rvert=1$ for every $F\in\mathcal N_{ij}$,
		the total contribution in this subcase to $\sum_{z\in Y}d_{\cH}(z)$ is
		\[
		\sum_{\{i,j\}\in\binom{[m]\setminus I}{2}}
		\sum_{F\in\mathcal N_{ij}}\lvert F\cap Y\rvert
		=\sum_{\{i,j\}\in\binom{[m]\setminus I}{2}}
		\lvert\mathcal N_{ij}\rvert
		\le6t\binom a2.
		\]
		
		\medskip
		\noindent\textbf{Case~2.3.} $\lvert F\cap E_i\rvert=\lvert F\cap E_j\rvert=1$ with $i\in [m]\setminus I$ and $j\in I$.
		
		Fix $i\notin I$ and let $
		B=\bigcup_{j\in I}\bigl(E_j\setminus c_j\bigr)$.
		Then $
		\lvert B\rvert=2g$.
		Let $\cF_i=\{F\in\cH\mid \lvert F\cap Y\rvert
		=
		\lvert F\cap E_i\rvert
		=
		\lvert F\cap B\rvert
		=1\}$.
		We claim that $\cF_i$ is intersecting.
		Otherwise, there exist $F_1,F_2\in\cF_i$ such that $F_1\cap F_2=\emptyset$.
		Choose $j,j'\in I$ such that $(F_1\cap B)\subseteq E_j$ and $(F_2\cap B)\subseteq E_{j'}$. If $j=j'$, then $
		\lvert(F_1\cup F_2)\cap Y\rvert=2$.
		Since $j\in I$, the graph $G_{j,c_j}$ contains five pairwise
		disjoint edges.  Hence there is $
		R\in G_{j,c_j}$
		such that $
		R\cap\bigl((F_1\cup F_2)\cap Y\bigr)=\emptyset$. Moreover, $
		F_1,F_2, c_j\cup R$
		are three pairwise disjoint edges. Together with the edges $E_x$ for $x\in[m]\setminus\{i,j\}$, they form an $(m+1)$-matching, a
		contradiction. If $j\neq j'$, we can choose $
		R_j\in G_{j,c_j}$
		such that $
		R_j\cap\bigl((F_1\cup F_2)\cap Y\bigr)=\emptyset$. Since $G_{j',c_{j'}}$ contains five
		pairwise disjoint edges and at most four vertices belong to $((F_1\cup F_2)\cap Y)\cup R_j$, we can choose $
		R_{j'}\in G_{j',c_{j'}}$
		such that $
		R_{j'}\cap
		\left(
		((F_1\cup F_2)\cap Y)\cup R_j
		\right)
		=\emptyset$. Consequently $
		F_1,
		F_2,
		c_j\cup R_j,
		c_{j'}\cup R_{j'}$
		are pairwise disjoint. Together with the edges $E_x$ for $x\in[m]\setminus\{i,j,j'\}$, they form an $(m+1)$-matching.  This contradiction proves that $\cF_i$ is
		intersecting. Next we consider $\lvert\cF_i\rvert$. Recall that $\lvert Y\rvert=t$, $\lvert E_i\rvert=3$ and $
		\lvert B\rvert=2g$.
		For fixed $i\notin I$, let $
		\mathcal R_i
		:=
		\left\{
		\{z,x,y\}
		\mid
		z\in Y,\ x\in E_i,\ y\in B
		\right\}$.
		Then $
		\cF_i\subseteq\mathcal R_i$ and $
		\lvert\mathcal R_i\rvert
		=
		\lvert Y\rvert\lvert E_i\rvert\lvert B\rvert
		=
		t\cdot3\cdot2g
		=
		6gt$. Suppose first that $g\ge2$. Since $
		\lvert E_i\rvert=3$, $
		\lvert B\rvert=2g\ge4$ and $
		\lvert Y\rvert=t\ge3$, we can choose uniformly at random a matching
		$M'\subseteq\mathcal R_i$ of size three.
		{
			Since $\cF_i$ is intersecting and $M'$ is a matching, we have $
			\lvert M'\cap\cF_i\rvert\le1$,
			and hence $
			\mathbb E(\lvert M'\cap\cF_i\rvert)\le1$.
			Every member of $\mathcal R_i$ has the same probability of belonging to
			$M'$. Since $M'$ always contains exactly three members of
			$\mathcal R_i$, for every $F\in\mathcal R_i$ we have
			\[
			\mathbb P(F\in M')
			=
			\frac{3}{\lvert\mathcal R_i\rvert}
			=
			\frac{3}{6gt}
			=
			\frac1{2gt}.
			\]
			Therefore, by linearity of expectation,
			\[
			\mathbb E(\lvert M'\cap\cF_i\rvert)
			=
			\sum_{F\in\cF_i}\mathbb P(F\in M')
			=
			\frac{\lvert\cF_i\rvert}{2gt}.
			\]
			Combining this with $\mathbb E(\lvert M'\cap\cF_i\rvert)\le1$, we obtain $
			\lvert\cF_i\rvert\le2gt$.
		}
		
		Now suppose that $g=1$. Then $
		\lvert B\rvert=2$
		and $
		\lvert\mathcal R_i\rvert
		=
		t\cdot3\cdot2
		=
		6t$.
		{
			Choose uniformly at random a matching $M^*\subseteq\mathcal R_i$ of size two.
			Again, since $\cF_i$ is intersecting, $
			\lvert M^*\cap\cF_i\rvert\le1$,
			so $
			\mathbb E(\lvert M^*\cap\cF_i\rvert)\le1$.
			Every member of $\mathcal R_i$ has the same probability of belonging to
			$M^*$ and $M^*$ contains exactly two members of $\mathcal R_i$. Hence,
			for every $F\in\mathcal R_i$,
			\[
			\mathbb P(F\in M^*)
			=
			\frac{2}{\lvert\mathcal R_i\rvert}
			=
			\frac{2}{6t}
			=
			\frac1{3t}.
			\]
			By linearity of expectation,
			\[
			\mathbb E(\lvert M^*\cap\cF_i\rvert)
			=
			\sum_{F\in\cF_i}\mathbb P(F\in M^*)
			=
			\frac{\lvert\cF_i\rvert}{3t}.
			\]
			Thus $
			\lvert\cF_i\rvert\le3t$.
		}
		
		Finally, if $g=0$, then $B=\emptyset$, and hence $
		\cF_i=\emptyset$.
		In all cases $
		\lvert\cF_i\rvert\le(2g+\chi_{\{g=1\}})t$.
		Since $\lvert F\cap Y\rvert=1$ for every $F\in\cF_i$,
		the total contribution in this subcase to $\sum_{z\in Y}d_{\cH}(z)$ is
		\[
		\sum_{i\in[m]\setminus I}\sum_{F\in\cF_i}\lvert F\cap Y\rvert
		=\sum_{i\in[m]\setminus I}\lvert\cF_i\rvert
		\le a(2g+\chi_{\{g=1\}})t.
		\]
		
		Consequently,
		\begin{equation}\label{eq-sum-upper}
			\sum_{z\in Y}d_{\cH}(z)
			\le
			tD_g(n)
			+a(10t-20)
			+6t\binom a2
			+a(2g+\chi_{\{g=1\}})t.
		\end{equation}
		If $a=0$, then $g=m$ and \eqref{eq-sum-upper} gives
		$\sum_{z\in Y}d_{\cH}(z)\le tD_m(n)$,
		contradicting \eqref{eq-sum-lower}. Hence we may assume that $a\ge1$.
		Since $g=m-a$ and $n=3m+t$,
		\begin{align*}
			&D_m(n)-D_g(n)-a\left(10-\frac{20}{t}+3(a-1)+2g+\chi_{\{g=1\}}\right)\\
			&\quad=\binom{n-g-1}{2}-\binom{n-m-1}{2}
			-a\left(10-\frac{20}{t}+3(a-1)+2g+\chi_{\{g=1\}}\right)\\
			&\quad=a(n-m-1)+\binom a2
			-a\left(10-\frac{20}{t}+3(a-1)+2g+\chi_{\{g=1\}}\right)=a\left(t-\frac a2-\frac{17}{2}-\chi_{\{g=1\}}+\frac{20}{t}\right).
		\end{align*}
		Since $a=m-g$ and $g$ is a nonnegative integer, $
		\frac a2+\chi_{\{g=1\}}
		=\frac{m-g}{2}+\chi_{\{g=1\}}\le\frac{m+1}{2}$.
		Combining this with $t\ge\max\left\{12,\frac m2+9\right\}$, we have $
		t-\frac a2-\frac{17}{2}-\chi_{\{g=1\}}+\frac{20}{t}
		\ge
		t-\frac m2-9+\frac{20}{t}>0$. Hence
		\[
		D_m(n)-D_g(n)
		>
		a\left(
		10-\frac{20}{t}+3(a-1)+2g+\chi_{\{g=1\}}
		\right).
		\]
		Multiplying by $t$ and using \eqref{eq-sum-upper}, we obtain $
		\sum_{z\in Y}d_{\cH}(z)<tD_m(n)$,
		contradicting \eqref{eq-sum-lower}. This completes the proof.
	\end{proof}
	
	\begin{proof}[\textbf{Proof of Theorem~\ref{thm-r3}}]
		Suppose to the contrary that $\nu(\cH)=m<s$.
		If $m=0$, then $\cH$ has no edge with
		$\sigma_3(\cH)=0\le3D_{s-1}(n)$, contrary to
		$\sigma_3(\cH)>3D_{s-1}(n)$. Thus $m\ge1$.
		If $s=2$, {Lemma~\ref{lem-ore-ekr}} with $r=3$ and $n\ge7$
		gives a contradiction.
		Next assume $s\ge3$. The assumed bound on $n$ can be written as $
		n\ge\max\left\{3s+9,\left\lceil\frac{7s+11}{2}\right\rceil\right\}$.
		Since $m=\nu(\cH)\le s-1$, we have $
		3m+12\le3s+9$ and $
		\frac{7m}{2}+9\le\frac{7s+11}{2}$.
		The assumed lower bound on $n$ therefore gives $
		n\ge\max\left\{3m+12,\frac{7m}{2}+9\right\}$.
		By Lemma~\ref{lem-r3-linear} we have $
		\sigma_3(\cH)\le3D_m(n)$.
		Since $D_k(n)$ is increasing in $k$ and $m\le s-1$, $
		D_m(n)\le D_{s-1}(n)$.
		Consequently, $
		\sigma_3(\cH)\le3D_m(n)\le3D_{s-1}(n)$,
		contradicting the assumption
		$\sigma_3(\cH)>3D_{s-1}(n)$. This completes the proof.
	\end{proof}
	\begin{remark}
		Theorem~\ref{thm-r3} implies Conjecture~\ref{conj-bpr} for $r=3$.
		Taking $r=3$ in Theorem~\ref{thm-main} gives the stronger bound
		$n\ge3s+1$ for every $s\ge s_0(3)$.
		Whether the same bound suffices for all remaining values
		$3\le s<s_0(3)$ remains open.
	\end{remark}
	\section[Concluding remarks]{Concluding remarks}\label{sec-additive}
	
	For each integer $r\ge3$, let $C_*(r)$ be the minimum nonnegative
	integer $C$ such that, for every $s\ge1$ and $n\ge rs+C$,
	an $n$-vertex $r$-graph $\cH$ with
	$\sigma_r(\cH)>rD_{s-1}(n)$ has $\nu(\cH)\ge s$.
	Set $C_*(r)=\infty$ if no such integer exists.
	For $r\ge3$ and $s\ge1$, let $n_{\mathrm{c}}(r,s)$ be the
	minimum integer $n_0\ge rs$ such that this implication holds for
	every $n\ge n_0$. Such an integer exists by
	Theorem~\ref{thm-bpr}~(i) if $s\ge2$; while the case $s=1$ is trivial.

\begin{lemma}\label{lem-obstruct}
Let  $r\geq3$,
$s\ge2$ and $1\le u\le c_r=\lfloor\alpha_r\rfloor$ are integers and $\cH_{r,s,u}$ be the disjoint union of the complete
$r$-graph on $rs-1$ vertices and $u$ isolated vertices, where $
\alpha_r=\frac{(r-1)^{r-1}}{r^{r-2}}$. Then, for all sufficiently large $s$, $\cH_{r,s,u}$ has order $n=rs-1+u$ and
$\nu(\cH_{r,s,u})=s-1$ with $
 \sigma_r(\cH_{r,s,u})=(r-u)\binom{rs-2}{r-1}>rD_{s-1}(n)$.
\end{lemma}
\begin{proof} From the construction of $\cH_{r,s,u}$,
 we have $n=rs-1+u$. Let $I$ be the set of $u$ isolated vertices.
By \eqref{eq-complete-match}, we have $
 \nu(\cH_{r,s,u})=\nu(\cH_{r,s,u}-I)
 =\left\lfloor\frac{rs-1}{r}\right\rfloor=s-1$.
Every vertex in $I$ has degree zero and every vertex outside
$I$ has degree $\binom{rs-2}{r-1}$. An $r$-set $S$ in $\cH_{r,s,u}$ is a
nonedge if and only if $S\cap I\ne\emptyset$. For every
$r$-set $S$ with $S\cap I\ne\emptyset$, $
 \sum_{v\in S}d_{\cH_{r,s,u}}(v)
 =(r-\lvert S\cap I\rvert)\binom{rs-2}{r-1}
 \ge(r-u)\binom{rs-2}{r-1}$.
Since $1\le u\le c_r=\lfloor\alpha_r\rfloor\le r-1$, we can choose an $r$-set $S$ with $I\subseteq S$.
Consequently, $
 \sigma_r(\cH_{r,s,u})=(r-u)\binom{rs-2}{r-1}$.

Note that $u\le c_r$. For fixed $r$ and
$1\le u\le r-1$, as $s\to\infty$, we have $
 \binom{rs-2}{r-1}=\frac{r^{r-1}}{(r-1)!}s^{r-1}+O_r(s^{r-2})$, $
 \binom{rs+u-2}{r-1}
 =\frac{r^{r-1}}{(r-1)!}s^{r-1}+O_r(s^{r-2})$ and $
 \binom{(r-1)s+u-1}{r-1}
 =\frac{(r-1)^{r-1}}{(r-1)!}s^{r-1}+O_r(s^{r-2})$.
Since $
 D_{s-1}(rs-1+u)
 =\binom{rs+u-2}{r-1}-\binom{(r-1)s+u-1}{r-1}$,
we obtain
\begin{align}
 \sigma_r(\cH_{r,s,u})-rD_{s-1}(rs-1+u)&=
 \frac{(r-u)r^{r-1}-r^r+r(r-1)^{r-1}}{(r-1)!}s^{r-1}
 +O_r(s^{r-2})\nonumber\\
 &=
 \frac{r^{r-1}}{(r-1)!}(\alpha_r-u)s^{r-1}
 +O_r(s^{r-2}).\label{e49}
\end{align}
Since $\alpha_r$ is not an integer,
$\alpha_r-u\ge\alpha_r-c_r>0$. Dividing \eqref{e49}
by $s^{r-1}$, the error term tends to zero and the remaining
term is positive. Thus
$\sigma_r(\cH_{r,s,u})>rD_{s-1}(n)$ for sufficiently large
$s$.
\end{proof}

 Taking $u=c_r$ in the above construction gives, for every
sufficiently large $s$, an $r$-graph on $rs+c_r-1$ vertices
satisfying $\sigma_r(\cH)>rD_{s-1}(rs+c_r-1)$ but having
matching number $s-1$. Thus the order bound in
Theorem~\ref{thm-main} cannot be decreased by one, and
$C_*(r)\ge c_r$. The case $u=2$ is the construction in
\cite{BPR2026}.

 For $s=2$, Lemma~\ref{lem-ore-ekr} gives the sufficient
order $n\ge2r+1$ for Conjecture \ref{conj-bpr}. The following result shows that this bound
is best possible.
\begin{corollary}\label{cor-s2}
 For every integer $r\ge3$, $n_{\mathrm{c}}(r,2)=2r+1$.
\end{corollary}
\begin{proof}
 By Lemma~\ref{lem-ore-ekr}, every $r$-graph $\cH$ on
$n\ge2r+1$ vertices with $\sigma_r(\cH)>rD_1(n)$ has a
matching of size two. Hence $n_{\mathrm{c}}(r,2)\le2r+1$.
Let $V$ be a vertex set of size
$2r$. Choose $A\subseteq V$ with $\lvert A\rvert=3$ and let $
 \cH=\left\{F\in\binom Vr\mid\lvert F\cap A\rvert\in\{0,2\}\right\}$.
If $F,F'\in\cH$ are disjoint, then $F\cup F'=V$ and $
 \lvert F\cap A\rvert+\lvert F'\cap A\rvert=\lvert A\rvert=3$.
This is impossible because both $\lvert F\cap A\rvert$ and $\lvert F'\cap A\rvert$ are even.
Hence $\nu(\cH)\le1$. Choosing two vertices of $A$ and
$r-2$ vertices of $V\setminus A$ gives an edge of $\cH$,
so $\nu(\cH)=1$.
Write $D=D_1(2r)=\binom{2r-2}{r-2}$.
For $v\in A$, an edge $F$ containing $v$ satisfies $|F\cap A|=2$. Thus, for $v\in A$, we have
\[
\begin{aligned}
 d_\cH(v)-D
 &=2\binom{2r-3}{r-2}-\binom{2r-2}{r-2}=\left(2-\frac{2r-2}{r}\right)\binom{2r-3}{r-2}
 =\frac2r\binom{2r-3}{r-2}>0.
\end{aligned}
\]
For $v\in V\setminus A$, an edge containing $v$ contains
either no vertex of $A$ or two vertices of $A$. Therefore $
 d_\cH(v)
 =\binom{2r-4}{r-1}+3\binom{2r-4}{r-3}
 =4\binom{2r-4}{r-3}$.
Consequently,
\[
\begin{aligned}
 d_\cH(v)-D
 &=\left(4-\frac{(2r-2)(2r-3)}{r(r-2)}\right)
       \binom{2r-4}{r-3}=\frac{2(r-3)}{r(r-2)}\binom{2r-4}{r-3}\ge0.
\end{aligned}
\]
Every nonedge $S$ satisfies $\lvert S\cap A\rvert\in\{1,3\}$.
The degree estimates above show that
\[
 \sum_{v\in S}d_\cH(v)-rD
 =\sum_{v\in S}(d_\cH(v)-D)
 \ge\sum_{v\in S\cap A}(d_\cH(v)-D)>0.
\]
Thus $\sigma_r(\cH)>rD_1(2r)$ with $\nu(\cH)=1$.
By the definition of $n_{\mathrm{c}}(r,2)$,
$n_{\mathrm{c}}(r,2)\ge2r+1$. Combining the two bounds
proves the result.
\end{proof}

 Theorem~\ref{thm-main} and Corollary~\ref{cor-s2} yield
a sufficient order for Conjecture \ref{conj-bpr} of the form $rs+C(r)$ for every $s\ge1$.
More precisely, we have the following result.
\begin{corollary}\label{cor-all-s}
For every integer $r\ge3$, there exists an integer $C(r)$
depending only on $r$ such that, for all integers $s\ge1$ and
$n\ge rs+C(r)$, every $n$-vertex $r$-graph $\cH$ satisfying
$\sigma_r(\cH)>rD_{s-1}(n)$ has $\nu(\cH)\ge s$.
\end{corollary}
\begin{proof}
Fix $r\ge3$, let $s_0(r)\ge3$ be as in Theorem~\ref{thm-main},
and let $
 C(r)=c_r+(r-1)(s_0(r)-3)$.
 Let $s\ge1$ and $n\ge rs+C(r)$, and let $\cH$ satisfy the stated Ore-degree condition.
If $\cH$ is complete, then
$\nu(\cH)=\lfloor n/r\rfloor\ge s$ by~\eqref{eq-complete-match}.
Thus assume that $\cH$ has a nonedge.
Since $s_0(r)\ge3$ and $c_r\ge1$, we have $C(r)\ge c_r\ge1$.
For $s=1$, the condition $\sigma_r(\cH)>rD_0(n)=0$
excludes the empty hypergraph, so $\nu(\cH)\ge1$.
For $s=2$, we have $n\ge2r+1$, and $\nu(\cH)\ge2$ by Corollary~\ref{cor-s2}.
For $s\ge s_0(r)$ with $n\ge rs+c_r$, the conclusion follows from
Theorem~\ref{thm-main}.

It remains to consider $3\le s<s_0(r)$. Let $q=s_0(r)-s$
and take a set $Z$ of $q$ new vertices disjoint from $V(\cH)$.
Define an $r$-graph  $
 \cH^+=\cH\cup
 \left\{F\in\binom{V(\cH)\cup Z}{r}\mid F\cap Z\ne\emptyset\right\}$ of order $n'=n+q$. Then $
 n'-rs_0(r)
 =n-rs-(r-1)q\ge c_r+(r-1)(s_0(r)-3)-(r-1)(s_0(r)-s)=c_r+(r-1)(s-3)\ge c_r$.
The hypergraphs $\cH^+$ and $\cH$ have the same nonedges.
For every $v\in V(\cH)$, the added edges containing $v$ in $\cH^+$ are
precisely the $r$-sets containing $v$ and intersecting $Z$.
Consequently, $
 d_{\cH^+}(v)=d_\cH(v)+B_r(n')-B_r(n)$. This implies $
 \sigma_r(\cH^+)=\sigma_r(\cH)+r\bigl(B_r(n')-B_r(n)\bigr)$.
We have $
 D_{s_0(r)-1}(n')
 =B_r(n')-\binom{n'-s_0(r)}{r-1}=D_{s-1}(n)+B_r(n')-B_r(n)$ from $n'-s_0(r)=n-s$,.
It follows that $
 \sigma_r(\cH^+)-rD_{s_0(r)-1}(n')
 =\sigma_r(\cH)-rD_{s-1}(n)>0$.
As $n'\ge rs_0(r)+c_r$, $\cH^+$ has a
matching $M^+$ of size $s_0(r)$ by Theorem~\ref{thm-main}.
Note that $
 \left\lvert\{F\in M^+\mid F\cap Z\ne\emptyset\}\right\rvert
 \le\sum_{F\in M^+}\lvert F\cap Z\rvert
 =\lvert V(M^+)\cap Z\rvert\le q$.
Thus $M=\{F\in M^+\mid F\cap Z=\emptyset\}$ is a matching
in $\cH$ with $\lvert M\rvert\ge s_0(r)-q=s$. This completes the proof.
\end{proof}

Let $r\ge3$, $s\ge1$ be integers and $\cH$ be an
$n$-vertex $r$-graph satisfying $n\ge rs+c_r$ and $\sigma_r(\cH)>rD_{s-1}(n)$.
For each fixed $r\ge3$, if $s\ge s_0(r)$, then $\nu(\cH)\geq s$ by Theorem~\ref{thm-main}; if $s=2$, then $\nu(\cH)\geq s$ by Corollary~\ref{cor-s2}; the case $s=1$ is trivial.
As in the proof of Theorem~\ref{thm-main}, adding new vertices
and all $r$-sets intersecting them as edges reduces the remaining cases
to finitely many pairs $(s,n)$.
Thus, to prove $C_*(r)=c_r$, it suffices to show that every
$r$-graph $\cH$ on $n$ vertices satisfying
$\sigma_r(\cH)>rD_{s-1}(n)$ has a matching of size $s$ for $
 3\le s<s_0(r)$ and $c_r\le n-rs\le c_r+r-2$.
This leads to the following question.
\begin{question}\label{q-exact-const}
For every integer $r\ge3$, is $
 C_*(r)=c_r=\left\lfloor\frac{(r-1)^{r-1}}{r^{r-2}}\right\rfloor$?
\end{question}

\section*{AI Declaration}
During the preparation of this work, the authors
used ChatGPT (OpenAI) to assist in exploring potential approaches to the proof of Theorem \ref{thm-r3}.
All AI-generated suggestions were verified and refined by the authors, who take full responsibility
for the correctness and originality of the paper.

\end{document}